\newcommand*{\ARXIV}{}%

\documentclass[10.5pt]{article}
\usepackage[displaymath, mathlines]{lineno}

\ifdefined\ARXIV
  \pdfoutput=1
\else
  \linenumbers
\fi

\usepackage{benstyle}

\newcommand*\patchAmsMathEnvironmentForLineno[1]{%
  \expandafter\let\csname old#1\expandafter\endcsname\csname #1\endcsname
  \expandafter\let\csname oldend#1\expandafter\endcsname\csname end#1\endcsname
  \renewenvironment{#1}%
     {\linenomath\csname old#1\endcsname}%
     {\csname oldend#1\endcsname\endlinenomath}}%
\newcommand*\patchBothAmsMathEnvironmentsForLineno[1]{%
  \patchAmsMathEnvironmentForLineno{#1}%
  \patchAmsMathEnvironmentForLineno{#1*}}%
\AtBeginDocument{%
\patchBothAmsMathEnvironmentsForLineno{equation}%
\patchBothAmsMathEnvironmentsForLineno{align}%
\patchBothAmsMathEnvironmentsForLineno{flalign}%
\patchBothAmsMathEnvironmentsForLineno{alignat}%
\patchBothAmsMathEnvironmentsForLineno{gather}%
\patchBothAmsMathEnvironmentsForLineno{multline}%
}

\usepackage[usenames,dvipsnames,svgnames,table]{xcolor}
\usepackage{mathtools,ifthen}
\usepackage{array,multirow,booktabs} 
\usepackage{lipsum}
\usepackage{bbm}

\newcommand{\anrev}[1]{{\leavevmode\color{BrickRed}{#1}}}
\renewcommand{\anrev}[1]{#1}

\newcommand{\mathd}{\mathrm{d}}
\newcommand{\dx}[1]{\mathd #1}

\newcommand{\edits}[1]{{\leavevmode\color{BrickRed}{#1}}}
\renewcommand{\edits}[1]{#1}

\newlength{\leftstackrelawd}
\newlength{\leftstackrelbwd}
\def\leftstackrel#1#2{\settowidth{\leftstackrelawd}%
{${{}^{#1}}$}\settowidth{\leftstackrelbwd}{$#2$}%
\addtolength{\leftstackrelawd}{-\leftstackrelbwd}%
\leavevmode\ifthenelse{\lengthtest{\leftstackrelawd>0pt}}%
{\kern-.5\leftstackrelawd}{}\mathrel{\mathop{#2}\limits^{#1}}}

\let\oldbibliography\thebibliography
\renewcommand{\thebibliography}[1]{%
  \oldbibliography{#1}%
  \setlength{\itemsep}{0pt}%
}

\begin{document}

\title{Compressed sensing with sparse corruptions: Fault-tolerant sparse collocation approximations}
\author{
  Ben Adcock\footnote{B. Adcock and A. Bao acknowledge the support of the Alfred P. Sloan Foundation and the Natural Sciences and Engineering Research Council of Canada through grant 611675.} \\ Department of Mathematics \\ Simon Fraser University \\ Burnaby, BC, Canada 
  \and 
  Anyi Bao\footnotemark[1] \\ Department of Mathematics \\ Simon Fraser University \\ Burnaby, BC, Canada \\[20pt]
  \and 
  John D.\ Jakeman\footnote{J.D.Jakeman's work was supported by DARPA EQUiPS.} \\ Computer Science Research Institute \\ Sandia National Laboratories \\ Albuquerque, NM, USA\\[10pt]
  \and 
  Akil Narayan\footnote{A. Narayan is partially supported by NSF DMS-1720416, AFOSR FA9550-15-1-0467, and DARPA EQUiPS N660011524053} \\ Department of Mathematics and \\ Scientific Computing and Imaging (SCI) Institute \\ University of Utah \\ Salt Lake City, UT, USA
}

\maketitle
\begin{abstract}
  The recovery of approximately sparse or compressible coefficients in a polynomial chaos expansion is a common goal in many modern parametric uncertainty quantification (UQ) problems. However, relatively little effort in UQ has been directed toward theoretical and computational strategies for addressing the sparse \textit{corruptions} problem, where a small number of measurements are highly corrupted. Such a situation has become pertinent today since modern computational frameworks are sufficiently complex with many interdependent components that may introduce hardware and software failures, some of which can be difficult to detect and result in a highly polluted simulation result. 
  
  In this paper we present a novel compressive sampling-based theoretical analysis for a regularized $\ell^1$ minimization algorithm that aims to recover sparse expansion coefficients in the presence of measurement corruptions. Our recovery results are uniform (the theoretical guarantees hold for all compressible signals and compressible corruptions vectors), and prescribe algorithmic regularization parameters in terms of a user-defined \textit{a priori} estimate on the ratio of measurements that are believed to be corrupted. We also propose an iteratively reweighted optimization algorithm that automatically refines the value of the regularization parameter, and empirically produces superior results. Our numerical results test our framework on several medium-to-high dimensional examples of solutions to parameterized differential equations, and demonstrate the effectiveness of our approach.
\end{abstract}


\section{Introduction}
The approximation of function values using point evaluations or samples is necessary in a wide number of applications. Much attention has been focused recently on the approximation technique of compressive sampling (CS): The ability to recover sparse linear representations of a function from a given dictionary. This is a particularly important problem in parametric uncertainty quantification (UQ) where the number of parameters translates into the number of variables on which an unknown function depends (the ``dimension" of the problem). It is common for dimension to be very large, and the number of degrees of freedom in classical approximation strategies generally grows exponentially with the dimension. This makes classical computational procedures for approximating functions infeasible for large dimensions. 

In contrast, compressive sampling seeks a sparse representation of a function using only a small number of samples or measurements, regardless of the parametric dimension. In a non-intrusive UQ pipeline, each function sample corresponds to a potentially large-scale simulation, and so minimizing the requisite number of samples is desirable. When functions are sparse or compressible in a given basis or dictionary, this reconstruction procedure has the potential to mitigate the exponentially debilitating curse of dimensionality. Algorithms in UQ that utilize compressive sampling have enjoyed great success in recent years \cite{KarniadakisUQCS,yang_alternating_2011,yin_minimization_2015,PengHamptonDoostantweighted,jakeman_generalized_2016,DoostanOwhadiSparse,JakemanEtAl_l1Enhance,GuoEtAlRandomizedQuad,NarayanZhouCCP}. For related theoretical contributions, see \cite{AdcockCSFunInterp,Adcockl1Pointwise,ChkifaDownwardsCS,HamptonDoostanCSPCE,RauhutWardWeighted,Rauhut,YanGuoXui_l1UQ}.

Missing from the sparse recovery UQ contributions above is a concrete strategy for fault-tolerant or resilient algorithms. Ensuring modeling resilience for UQ in the presence of system failures is essential for credible prediction on new and emerging massively parallel systems. Fault-tolerant algorithms in general have become necessary in computational science since node failures on distributed architectures can yield corrupted data (the frequency of which increases as the number of processors increases), or algorithmic run-time software failures can result in polluted simulation results. These failures can generate polluted measurements in unpredictable and sometimes undetectable ways \cite{bridges_fault-tolerant_2012}. 

Faults can occur  due  to complex  combination  of  internal  and  external  conditions  that are difficult to reproduce. For example, bits may suffer random corruption, or physical defects in hardware may cause data faults. Corruption errors during model simulation can be grouped into two main types, soft and hard. In this paper, we consider hard faults as errors that cause the simulation to terminate prematurely and/or return obvious, automatically detectable error values such as NaN or Inf.  Hard faults by this definition are easy to identify and mark for discard, thus obviating or ameliorating the need for fault-tolerant algorithms. 

In contrast, soft failures are essentially random systematic corruption of results that are not easily identifiable. These soft failures pose challenges in UQ: A soft failure  will not cause obvious failure in fault-intolerant UQ methods; however, incorrect model values caused by soft failures can significantly degrade an approximation. It is in this case that we require the development of robust and resilient algorithms that can, ideally, deliver constant levels of performance when faced with a few highly corrupted data points.

To address this issue, fault-tolerant algorithms for UQ have been investigated in the context of multilevel Monte Carlo algorithms \cite{pauli_fault_thesis_2014,pauli_intrinsic_2015,pauli_fault_2014}, and in overdetermined least-squares polynomial recovery problems \cite{shin_correcting_2016}. To the best of our knowledge, there is no comprehensive research in the UQ literature on fault-tolerant sparse recovery algorithms, and in the compressive sampling literature only a handful of papers \cite{LaskaEtAlCorrupt,LiCorruptionsConstrApprox,NguyenTranCorrupt,stankovic_missing_2014,StuderEtAlCorrupt,DSuCSCorruptFourier,DSuCorrupted,WrightMaCorruption} deal with the problem of corrupted measurements.

The operative distinction in the problem we consider in this paper is a hardware or software fault resulting in occasional large-magnitude errors; we call this the problem of \textit{corruptions}. Existing CS algorithms are known to be stable with respect to small noise perturbations, but cannot handle sparse corruptions, i.e., situations when a small number of samples are highly corrupted with the corruption magnitudes much larger than typical noise. In this paper we present novel theory and application studies of a sparse corruptions algorithm for CS. The algorithm we use was considered in \cite{LiCorruptionsConstrApprox}, but we present more general theoretical guarantees on recovery, including practical guidance for the choice of algorithmic regularization parameters. 

For fault-tolerance in the context of the sparse recovery problem, the recovery properties of an ideal resilient algorithm would be agnostic to large-magnitude corruptions in a small number of function samples. As described above, these corruptions can arise due to unknown failure modes in computational models or because of large but intermittent measurement errors. Development of mathematical theory for the corrupted compressive sampling problem, and investigation of a corresponding resilient algorithm for sparse recovery of expansion coefficients are the central goals of this paper. The target applications we investigate are exemplars of a common task in UQ: recovery of approximately sparse expansion coefficients in an orthogonal polynomial (polynomial chaos) basis.

The theory and algorithms developed in this paper have the following features: 
\begin{itemize}
  \item The compressive sampling recovery theorems are uniform with respect to the function and the corruptions. That is, the recovery guarantees hold over all compressible functions having sparsely corrupted measurements for a single random sampling of measurements.
  \item The algorithm involves a tunable regularization parameter $\lambda$, and a theoretically optimal choice of this parameter is explicitly determined by our analytical results. This theoretically optimal value is defined only by the \textit{ratio} of measurement corruptions to signal sparsity. Since signal sparsity is frequently comparable to the number of measurements, this optimal $\lambda$ loosely translates into the fraction of measurement samples that are corrupted. From a user's point of view, our analysis thus suggests a value of $\lambda$ having knowledge only of the ratio of measurements believed to be corrupted.
  \item In experiments, we observe that optimal values of the regularization parameter are non-trivially dependent on the number of measurements, the signal sparsity, and the number of corruptions. We thus propose an iteratively reweighted algorithm for recovery that learns values of the regularization parameter. Our experiments suggest that these learned algorithmic parameters perform better than the value defined by our theoretical results, and thus this reweighted algorithm is more useful in practice. 
  \item The location and magnitude of the corruptions amongst the collection of function samples can be unknown, but the algorithm recovers those locations and the corresponding corruption values.
  \item The algorithm is robust to small, but non-sparse measurement errors -- e.g.\ due to noise, truncation of an infinite polynomial expansion or numerical error in computing function samples -- and moreover is \textit{noise-blind}.  That is to say, it requires no \textit{a priori} upper bound on such errors.
  \item \edits{The optimization problem we solve to compute solutions is from \cite{LiCorruptionsConstrApprox}, but our work is both a theoretical and practical advancement over the results in that reference. In order to show the solution computed is indeed the original sparse solution, \cite{LiCorruptionsConstrApprox} uses conditions on the restricted isometry constant (RIC) of the measurement matrix. Our results are a significant relaxation of previously reported conditions on the RIC (compare conditions on $\delta_{2s, 2k}$ in Lemma 2.3 of \cite{LiCorruptionsConstrApprox} versus our Theorem \ref{t:RIP_stable_robust}, equation \eqref{delta_cond}, and the discussion in Section \ref{ss:lambda-strategy}). The results for general sensing matrices in \cite{LiCorruptionsConstrApprox} are nonuniform with respect to the signal and corruptions support, and require certain models for the signal and corruptions; our results are uniform and require no model for the signal or corruptions, other than compressibility. Finally, our paper is also devoted to numerical investigation of the performance of the method, including practical guidance for choosing the regularization parameter $\lambda$; such thorough investigations are absent in \cite{LiCorruptionsConstrApprox}.}
\end{itemize}
We first introduce notation and summarize the main mathematical statements of this paper in Section \ref{sec:notation}. This is followed in Section \ref{sec:theory} by our theoretical analysis. Section \ref{sec:results} presents numerical results to complement our theoretical analysis and verify the practical efficacy of the algorithm.

\section{Model problem and main results}\label{sec:notation}

Let $f: \bbR^d \rightarrow \bbR$ denote an unknown function, and let $\left\{ \phi_j\right\}_{j=1}^N$ be a given dictionary of functions, $\phi_j : \bbR^d \rightarrow \bbR$. For example, the functions $\phi_j$ are frequently multivariate polynomial chaos basis elements; our capstone numerical examples will show results from such a basis. In scenarios of interest, the size $N$ of the dictionary is very large.

The ultimate goal is to recover coefficients $x_j$ that determine the approximation
\begin{align}\label{eq:f-definition}
  f(\xi) = \sum_{j=1}^N x_j \phi_j(\xi) + n(\xi),
\end{align}
using samples of $f$, where $n(\xi)$ is an assumed small discrepancy term between the exact function and its $N$-term linear approximation in $\phi_j$\footnote{Our notation suggests that $n = n(\xi)$ depends explicitly and deterministically on $\xi$; however, our theory encompasses the case when $n$ is a stochastic variable or process, e.g., independent Gaussian random variable additive perturbations of the measurements.}. For the purposes of exposition we assume $| n(\xi)| \leq \epsilon$ for some known uniform noise bound $\epsilon$; we will show later that lack of \textit{a priori} knowledge for this bound only affects theoretical results in benign ways. As described above, we assume the vector $x = \left(x_1, \ldots, x_N \right)^T \in \bbR^N$ to be compressible. Sparsity or compressibility of a vector can be quantified via its best $s$-term approximation error,
\begin{align*}
  \sigma_{s}(x)_p = \inf_{\|\tilde{x}\|_0 \leq s} \left\| x - \tilde{x} \right\|_p,
\end{align*}
where $\|\cdot\|_p$ is the standard $\ell^p$ norm on vectors; for $p=0$, $\| x \|_0$ is the sparsity of $x$, i.e., the number of non-zero elements in the vector.

With $\left\{\xi_1, \ldots, \xi_m\right\} \subset \bbR^d$ a collection of samples of $\xi$, we have the corresponding corrupted function measurements,
\begin{align*}
  y_k &= f(\xi_k) + c_k = \sum_{j=1}^N x_j \phi_j\left(\xi_k\right) + n\left(\xi_k\right) + c_k, & k&=1, \ldots m,
\end{align*}
where the corruption vector $c = \left(c_1, \ldots c_m \right)^T \in \bbR^m$ is assumed to be $k$-sparse but can have large entries. 
To enforce an underdetermined system, we assume $m < N$. Defining the rectangular matrix $A$ with entries $(A)_{j,k} = \phi_k\left(\xi_j\right)$, then the unknown vectors $x$ and $c$ satisfy the underdetermined linear system
\begin{align}\label{eq:linear-system}
  y = A x + c + n \in \bbR^m.
\end{align}
In order to compute the solution $(x,c)$ having knowledge of only $A$ and $y$, we consider the following model problem (see also \cite{LiCorruptionsConstrApprox} and references therein):
\be{
\label{l1_lambda_recovery}
\min_{z \in \bbC^N, d \in \bbC^m} \| z \|_{1} + \lambda \| d \|_{1}\ \mbox{subject to $\| A z + d - y \|_{2} \leq \epsilon \sqrt{m}$}.
}
Let $(\hat{x},\hat{c})$ be a minimizer of this problem, where $\hat{x} \in \bbR^N$ and $\hat{c} \in \bbR^m$.  Our objective is to obtain conditions on $A$ (in particular, on the number of measurements $m$) and $\lambda$ such that the error
\bes{
\| \hat{x} - x \|_{2} + \| c - \hat{c} \|_{2}
}
can be bounded by the best approximation numbers $\sigma_{s}(x)_1$ and $\sigma_{k}(c)_1$, and the noise magnitude $\epsilon$. 

\subsection{Main results}

In all that follows, the statement $a \lesssim b$ means $a \leq C b$ for some universal constant $C$. Our first main result shows that stable and robust recovery of $x$ and $c$ is implied by a certain modification of the classical Restricted Isometry Property (RIP) which incorporates the sparse corruptions term (Definition \ref{d:RIPcorruptions}).  Specifically, Theorem \ref{t:RIP_stable_robust} establishes that if the matrix $A$ satisfies the RIP for the corruptions problem of order $(2s,2k)$ (see Definition \ref{d:RIPcorruptions}) with constant $\delta_{2s,2k}$ satisfying
\begin{align}\label{eq:delta_intro}
\delta_{2s,2k} &< \frac{1}{\sqrt{1+\left ( \frac{1}{2\sqrt{2}} + \sqrt{\eta} \right)^2}}, & \eta = \frac{s + \lambda^2 k}{\min\left\{ s, \lambda^2 k \right\}},
\end{align}
then the following error bounds hold:
\begin{subequations}\label{eq:recovery-summary}
\begin{align}\label{eq:recovery-summary-l1}
  \left\| x - \hat{x} \right\|_1 + \lambda \left\| c - \hat{c} \right\|_1 &\lesssim \sigma_s(x)_1 + \lambda \sigma_k(c)_1 + \epsilon \sqrt{s + \lambda^2 k},\\\label{eq:recovery-summary-l2}
  \left\| x - \hat{x} \right\|_2 + \left\| c - \hat{c} \right\|_2 &\lesssim \left(1 + \eta^{1/4}\right) \left( \frac{\sigma_s(x)_1}{\sqrt{s}} + \frac{\sigma_k(c)_1}{\sqrt{k}} + \epsilon \right).
\end{align}
\end{subequations}
Our second main result (Theorem \ref{t:BOS-RIP}) provides explicit conditions on $m$, $s$ and $k$ for \R{eq:delta_intro} to hold for matrices of so-called bounded orthonormal system \cite[Chpt.\ 12]{FoucartRauhutCSbook}.  Specifically, suppose that $\{ \phi_j \}^{N}_{j=1}$ is an $L^2_{\dx{\nu}}(D)$-orthonormal system, where $\nu$ is a probability measure and $D \subset \bbR^d$ its support.  Define
\begin{align*}
  K \coloneqq \max_{j=1, \ldots, N} \sup_{\xi \in D} \left|\phi_j(\xi)\right| < \infty,
\end{align*}
and let $A = \left \{ \phi_j(\xi_i) \right \}^{m,N}_{i,j=1}$ where $\xi_1,\ldots,\xi_m$ are drawn i.i.d.\ according to $\nu$.  If 
\begin{align}\label{eq:m-bound}
\begin{aligned}
  m & \gtrsim \delta^{-2} \cdot K^2 \cdot s \cdot \left( \log^3(2s) \cdot \log(2 N) + \log \epsilon^{-1} \right),
  \\
    m & \gtrsim \delta^{-2} \cdot K \cdot s \cdot k,
  \end{aligned}
\end{align}
then with probability at least $1-\epsilon$, the restricted isometry constant $\delta_{2s,2k}$ of the scaled matrix $\frac{1}{\sqrt{m}} A$ satisfies $\delta_{2s,2k} \leq \delta$.

One can see from these estimates that optimizing $\eta$ over values of $\lambda$ yields a minimum value of $\eta = 2$ when $\lambda^2 = s/k$. Assuming $s \sim m$, this provides a concrete determination of the parameter $\lambda$ for use in \eqref{l1_lambda_recovery} having knowledge only of the ratio of corrupted measurements.  We note in passing that we do not believe that the second condition in \R{eq:m-bound} is sharp in the dependence on the product $s \cdot k$.  Improvement of this to a condition of the form
\be{
\label{optimal-conjecture}
m \gtrsim \delta^{-2} \cdot K \cdot k,
}
is left as a topic for future work.  Note that such a condition is known for Gaussian random matrices.  Moreover, a nonuniform recovery result with the scaling \R{optimal-conjecture} for exactly sparse coefficients $x$ and corruptions $c$ having random sign patterns was given in \cite{LiCorruptionsConstrApprox}.  See Section \ref{sec:rip-2} for further discussion.

It is common in compressed sensing to assume some \textit{a priori} known noise bound $\epsilon$ based on the user's knowledge of measurement noise or truncation error.  Although there are some results that circumvent this assumption \cite{AdcockCSFunInterp,Adcockl1Pointwise}, they typically yield somewhat weaker recovery guarantees. However, in the context of the sparse corruptions theory presented above, such prior knowledge of $\epsilon$ is not necessary for stable recovery: The error introduced by an unknown noise $\epsilon$ can be passed into theoretical estimates as a penalty of size $\epsilon$. To see this, note that if we define $c' \coloneqq \frac{1}{\sqrt{m}}(c + n)$, then the system $y = A x + c +n$ can be written as $\frac{1}{\sqrt{m}} y = \frac{1}{\sqrt{m}} A x + c'$. Solving \eqref{l1_lambda_recovery} by setting $\epsilon = 0$ results in the $\epsilon=0$ version of the estimate \eqref{eq:recovery-summary-l2} with $c'$ replacing $c$. However, the normalized best $k$-term approximation error to $c'$ appearing in \eqref{eq:recovery-summary-l2} is stable with respect to noise perturbations:
\begin{align*}
  \frac{\sigma_{k}(c')_{1}}{\sqrt{k}} \leq \frac{1}{\sqrt{k m}} \left( \sigma_{k}(c)_1 + \| n \|_{1}\right) \leq \frac{\sigma_{k}(c)_1}{\sqrt{k m}} + \sqrt{\frac{m}{k}} \epsilon.
\end{align*}
Here $\epsilon \geq \| n \|_{\infty}$ is any bound for the perturbation $n$ in the uniform norm. Using \eqref{optimal-conjecture}, we see that $\sqrt{\frac{m}{k}} \epsilon \lesssim \epsilon$, which is on the same order as the estimate \eqref{eq:recovery-summary-l2} that uses \textit{a priori} knowledge of $\epsilon$. A similar argument holds for the bound \eqref{eq:recovery-summary-l1}. 


While our theoretical results are thus insensitive to ignorance about small noise levels, we caution that it is always a good idea to use such information in practical recovery algorithms if available, e.g.\ as the result of cross validation.  See, for example, \cite{DoostanOwhadiSparse,KarniadakisUQCS,JakemanEtAl_l1Enhance}.

\subsection{Remarks on numerical results}
We postpone presenting numerical results until the end of this paper in Section \ref{sec:results}. However, some remarks on our findings are pertinent here in the context of the previous section's theory. First, the optimal value of $\lambda^2 = s/k$ that is suggested by \eqref{eq:delta_intro} does not appear to be the computationally optimal value of $\lambda$. \edits{That this fixed value of $\lambda$ is not the best is not surprising since the bounds \eqref{eq:recovery-summary} are derived using some loose inequalities. However, such bounds can be useful in understanding qualitative trends.} Results from our experimentation do suggest that large values of $\lambda$ more reliably recover corruptions when $s/k$ is large (see Figures \ref{f:model1} and \ref{f:model2}). This general trend in numerical results is consistent with the behavior of $\eta$ in \eqref{eq:delta_intro} as a function of $\lambda$ when $s/k$ is large.

We address this discrepancy between the theory and empirical results by propose an iteratively reweighted $\ell^1$ optimization scheme (see \cite{CandesWakinBoydReweighted}) that learns and updates the value of $\lambda$. Our results show that this proposed algorithm performs much better in practice than algorithms that fix $\lambda$. \edits{However, we do not present any theory to support the observed superiority of reweighted $\ell^1$ optimization schemes for the corruptions problem.}

Many of our capstone numerical examples are from applications using polynomial chaos expansions, where the compressible function has an expansion in a multivariate orthogonal polynomial basis. To simplify the presentation of our results, we focus on such examples where the basis is a tensor-product Legendre polynomial or Chebyshev polynomial system. Much recent work has shown that randomly generating measurements using samples from standard distributions (e.g., the uniform distribution) can accurately and near-optimally recover orthogonal polynomial expansions from such basis sets \cite{Rauhut,JakemanEtAl_l1Enhance,YanGuoXui_l1UQ}. Recovery in more general polynomial spaces has been investigated \cite{HamptonDoostanCSPCE,jakeman_generalized_2016,GuoEtAlRandomizedQuad}, but these methods usually rely on sophisticated sampling strategies and optimal sampling schemes are still an active area of research.

\section{Theory for the sparse corruptions problem}\label{sec:theory}


\begin{table}
  \begin{center}
  \resizebox{\textwidth}{!}{
    \renewcommand{\tabcolsep}{0.4cm}
    \renewcommand{\arraystretch}{1.3}
    {\scriptsize
    \begin{tabular}{@{}cp{0.8\textwidth}@{}}
      \toprule
      $m$ & number of measurements \\
      $N$ & length of sparse vector \\
      $x$ & sparse vector in $\bbC^N$ \\
      $c$ & corruptions vector in $\bbC^m$\\
      $A$ & $m \times N$ measurement matrix \\
      $n$ & noise vector in $\bbC^m$ \\
      $\epsilon$ & noise bound \\
      $\lambda$ & non-negative weighting parameter for the corruptions vector \\
      $\hat{x}$, $\hat{c}$ & solutions of the optimization problem \\
      $S$ & subset of $\left\{1, \ldots, N \right\}$, indices corresponding to $x$ \\
      $T$ & subset of $\left\{1, \ldots, m \right\}$, indices corresponding to $c$ \\
      $s$ & sparsity of $x$ \\
      $k$ & sparsity of $c$ \\
      $\Sigma_{s}$ & set of $s$-sparse vectors in $\bbC^N$ \\
      $\Sigma_{k}$ & set of $k$-sparse vectors in $\bbC^m$ \\
      $\sigma_{s}(x)_1$ & best $s$-term approximation error, measured in the $\ell^1$ norm \\
      $\sigma_{k}(c)_1$ & best $k$-term approximation error, measured in the $\ell^1$ norm \\
    \bottomrule
    \end{tabular}
  }
    \renewcommand{\arraystretch}{1}
    \renewcommand{\tabcolsep}{12pt}
  }
  \end{center}
  \caption{Notation used throughout this article.}\label{tab:notation}
\end{table}

We recall and summarize our notation for the sparse corruptions problem in Table \ref{tab:notation}. Our previous discussion was framed for real-valued signals $x$ and measurements $y$, but we now generalize to the complex-valued setting. This adds generality with no additional mathematical difficulty.

We follow a familiar path for deriving conditions on $m$ such that $\ell^1$ optimization problems recover sparse solutions (see, for example, \cite{FoucartRauhutCSbook}). Section \ref{sec:robust-nsp} defines an appropriate robust Null Space Property (NSP) for the matrix $A$ in the sparse corruptions setting. Under this property, we show that the recovery estimates \eqref{eq:recovery-summary} hold. In order to construct matrices $A$ that satisfy the robust NSP, Section \ref{sec:rip} generalizes the concept of the Restricted Isometry Property (RIP) for matrices to the sparse corruptions setting. That section shows that matrices satisfying the RIP for the sparse corruptions problem also satisfy the robust NSP.  Sections \ref{sec:rip-2} and \ref{sec:bos} show that if the dictionary elements $\phi_j$ form a bounded orthonormal system, then under the condition \eqref{eq:m-bound}, the matrix $A$ satisfies the RIP with high probability.  Finally, using these various results, we discuss a theoretically-optimal choice for $\lambda$ in Section \ref{ss:lambda-strategy}.

\subsection{The Robust Null Space Property for the sparse corruptions problem}\label{sec:robust-nsp}




The following two definitions are generalizations of robust null space properties (cf. \cite[Definition 4.17]{FoucartRauhutCSbook} and \cite[Definition 4.21]{FoucartRauhutCSbook}, respectively), and prescribe classes of matrices whose kernels do not contain sparse vectors.

\defn{
Let $1 \leq s \leq N$, $1 \leq k \leq m$ and $\lambda > 0$.  A matrix $A \in \bbC^{m \times N}$ satisfies the $\ell^1$-robust null space property of order $(s,k)$ with weight $\lambda$ if there exist constants $0 < \rho <1$ and $\tau > 0$ such that
\bes{
\| x_{S} \|_{1} + \lambda \| c_{T} \|_{1} \leq \rho \left ( \| x_{S^c} \|_{1} + \lambda \| c_{T^c} \|_{1} \right ) + \tau \| A x + c \|_{2},\quad \forall x \in \bbC^N,\ c \in \bbC^m,
}
for all sets $S \subseteq \{1,\ldots,N\}$ and $T \subseteq \{1,\ldots,m\}$ with $| S | \leq s$ and $|T| \leq k$. Above, $S^c$ is the complement of $S$ in $\{1, \ldots, N\}$, and similarly for $T^c$.
}

\defn{
Let $1 \leq s \leq N$, $1 \leq k \leq m$ and $\lambda > 0$.  A matrix $A \in \bbC^{m \times N}$ satisfies the $\ell^2$-robust null space property of order $(s,k)$ with weight $\lambda$ if there exist constants $0 < \rho <1$ and $\tau > 0$ such that
\be{
\label{l2_rNSP_def}
\sqrt{\| x_{S} \|^2_{2} + \| c_{T} \|^2_{2}} \leq \frac{\rho}{\sqrt{s+\lambda^2 k}} \left ( \| x_{S^c} \|_{1} + \lambda \| c_{T^c} \|_{1} \right ) + \tau \| A x + c \|_{2},\quad \forall x \in \bbC^N,\ c \in \bbC^m,
}
for all sets $S \subseteq \{1,\ldots,N\}$ and $T \subseteq \{1,\ldots,m\}$ with $| S | \leq s$ and $|T| \leq k$.
}



These definitions yield the following two results:

\lem{
\label{l:21_RNSP}
If $A \in \bbC^{m \times N}$ satisfies the $\ell^2$-robust null space property of order $(s,k)$ with weight $\lambda >0$ and constants $0 < \rho <1$, $\tau >0$ then it satisfies the $\ell^1$-robust null space property of order $(s,k)$ with weight $\lambda >0$ and constants $\rho$, $ \tau\sqrt{s+\lambda^2 k}$.
}
\prf{
Observe that
\bes{
\| x_{S} \|_{1} + \lambda \| c_{T} \|_{1} \leq \sqrt{s} \| x_{S} \|_{2} + \lambda \sqrt{k} \| c_{T} \|_{2} \leq \sqrt{s+\lambda^2 k} \sqrt{\| x_{S} \|^2_{2} + \| c_{T} \|^2_{2} }.
}
We now use the definition of the $\ell^2$-robust null space property.
}

\thm{
\label{t:rNSP_stable_robust}
Let $1 \leq s \leq N$, $1 \leq k \leq m$ and $\lambda > 0$ and suppose that $A \in \bbC^{m \times N}$ satisfies the $\ell^2$-robust null space property of order $(s,k)$ with weight $\lambda$.  Let $x \in \bbC^N$, $c \in \bbC^m$, $y \in \bbC^m$ and $\epsilon > 0$ be such that $\| A x + c - y \|_{2} \leq \epsilon$, and suppose that  $(\hat{x},\hat{c})$ is a minimizer of
\bes{
\min_{z \in \bbC^N, d \in \bbC^m} \| z \|_{1} + \lambda \| d \|_{1}\ \mbox{subject to $\| A z + d - y \|_{2} \leq \epsilon$}.
}
Then
\be{
\label{l1_err_bound}
\| x - \hat{x} \|_{1} + \lambda \| c - \hat{c} \|_{1} \leq C_1 \left( \sigma_{s}(x)_1 + \lambda \sigma_{k}(c)_1 \right) + C_2\sqrt{s+\lambda^2 k} \epsilon,
}
and
\be{
\label{l2_err_bound}
\| x - \hat{x} \|_{2} + \| c - \hat{c} \|_{2} \leq C_3 \left ( 1 + \eta^{1/4} \right )\left( \frac{\sigma_{s}(x)_1}{\sqrt{s}} + \frac{\sigma_{k}(c)_1}{\sqrt{k}} \right) + C_4 \left ( 1 + \eta^{1/4} \right ) \epsilon,
}
where the constants $C_1,C_2,C_3,C_4$ depend on $\rho$ and $\tau$ only and $\eta$ is given by
\be{
\label{eta_def}
\eta = \eta_{s,k}(\lambda) = \frac{s+\lambda^2k}{\min \{ s , \lambda^2 k \} }.
}
}
%

\prf{
We first prove \R{l1_err_bound}.
Lemma \ref{l:21_RNSP} implies that $A$ satisfies the $\ell^1$-robust null space property.  Let $S \subseteq \{1,\ldots,N\}$, $|S| \leq s$ and $T \subseteq \{1,\ldots,m\}$, $|T| \leq k$ be such that $\| x_{S^c} \|_{1} = \sigma_{s}(x)_1$ and $\| c_{T^c} \|_{1} = \sigma_{k}(c)_1$.  Then, if $v = x - \hat{x}$ and $e = c - \hat{c}$ we have
\eas{
\| x \|_{1} +\lambda \| c \|_{1} + \| v_{S^c} \|_{1} + \lambda \| e_{T^c} \|_{1} &\leq 2 \| x_{S^c} \|_{1} + \| x_{S} \|_{1} + \lambda \left ( 2 \| c_{T^c} \|_{1} + \| c_{T} \|_1 \right ) + \| \hat{x}_{S^c} \|_{1} + \lambda \| \hat{c}_{T^c} \|_{1}
\\
& \leq 2 \| x_{S^c} \|_{1} + \| v_{S} \|_{1} + \| \hat{x} \|_{1} + \lambda \left ( 2 \| c_{T^c} \|_{1} + \| e_{T} \|_{1} + \| \hat{c} \|_{1} \right ).
}
Rearranging now gives
\eas{
\| v_{S^c} \|_{1} + \lambda \| e_{T^c} \|_{1} \leq & \left ( 2 \| x_{S^c} \|_{1} + \| v_{S} \|_{1} \right ) + \lambda \left ( 2 \| c_{T^c} \|_{1} + \| e_{T} \|_{1} \right )
\\
& + \left ( \| \hat{x} \|_{1} + \lambda \| \hat{c} \|_{1} \right ) - \left ( \| x \|_{1} + \lambda \| c \|_{1} \right )
\\
& \leq 2 \left ( \| x_{S^c} \|_{1} + \lambda \| c_{T^c} \|_{1} \right ) + \left ( \| v_{S} \|_{1} + \lambda \| e_{T} \|_{1} \right ),
}
where in the second inequality we note that $\| x \|_{1} + \lambda \| c \|_{1} \geq \| \hat{x} \|_{1} + \lambda \| \hat{c} \|_{1}$ since $(x,c)$ is feasible and $(\hat{x},\hat{c})$ is a minimizer.    The $\ell^1$-robust null space property now implies that
\bes{
\| v_{S^c} \|_{1} + \lambda \| e_{T^c} \|_{1} \leq \frac{2}{1-\rho} \left ( \| x_{S^c} \|_{1} + \lambda \| c_{T^c} \|_{1} \right ) + \frac{\tau\sqrt{s+\lambda^2 k}}{1-\rho} \| A v + e \|_{2},
}
and since $\| x_{S^c} \|_{1} = \sigma_{s}(x)_1$, $\| c_{T^c} \|_{1} = \sigma_{k}(c)_1$ and
\be{
\label{Ave_eps_bound}
\| A v + e \|_{2} \leq \| A \hat{x} + \hat{c} - y \|_{2} + \| A x + c - y \|_{2} \leq 2 \epsilon,
}
we deduce that
\be{
\label{ve_ST_comp_bound}
\| v_{S^c} \|_{1} + \lambda \| e_{T^c} \|_{1} \leq\frac{2}{1-\rho} \left ( \sigma_{s}(x)_1+ \lambda \sigma_{k}(c)_1 \right ) + \frac{2 \tau}{1-\rho} \sqrt{s+\lambda^2 k}\epsilon.
}
Finally, to complete the proof of \R{l1_err_bound} we argue as follows:
\eas{
\| v \|_{1} + \lambda \| e\|_{1} & \leq \| v_{S} \|_{1} + \lambda \| e_{T} \|_{1} + \| v_{S^c} \|_{1} + \lambda \| e_{T^c} \|_{1}
\\
& \leq (1+\rho) \left ( \| v_{S^c} \|_{1} + \lambda \| e_{T^c} \|_{1} \right ) + \tau \sqrt{s+\lambda^2 k} \| A v + e \|_{2}
\\
& \leq 2 \frac{1+\rho}{1-\rho} \left (  \sigma_{s}(x)_1+ \lambda \sigma_{k}(c)_1 \right ) + \frac{4}{1-\rho} \tau \sqrt{s+\lambda^2 k} \epsilon.
}
Here, we use the $\ell^1$-robust null space property in the second step, and \R{Ave_eps_bound} and \R{ve_ST_comp_bound} in the third step.

We now consider \R{l2_err_bound}.  Writing $v = x - \hat{x}$ and $e = c - \hat{c}$ as before, let $S$ be the index of the largest $s$ elements of $v$ in absolute value and $T$ be the index set of the largest $k$ elements of $e$ in absolute value.  Define
\bes{
\theta_{v} = \min_{i \in S} | v_i|,\quad \theta_{e} = \min_{j \in T} |e_j |,\qquad \theta = \max \{ \theta_v , \theta_e / \lambda \}.
}
Then
\bes{
\| v_{S^c} \|^2_{2} + \| e_{T^c} \|^2_2 = \sum_{i \notin S} | v_i |^2 + \sum_{j \notin T } |e_j |^2 \leq \theta_v \sum_{i \notin S} | v_i | + \theta_e \sum_{j \notin T } |e_j | \leq \theta \left ( \| v_{S^c} \|_{1} + \lambda \| e_{T^c} \|_{1} \right ).
}
Now observe that $\theta_{v} \leq \| v_{S} \|_{2} / \sqrt{s}$ and $\theta_{e} \leq \| e_{T} \|_{2} / \sqrt{k}$, and therefore
\bes{
\theta \leq \frac{\sqrt{\|v_S\|^2_2+\|e_T\|^2_2}}{\min \{ \sqrt{s} , \lambda \sqrt{k} \}} \leq \frac{1}{\min \{ \sqrt{s} , \lambda \sqrt{k} \}} \left ( \frac{\rho}{\sqrt{s+\lambda^2 k}} \left ( \| v_{S^c} \|_{1} + \lambda \| e_{T^c} \|_1 \right ) + 2 \tau \epsilon \right ),
}
where in the second step we use the $\ell^2$-robust null space property and \R{Ave_eps_bound}.  Combining this with the previous estimate and using the definition of $\eta$ gives
\eas{
\| v_{S^c} \|^2_{2} + \| e_{T^c} \|^2_2  \leq &  \frac{1}{\min \{ \sqrt{s} , \lambda \sqrt{k} \}} \left ( \frac{\rho}{\sqrt{s+\lambda^2 k}} \left ( \| v_{S^c} \|_{1} + \lambda \| e_{T^c} \|_1 \right )^2 + 2 \tau \epsilon \left ( \| v_{S^c} \|_{1} + \lambda \| e_{T^c} \|_1 \right ) \right )
\\
= &
  \sqrt{\eta} \left[ \rho w^2 + 2 \tau \epsilon w \right],
}
where we have defined the non-negative scalar $w$ as
  \begin{align*}
    w \coloneqq \frac{\left\| v_{S^c} \right\|_1 + \lambda \left\| e_{T^c} \right\|_1}{\sqrt{s + \lambda^2 k}}
  \end{align*}
  Completing the square with respect to $w$ under the brackets yields
  \begin{align*}
    \| v_{S^c} \|^2_{2} + \| e_{T^c} \|^2_2 \leq \rho \sqrt{\eta} \left[ \left(w + \frac{\tau \epsilon}{\sqrt{\rho}} \right)^2 - \frac{\tau^2 \epsilon^2}{\rho} \right] 
                                            \leq \rho \sqrt{\eta} \left(w + \frac{\tau \epsilon}{\sqrt{\rho}}\right)^2
  \end{align*}
Using the $\ell^2$-robust NSP on the pair $(v,e)$ along with the above estimate, we have
\begin{align}\nonumber
  \frac{1}{\sqrt{2}} \left( \| v\|_2 + \|e\|_2 \right) \leq \sqrt{\| v \|^2_{2} + \| e \|^2_{2}} &= \sqrt{\| v_S \|^2_2 + \| e_T \|^2_2 + \| v_{S^c} \|^2_2 + \| e_{T^c} \|^2_2}
\\\nonumber
&\leq \sqrt{\| v_S \|^2_2 + \| e_T \|^2_2} + \sqrt{\| v_{S^c} \|^2_2 + \| e_{T^c} \|^2_2}
\\\nonumber
&\leq \rho w + 2 \tau \epsilon + \sqrt{\rho} \eta^{1/4} \left(w + \frac{\tau \epsilon}{\sqrt{\rho}}\right) \\\label{eq:rNSP_stable_robust-temp}
&= \sqrt{\rho} \left( \sqrt{\rho} + \eta^{1/4}\right) w + \tau \left(2 + \eta^{1/4} \right) \epsilon
\end{align}
We note that
\begin{align*}
  w &= \frac{ \left\| v_{S^c} \right\|_1 + \lambda \left\| e_{T^c} \right\|_1}{\sqrt{s + \lambda^2 k}} \leq \frac{ \left\| v \right\|_1 + \lambda \left\| e \right\|_1}{\sqrt{s + \lambda^2 k}} \\
    &\leftstackrel{\eqref{l1_err_bound}}{\leq} C_1 \left[ \frac{\sigma_s(x)_1}{\sqrt{s + \lambda^2 k}} + \lambda \frac{\sigma_k(c)_1}{\sqrt{s + \lambda^2 k}} \right] + C_2 \epsilon
  \leq C_1 \left[ \frac{\sigma_s(x)_1}{\sqrt{s}} + \frac{\sigma_k(c)_1}{\sqrt{k}} \right] + C_2 \epsilon
\end{align*}
Combining the above with \eqref{eq:rNSP_stable_robust-temp} proves \eqref{l2_err_bound}.
}

\subsection{The Restricted Isometry Property for the sparse corruptions problem}\label{sec:rip}

The robust NSP is typically difficult to prove directly.  Hence we now introduce the Restricted Isometry Property (RIP) for the sparse corruptions problem, and show that it implies the robust NSP.  Note that this has been defined previously in \cite[Defn.\ 2.1]{LiCorruptionsConstrApprox}.

\defn{
\label{d:RIPcorruptions}
Let $1 \leq s \leq N$, $1 \leq k \leq m \leq N$ and $A \in \bbC^{m \times N}$.  The $(s,k)^{\rth}$ Restricted Isometry Constant (RIC) $\delta = \delta_{s,k}$ of the matrix $A$ is the smallest constant such that
\bes{
(1-\delta) \left ( \| x \|^2_{2} + \| c \|^2_{2} \right ) \leq \| A x  + c \|^2_{2} \leq (1+\delta) \left ( \| x \|^2_{2} + \| c \|^2_{2} \right )
}
for all $x \in \Sigma_{s}$ and $c \in \Sigma_{k}$.  If $0 < \delta_{s,k} < 1$ then we say that $A$ has the Restricted Isometry Property (RIP) of order $(s,k)$.
}

Our first result is the following:

\lem{
\label{l:RIP_implies_rNSP}
Let $1 \leq s \leq N$, $1 \leq k \leq m \leq N$, $\lambda > 0$ and $A \in \bbC^{m \times N}$.  If $A$ satisfies the RIP of order $(2s,2k)$ with constant
\be{
\label{delta_cond}
\delta_{2s,2k} < \frac{1}{\sqrt{1 + \left ( \frac{1}{2 \sqrt{2}} + \sqrt{\eta} \right )^2 }},
}
where $\eta$ is as in \R{eta_def}, then $A$ satisfies the $\ell^2$-robust NSP of order $(s,k)$ with weight $\lambda$ and constants $0 < \rho < 1$ and $\tau > 0$ depending only on $\delta_{2s,2k}$.
}
The proof of this result is given next. Combining this lemma with Theorem \ref{t:rNSP_stable_robust} now yields our main result:
\thm{
\label{t:RIP_stable_robust}
Let $1 \leq s \leq N$, $1 \leq k \leq m$ and $\lambda > 0$ and suppose that $A \in \bbC^{m \times N}$ satisfies the RIP of order $(2s,2k)$ with constant $\delta_{2s, 2k}$ satisfying \eqref{delta_cond}
and $\eta$ as in \eqref{eta_def}.
Let $x \in \bbC^N$, $c \in \bbC^m$, $y \in \bbC^m$ and $\epsilon > 0$ be such that $\| A x + c - y \|_{2} \leq \epsilon$, and suppose that  $(\hat{x},\hat{c})$ is a minimizer of
\bes{
\min_{z \in \bbC^N, d \in \bbC^m} \| z \|_{1} + \lambda \| d \|_{1}\ \mbox{subject to $\| A z + d - y \|_{2} \leq \epsilon$},
}
Then
\begin{align*}
  \| x - \hat{x} \|_{1} + \lambda \| c - \hat{c} \|_{1} &\leq C_1 \left( \sigma_{s}(x)_1 + \lambda \sigma_{k}(c)_1 \right) + C_2\sqrt{s+\lambda^2 k} \epsilon,\\ 
  \| x - \hat{x} \|_{2} + \| c - \hat{c} \|_{2} &\leq C_3 \left ( 1 + \eta^{1/4} \right )\left( \frac{\sigma_{s}(x)_1}{\sqrt{s}} + \frac{\sigma_{k}(c)_1}{\sqrt{k}} \right) + C_4 \left ( 1 + \eta^{1/4} \right ) \epsilon,
\end{align*}
where the constants $C_1,C_2,C_3,C_4$ depend on $\delta_{2s,2k}$ only.
}


We now prove Lemma \ref{l:RIP_implies_rNSP}.  We first require the following:
\lem{
\label{l:disjoint_inner_product_RIP}
Let $1 \leq s \leq N$, $1 \leq k \leq m \leq N$, and let $A \in \bbC^{m \times N}$ satisfy the RIP of order $(2s, 2k)$ with constant $\delta_{2s, 2k}$.  Suppose that $x \in \Sigma_{s}$ and $c \in \Sigma_{k}$ are such that
\bes{
\nm{A x + c}^2_{2} - \left ( \| x \|^2_{2} + \| c \|^2_{2} \right ) = t  \left ( \| x \|^2_{2} + \| c \|^2_{2} \right ),
}
for some $t$ with $0 \leq |t| \leq \delta_{2s,2k}$. If $z \in \Sigma_s$ and $d \in \Sigma_k$ are orthogonal to $x$ and $c$, respectively, then
\bes{
\left | \ip{A x + c}{A z + d} \right | \leq \sqrt{\delta^2_{2s,2k} - t^2} \sqrt{\| x \|^2_{2} + \| c \|^2_{2}} \sqrt{\| z \|^2_{2} + \| d\|^2_{2}}.
}
}
\prf{
Assume that $\| x \|^2_{2} + \| c \|^2_{2} = \| z \|^2_{2} + \| d\|^2_{2} = 1$ without loss of generality.  Let $\alpha,\beta \in \bbR$ and $\gamma \in \bbC$ and notice that $\alpha x + \gamma z, \beta x - \gamma z \in \Sigma_{2s}$ and $\alpha c + \gamma d , \beta c - \gamma d \in \Sigma_{2k}$.  Therefore
\eas{
\nm{A (\alpha x + \gamma z ) +( \alpha c + \gamma d) }^2_{2} &\leq \left ( 1 + \delta_{2s,2k} \right ) \left ( \nm{\alpha x + \gamma z }^2_2 + \nm{\alpha c + \gamma d }^2_{2} \right )
\\
& = \left ( 1 + \delta_{2s,2k} \right ) \left (  \alpha^2 \left ( \nm{x}^2_2 + \nm{c}^2_2 \right ) + |\gamma|^2 \left ( \nm{z}^2_2 + \nm{d}^2_2 \right ) \right )
\\
& = \left ( 1 + \delta_{2s,2k} \right )  \left ( \alpha^2 + | \gamma |^2 \right ).
}
Note that in the second step we use orthogonality of the vectors $x$ and $z$ and $c$ and $d$.  Similarly,
\bes{
\nm{A (\beta x - \gamma z ) +( \beta c - \gamma d) }^2_{2} \geq \left ( 1 - \delta_{2s,2k} \right ) \left ( \beta^2 + |\gamma|^2 \right ).
}
Subtracting the second equation from the first gives
\ea{
\nm{A (\alpha x + \gamma z ) +( \alpha c + \gamma d) }^2_{2} -& \nm{A (\beta x - \gamma z ) +( \beta c - \gamma d) }^2_{2} \nn
\\
& \leq \left ( 1 + \delta_{2s,2k} \right )  \left ( \alpha^2 + | \gamma |^2 \right ) - \left ( 1 - \delta_{2s,2k} \right ) \left ( \beta^2 + |\gamma|^2 \right ) \nn
\\
& = \delta_{2s,2k} \left ( \alpha^2 + \beta^2 + 2 | \gamma |^2 \right ) + \alpha^2 - \beta^2. \label{diff_upper_bd}
}
On the other hand
\eas{
\nm{A (\alpha x + \gamma z ) +( \alpha c + \gamma d) }^2_{2} - &\nm{A (\beta x - \gamma z ) +( \beta c - \gamma d) }^2_{2}
\\
=& \ \alpha^2 \nm{A x + c }^2_{2} + | \gamma |^2 \nm{A z + d }^2_{2} + 2 \Re \ip{\alpha(Ax+c)}{\gamma(A z + d ) }
\\
& \ - \beta^2 \nm{A x + c }^2_{2} - | \gamma |^2 \nm{A z + d }^2_{2} + 2 \Re \ip{\beta(Ax+c)}{\gamma(Az+d)}
\\
=& \left ( \alpha^2 - \beta^2 \right ) \nm{A x + c }^2_{2} + 2 (\alpha+\beta) \Re \left ( \bar{\gamma} \ip{Ax+c}{Az+d} \right )
\\
= & \left ( \alpha^2 - \beta^2 \right ) (1+t)  +  2 (\alpha+\beta) \Re \left ( \bar{\gamma} \ip{Ax+c}{Az+d} \right ).
}
Combining this with \R{diff_upper_bd} gives
\bes{
\left ( \alpha^2 - \beta^2 \right ) (1+t)  +  2 (\alpha+\beta) \Re \left ( \bar{\gamma} \ip{Ax+c}{Az+d} \right ) \leq \delta_{2s,2k} \left ( \alpha^2 + \beta^2 + 2 | \gamma |^2 \right ) + \alpha^2 - \beta^2.
}
Now let $\gamma$ be such that $| \gamma | = 1$ and $\Re \left ( \bar{\gamma} \ip{Ax+c}{Az+d} \right ) = | \ip{Ax+c}{Az+d}|$.  Then, after rearranging, we get
\bes{
| \ip{Ax+c}{Az+d}| \leq \frac{\left ( \delta_{2s,2k} - t\right ) \alpha^2 +\left ( \delta_{2s,2k} +t\right ) \beta^2 + 2 \delta_{2s,2k} }{2(\alpha+\beta)} .
}
We now seek values $\alpha$ and $\beta$ which minimize the right-hand side of this expression.  If $t = \delta_{2s,2k}$ then the minimal value $0$ is attained by setting $\beta =0 $ and letting $\alpha \rightarrow \infty$.  Conversely, if $t < \delta_{2s,2k}$ the minimal value is attained when $\alpha = \sqrt{\frac{\delta_{2s,2k}+t}{\delta_{2s,2k}-t}}$ and $\beta = \frac{1}{\alpha}$.  This gives
\bes{
| \ip{Ax+c}{Az+d}| \leq \sqrt{\delta^2_{2s,2k} - t^2},
}
which completes the proof.
}

\prf{[Proof of Lemma \ref{l:RIP_implies_rNSP}]
Let $x \in \bbC^N$ and $c \in \bbC^m$.  To prove the $\ell^2$-robust NSP for $A$ it is enough to show that \R{l2_rNSP_def} holds when $S = S_0$ is the index set of the $s$ largest coefficients of $x$ in absolute value and $T = T_0$ is the set of the $k$ largest values of $c$ in absolute value.  Given $S_0$, let $S_1$ be the index set of the next $s$ largest coefficients of $x$ in absolute value, $S_2$ be the index set of the next $s$ largest coefficients and so on.  Define $T_1,T_2,\ldots$ in a similar way.  We now have the following:
\ea{
\nm{A x_{S_0} + c_{T_0}}^2 &= \ip{A x_{S_0} + c_{T_0} }{A x_{S_0} + c_{T_0} } \nn
\\
& = \ip{A x_{S_0} + c_{T_0} }{A x + c} - \sum_{j \geq 1} \ip{A x_{S_0} + c_{T_0} }{A x_{S_j} + c_{T_j}}.\label{RIP_sum_split}
}
Let $0 \leq |t| \leq \delta_{2s,2k}$ be such that
\be{
\label{t_def}
\nm{A x_{S_0} + c_{T_0}}^2_2 = (1+t) \left ( \nm{x_{S_0}}^2_2 + \nm{c_{T_0}}^2_2 \right ),
}
and note that this gives
\be{
\label{RIP_sum_split_term1}
\left | \ip{A x_{S_0} + c_{T_0} }{A x + c} \right | \leq \sqrt{1+t} \sqrt{ \nm{x_{S_0}}^2_2 + \nm{c_{T_0}}^2_2 } \nm{A x + c}_{2}.
}
For the second term of \R{RIP_sum_split}, we use the disjointness of $S_0$ and $S_j$ and $T_0$ and $T_j$ for $j \geq 1$ in combination with Lemma \ref{l:disjoint_inner_product_RIP} to get
\ea{
\left | \sum_{j \geq 1} \ip{A x_{S_0} + c_{T_0} }{A x_{S_j} + c_{T_j}} \right | &\leq \sqrt{\delta^2_{2s,2k} - t^2} \sqrt{ \nm{x_{S_0}}^2_2 + \nm{c_{T_0}}^2_2 }  \sum_{j \geq 1}  \sqrt{ \nm{x_{S_j}}^2_2 + \nm{c_{T_j}}^2_2 } \nn
\\
& \leq \sqrt{\delta^2_{2s,2k} - t^2} \sqrt{ \nm{x_{S_0}}^2_2 + \nm{c_{T_0}}^2_2 }  \left ( \sum_{j \geq 1} \| x_{S_j} \|_{2} + \sum_{j \geq 1} \nm{c_{T_j}}_2 \right ). \label{RIP_sum_split_term2}
}
Let $x^{+}_{j}$ and $x^{-}_{j}$ be the largest entries of $x_{S_j}$ in absolute value.  Then, by \cite[Lem.\ 6.14]{FoucartRauhutCSbook}, we have
\eas{
\sum_{j \geq 1} \| x_{S_j} \|_{2} & \leq \sum_{j \geq 1} \left ( \frac{\nm{x_{S_j}}_1}{\sqrt{s}} + \frac{\sqrt{s}}{4} \left ( x^{+}_{j} - x^{-}_{j} \right ) \right )
\\
& \leq \frac{\nm{x_{S^c_0}}_1}{\sqrt{s}} + \frac{\sqrt{s}}{4} \sum_{j \geq 1} \left ( x^{+}_{j} - x^{+}_{j+1} \right )
\leq \frac{\nm{x_{S^c_0}}_1}{\sqrt{s}} + \frac{\sqrt{s}}{4} x^{+}_{1}
\leq \frac{\nm{x_{S^c_0}}_1}{\sqrt{s}} + \frac{1}{4} \| x_{S_0} \|_{2}.
}
Similarly,
\bes{
\sum_{j \geq 1} \nm{c_{T_j}}_2 \leq \frac{\nm{c_{T^c_0}}_1}{\sqrt{k}} + \frac{1}{4} \nm{c_{T_0}}_2 \leq  \frac{\lambda \nm{c_{T^c_0}}_1}{\lambda \sqrt{k}} + \frac{1}{4} \nm{c_{T_0}}_2 ,
}
which gives
\bes{
\sum_{j \geq 1} \| x_{S_j} \|_{2} + \sum_{j \geq 1} \nm{c_{T_j}}_2 \leq \frac{1}{\min \left \{ \sqrt{s} , \lambda \sqrt{k} \right \} } \left ( \nm{x_{S^c_0}}_1 + \lambda \nm{c_{T^c_0}}_1 \right ) + \frac14 \left ( \| x_{S_0} \|_{2} + \nm{c_{T_0}}_2 \right ).
}
Therefore, combining this with \R{RIP_sum_split}, \R{t_def}, \R{RIP_sum_split_term1} and \R{RIP_sum_split_term2} yields
\eas{
(1+t) &\sqrt{\nm{x_{S_0}}^2_2 + \nm{c_{T_0}}^2_2} \leq  \sqrt{1+t} \nm{A x + c}_{2}
\\
& + \sqrt{\delta^2_{2s,2k} - t^2} \left ( \frac{1}{\min \left \{ \sqrt{s} , \lambda \sqrt{k} \right \} } \left ( \nm{x_{S^c_0}}_1 + \lambda \nm{c_{T^c_0}}_1 \right ) + \frac14 \left ( \| x_{S_0} \|_{2} + \nm{c_{T_0}}_2 \right ) \right ).
}
Consider the function $g(t) = \frac{\delta^2_{2s,2k} - t^2}{(1+t)^2}$, where $0 \leq t \leq \delta_{2s,2k}$.  This function attains its maximum value at $t = - \delta^2_{2s,2k}$ and takes value $\frac{\delta^2_{2s,2k}}{1-\delta^2_{2s,2k}}$ there.  Additionally $\frac{1}{\sqrt{1+t}} \leq \frac{1}{\sqrt{1-\delta_{2s,2k}}}$.  Hence we get
\eas{
\sqrt{\nm{x_{S_0}}^2_2 + \nm{c_{T_0}}^2_2} &\leq \frac{1}{\sqrt{1-\delta_{2s,2k}}} \nm{Ax+c}_2
\\
 +& \frac{\delta_{2s,2k}}{\sqrt{1-\delta^2_{2s,2k}}} \left ( \frac{1}{\min \left \{ \sqrt{s} , \lambda \sqrt{k} \right \} } \left ( \nm{x_{S^c_0}}_1 + \nm{c_{T^c_0}}_1 \right ) + \frac14 \left ( \| x_{S_0} \|_{2} + \nm{c_{T_0}}_2 \right ) \right ).
}
After noting that $\| x_{S_0} \|_{2} + \nm{c_{T_0}}_2 \leq \sqrt{2} \sqrt{\nm{x_{S_0}}^2_2 + \nm{c_{T_0}}^2_2}$ and rearranging, we obtain
\bes{
  \sqrt{\nm{x_{S_0}}^2_2 + \nm{c_{T_0}}^2_2} \leq \frac{\rho}{\sqrt{s + \lambda^2 k}} \left ( \nm{x_{S^c_0}}_1 + \nm{c_{T^c_0}}_1 \right ) + \tau \nm{A x + c}_{2},
}
where
\be{\label{eq:rho-tau-def}
\rho = \frac{2 \sqrt{2} \delta_{2s,2k}}{2\sqrt{2} \sqrt{1-\delta^2_{2s,2k}} - \delta_{2s,2k}} \sqrt{\eta} ,\quad \tau = \frac{2 \sqrt{2} \sqrt{1+\delta_{2s,2k}} }{2\sqrt{2} \sqrt{1-\delta^2_{2s,2k}} - \delta_{2s,2k}}.
}
To complete the proof we note that $\rho, \tau > 0$ provided $\delta_{2s,2k} < \sqrt{8/9}$.  This holds by assumption, since $\eta \geq 2$ and therefore the condition \R{delta_cond} implies that $\delta_{2s,2k} < \sqrt{8/33} < \sqrt{8/9}$.  Also, after rearranging we see that $\rho < 1$ if
\bes{
\left ( 1 + \left ( \frac{1}{2 \sqrt{2}} + \sqrt{\eta} \right )^2 \right ) \delta^2_{2s,2k} < 1,
}
which again holds by assumption.
}

\rem{
\label{r:RIPinlevels}
The RIP for the sparse corruptions problem is a special case of the RIP in levels (RIPL), introduced in \cite{BastounisHansen}.  The RIPL applies to vectors that are sparse in levels; namely, having different amounts of sparsity in different (but fixed) sections of the vector.  In the context of the sparse corruptions problem, this corresponds to the concatenated vector $z = [ x ; c ]$, which is $s$-sparse in its first $N$ entries and $k$-sparse in its remaining $m$ entries.  As a general tool, sparsity in levels has been used in the context of compressive imaging \cite{AHPRBreaking,OptimalSamplingQuest,AsymptoticCS}, radar \cite{Dorsch2016} and multi-sensor acquisition \cite{AdcockChunParallel}.  It is interesting that the same model also occurs naturally in the, seemingly unrelated, sparse corruptions problem.  We note in passing that Theorems \ref{t:rNSP_stable_robust} and \ref{t:RIP_stable_robust} follow a similar approach to that of \cite{BastounisHansen} with some changes made to incorporate the weighted optimization problem.  
}

\subsection{Matrices that satisfy the RIP for sparse corruptions}\label{sec:rip-2}

We first recall the classical RIP for sparse vectors:

\defn{
Let $1 \leq s \leq N$ and $A \in \bbC^{m \times N}$.  The $s^{\rth}$ Restricted Isometry Constant (RIC) $\delta = \delta_{s}$ of the matrix $A$ is the smallest constant such that
\bes{
(1-\delta)\nm{x}^2_2 \leq \| A x \|^2_{2} \leq (1+\delta) \nm{x}^2_2,
}
for all $x \in \Sigma_{s}$.  If $0 < \delta_{s} < 1$ then we say that $A$ has the Restricted Isometry Property (RIP) of order $s$.
}

To distinguish it from the RIP for the sparse corruptions problem (Definition \ref{d:RIPcorruptions}), we shall refer to this as the \textit{RIP for sparse vectors}.

\lem{
\label{l:deltask_deltas_sigmask}
Let $1 \leq s \leq N$, $1 \leq k \leq m$, $A \in \bbC^{m \times N}$ and define
\be{
\label{sigma_def}
\sigma_{s,k} = \max_{\substack{S \subseteq \{1,\ldots,N\}, |S| = s \\ T \subseteq \{1,\ldots,m \}, |T| = k }} \| A_{S,T} \|_{2},
}
where $A_{S,T} \in \bbC^{|T| \times |S|}$ is the submatrix of $A$ with entries $\{ A_{ij} \}_{i \in T, j \in S}$.  Suppose that $A$ has the RIP for sparse vectors with constant $\delta_s$ and that $\sigma_{s,k} < \sqrt{1-\delta_{s}}$.  Then $A$ has the RIP of order $(s,k)$ for the sparse corruptions problem with constant
\bes{
\delta_{s,k} =  \frac{\delta_s + \sqrt{\delta^2_s+4 \sigma^2_{s,k}}}{2}.
}
In other words,
\bes{
\left ( 1 - \delta_{s,k} \right ) \left ( \nm{x}^2_{2} + \nm{c}^2_2 \right ) \leq \nm{A x + c}^2_{2} \leq \left ( 1 + \delta_{s,k} \right ) \left ( \nm{x}^2_{2} + \nm{c}^2_2 \right )
}
for all $x \in \Sigma_{s}$ and $c \in \Sigma_k$.
}
\prf{
Let $x \in \Sigma_{s}$ and $c \in \Sigma_k$ and write $S = \supp(x)$ and $T = \supp(c)$.  Then
\bes{
\nm{A x + c}^2_{2} = \nm{A x}^2_{2} + \nm{c}^2_{2} + 2 \Re \ip{A_{S,T} x}{c}.
}
By Young's inequality
\bes{
2\left | \ip{A_{S,T} x}{c} \right | \leq 2 \| A_{S,T} \|_{2} \nm{x}_2 \nm{c}_2 \leq \| A_{S,T} \|_{2} \left ( \nm{x}^2_2 / \epsilon + \epsilon \nm{c}^2_2 \right ),
}
for any $\epsilon > 0$.  Hence
\bes{
\left ( 1 - \delta_s - \sigma_{s,k}/\epsilon \right ) \nm{x}^2_2 + \left ( 1 - \sigma \epsilon \right ) \nm{c}^2_{2} \leq \nm{A x + c }^2_{2} \leq \left ( 1 + \delta_s + \sigma_{s,k}/\epsilon \right ) \nm{x}^2_2 + \left ( 1 + \sigma \epsilon \right ) \nm{c}^2_{2}.
}
Solving the equation $\delta_s + \sigma_{s,k}/\epsilon = \sigma_{s,k} \epsilon$ yields the value $\epsilon = \frac{\delta_s + \sqrt{\delta^2_s + 4 \sigma^2_{s,k}}}{2 \sigma}$, and substituting this value of $\epsilon$ into the previous expression yields the proof.
}

This result shows that any matrix satisfying the RIP for sparse vectors also satisfies the RIP for the sparse corruptions problem, provided the all $k \times s$ submatrices have small spectral norm.

\subsubsection{Gaussian random matrices}
Gaussian random matrices in the context of the sparse corruptions problem were considered in \cite{LiCorruptionsConstrApprox}.  The following result essentially recaps the main result for this case given therein.  We include a short proof for completeness:

\thm{
\label{t:GaussCorruption}
Let $0 < \delta , \epsilon < 1$, $1 \leq s \leq$, $1 \leq k \leq m$ and suppose that
\ea{
m &\gtrsim \delta^{-2} \left ( s \cdot \log(2N/s) + \log(2 \epsilon^{-1}) \right ), \label{mGaussRIP}
\\
m &\gtrsim \delta^{-2} \cdot k \cdot \log(\delta^{-1}). \label{mGaussSigma}
}
Let $A \in \bbC^{m \times N}$ be a matrix whose entries are independent Gaussian random variables with mean zero and variance $1$.  Then with probability at least $1-\epsilon$, the matrix $\frac{1}{\sqrt{m}} A$ has the RIP for the sparse corruptions problem of order $(s,k)$ with constant $\delta_{s,k} \leq \delta$.
}
\prf{
Lemma \ref{l:deltask_deltas_sigmask} asserts that $A$ has the RIP of order $(s,k)$ for the sparse corruptions problem with constant $\delta_{s,k} \leq \delta$ provided (i) $A$ has the RIP of order $s$ with $\delta_{s} \leq \delta/\sqrt{2}$ and (ii) the constant $\sigma_{s,k}$ defined in \R{sigma_def} satisfies $\sigma_{s,k} \leq \delta/(2 \sqrt{2})$.  Hence, by the union bound it suffices to show that \R{mGaussRIP} and \R{mGaussSigma} imply both (i) and (ii) separately with probabilities at least $1-\epsilon/2$.  Due to a standard result in compressed sensing (see, for example, \cite[Thm.\ 9.2]{FoucartRauhutCSbook}), property (i) holds with probability at least $1-\epsilon/2$ whenever the condition \R{mGaussRIP} is satisfied.  We now consider property (ii).  First, notice that $\sigma_{s,k}$ is increasing in $k$.  Therefore, we may assume that $ k \asymp \delta^2 \cdot m$, i.e.\ $k \gtrsim \delta^{2} \cdot m$ and $k \lesssim \delta^{2} \cdot m$.  Fix subsets $S \subseteq \{1,\ldots,N\}$ and $T \subseteq \{1,\ldots,m\}$ with $|S| = s$ and $|T|  =k$.  Then, due to a known result for singular values of random Gaussian matrices (see, for example, \cite[Cor.\ 5.35]{Vershynin:bookCh}), we have
\bes{
\bbP \left ( \nm{A_{S,T}}_2 \geq \sqrt{s} + \sqrt{k} + t \right ) \leq 2 \exp(-t^2/2).
}
The conditions \R{mGaussRIP} and \R{mGaussSigma} imply that $\sqrt{s/m} \leq \delta  /(6\sqrt{2})$ and $\sqrt{k/m} \leq \delta  /(6\sqrt{2})$.  Hence, by the union bound
\eas{
\bbP \left ( \sigma_{s,k} > \delta / (2 \sqrt{2}) \right )  \leq \left ( \begin{array}{c} N \\ s \end{array} \right ) \left ( \begin{array}{c} m \\ k \end{array} \right ) \exp(-m\delta^2/48) \leq \left ( \frac{\E N}{s} \right )^s \left ( \frac{\E m}{k} \right )^k \exp(-m\delta^2/48).
}
In particular, $\bbP \left ( \sigma_{s,k} > \delta / (2 \sqrt{2}) \right )  \leq \epsilon/2$ provided
\bes{
m \geq 48 \cdot \delta^{-2} \left ( s \log(\E N/s) + k \log(\E m/k) + \log(2 \epsilon^{-1}) \right ).
}
Since $k \asymp \delta^2 \cdot m$, we have $\log(\E m / k) \lesssim \log(2 \delta^{-1})$.  Hence this condition is implied by \R{mGaussRIP} and \R{mGaussSigma}.  This establishes property (ii) and completes the proof.
}

This result asserts that Gaussian random matrices can recover a fixed fraction $k/m \leq c$ of corruptions (see \R{mGaussSigma}) and (up to constants) the same level of sparsity $s$ as in the uncorrupted case (see \R{mGaussRIP}).

\subsubsection{Bounded orthonormal systems}\label{sec:bos}

Gaussian random matrices, while mathematically appealing, are of little relevance to multivariate approximation using Polynomial Chaos expansions.  In this case, a more suitable framework is that of bounded orthonormal systems (see, for example, \cite[Chpt.\ 12]{FoucartRauhutCSbook}):

Let $D$ be a domain with a probability measure $\nu$ and $\phi_1,\ldots,\phi_N$ be an orthonormal system of complex-value functions in $L^2(D)$.  Recall that this system is bounded if
\bes{
\nm{\phi_i }_{L^\infty} = \sup_{\xi \in D} | \phi_i(\xi) | \leq K
}
Given such a system, we construct the measurement matrix $A$ as
\be{
\label{ABOS}
A = \frac{1}{\sqrt{m}} \left \{ \phi_j(\xi_i) \right \}^{m,N}_{i=1,j=1} \in \bbC^{m \times N},
}
where $t_i$ are drawn independently at random according to the probability measure $\nu$.  

\thm{
\label{t:RIP_BOS}
Let $A \in \bbC^{m \times N}$ be the matrix of a bounded orthonormal system, $1 \leq s \leq N$ and $0 < \delta , \epsilon < 1$.  If
\bes{
m \gtrsim \delta^{-2} \cdot s \cdot \left ( \log^3(2s) \cdot \log(2N) + \log(\epsilon^{-1}) \right ),
}
then $A$ satisfies the RIP for sparse vectors with probability at least $1-\epsilon$.
}

We remark in passing that the logarithmic dependence in $s$ can be improved by one power, at the expense of a larger factor in $\delta^{-1}$ \cite{ChkifaDownwardsCS}.  However, this may not be best for the purposes of this paper, since in view of Theorem \ref{t:RIP_stable_robust}, $\delta^{-2}$ scales linearly in the parameter $\eta$ (see next).   

The following lemma estimates the constant $\sigma_{s,k}$ for matrices of the form \R{ABOS}:
\lem{
\label{l:BOSsigmaEst}
Let $A \in \bbC^{m \times N}$ be the matrix of a bounded orthonormal system, $1 \leq s,k \leq N$ and $\sigma_{s,k}$ be as in \R{sigma_def}.  Then
\bes{
\sigma_{s,k} \leq \sqrt{\frac{K^2 s k}{m}}.
}
}
\prf{
Fix subsets $S \subseteq \{1,\ldots,N\}$, $|S|=s$ and $T \subseteq \{1,\ldots,m\}$, $|T|=k$ and let $x \in \bbC^N$ and $c \in \bbC^m$ with $\supp(x) = S$ and $\supp(c) = T$.  Then
\eas{
\left | c^* A x \right |^2 &= \frac{1}{\sqrt{m}} \left | \sum_{i \in T} \overline{c_i} \sum_{j \in S} \phi_j(t_i) x_j \right | 
\\
& \leq \frac{1}{\sqrt{m}} \max_{i =1,\ldots,m} \left | \sum_{j \in S} \phi_j(t_i) x_j \right | \sum_{i \in T} | c_i | 
\leq \frac{K}{\sqrt{m}}  \| x \|_{1} \| c \|_{1}
\leq \sqrt{\frac{K^2 s k}{m}} \| x \|_{2} \| c \|_{2}.
}
Hence $\| P_{T} A P_{S} \|_{2} \leq \sqrt{\frac{K^2 s k}{m}}$.  This now gives the result.
}

With this in hand, we now deduce the following result:
\thm{
\label{t:BOS-RIP}
Let $1 \leq s \leq N$, $1 \leq k \leq m$, $0 < \delta,\epsilon < 1$ and suppose that
\be{
\label{mRIP1}
m \gtrsim \delta^{-2} \cdot K^2 \cdot s \cdot \left ( \log^3(2s) \cdot \log(2N) + \log(\epsilon^{-1}) \right ),
}
and 
\bes{
\label{mRIP2}
m \geq 8 \cdot \delta^{-2} \cdot K^2 \cdot s \cdot k.
}
Then, with probability at least $1-\epsilon$, $A$ has the RIP of order $(s,k)$ for the sparse corruptions problem with constant $\delta_{s,k} \leq \delta$.
}
\prf{
Theorem \ref{t:RIP_BOS} and \R{mRIP1} imply that $A$ has the RIP of order $s$ with $\delta_{s} \leq \delta/\sqrt{2}$ with probability at least $1-\epsilon$.  Moreover, Lemma \ref{l:BOSsigmaEst} and \R{mRIP2} imply that $\sigma_{s,k} \leq \delta/(2 \sqrt{2})$.  We now apply Lemma \ref{l:deltask_deltas_sigmask}.
%
}

\rem{
This result asserts that the number of corruptions that can be tolerated is a fraction of $m/s$.  This is inferior to the case of Gaussian random measurements, where Theorem \ref{t:GaussCorruption} gives that a fraction of $m$ corruptions are permitted.  We conjecture, however, that a similar estimate can be proved for the bounded orthonormal systems case -- indeed, a nonuniform recovery result of this form was proved in \cite{LiCorruptionsConstrApprox} for the case of exactly sparse coefficients $x$ and corruptions $c$ with random sign sequences -- albeit with a substantially more sophisticated argument than the proof of Theorem \ref{t:GaussCorruption}.  In particular, while estimates for the singular values of matrices of bounded orthonormal systems are known \cite{Vershynin:bookCh}, they are more stringent than those for Gaussian random matrices.  Using these estimates and arguing via the union bound (as in the proof of Theorem \ref{t:GaussCorruption}) unfortunately results in an estimate similar to \R{mRIP2}.  We also note in passing that while there exist RIP estimates for quite general matrices under the sparsity in levels model \cite{LiAdcockRIP} (see Remark \ref{r:RIPinlevels}), these unfortunately do not apply to the setup of the sparse corruptions problem.
We therefore leave the problem of improving \R{mRIP2} for future work.  
}



\subsection{Strategy for choosing $\lambda$}\label{ss:lambda-strategy}
Regardless of the matrix $A$, our main theorems (Theorems \ref{t:RIP_stable_robust} and \ref{t:RIP_BOS}) suggest an optimal strategy for choosing the parameter $\lambda$.  Notice that the restricted isometry constant $\delta$ enters into the measurement condition in Theorem \ref{t:RIP_BOS} as $\delta^{-2}$.  Since Theorem \ref{t:RIP_stable_robust} requires that \R{delta_cond} holds, the measurement condition contains a factor that is at least as large as 
\bes{
 1 + \left ( \frac{1}{2 \sqrt{2}} + \sqrt{\eta} \right )^2.
}
We wish to minimize this factor so as to reduce the measurement condition as much as possible.
This can be done by minimizing $\eta$, which in turn yields the theoretically-optimal scaling
\be{
\label{eq:lambda-opt}
\lambda = \sqrt{\frac{s}{k}}.
}
Notice that this gives the value $\eta = 2$.  In particular, the condition \R{delta_cond} becomes
\be{\label{eq:delta-condition-1}
\delta_{2s,2k} < \sqrt{8/33} \approx 0.492,
}
with right-hand side independent of $s$ and $k$.  We remark in passing that the choice \R{eq:lambda-opt} is implicitly made in \cite{LiCorruptionsConstrApprox}.  However, the condition given in \cite[Lem.\ 2.3]{LiCorruptionsConstrApprox} is $\delta_{2s,2k} < 1/18 \approx 0.056$  which is significantly more stringent than \R{eq:lambda-opt}.  Moreover, \cite{LiCorruptionsConstrApprox} only considers exact sparsity, whereas Theorem \ref{t:RIP_stable_robust} also treats the case of stable recovery of inexactly sparse coefficients and corruptions.

%

\section{Numerical experiments}\label{sec:results}

We divide our numerical results into two main sections. The goal of Section \ref{sec:results-algorithm} is to study the behavior of numerical algorithms in the context of the theoretical estimates presented earlier. In particular, we investigate the influence that the regularization parameter $\lambda$ has on recovery properties. We confine these investigations to problems with manufactured sparsity so that systematic studies may be carried out. The lessons learned from these studies allow us to formulate and propose an iteratively reweighted alternative to the one-time optimization \eqref{l1_lambda_recovery}. \edits{Note that none of our theoretical error estimates apply to algorithms with weighted norms. However, weighted $\ell^1$ schemes can provide empirically superior results, e.g., \cite{KarniadakisUQCS}. Thus, we explore weighted algorithms because their use is natural from a practical point of view, but is not in the scope of our theoretical analysis.} Our simulations in this section use the SPGL1 package \cite{spgl1:2007,BergFriedlander:2008}.

The second collection of results, Section \ref{sec:results-pce}, focuses on more practical scenarios in scientific computing, dealing with recovery of sparse or compressible polynomial Chaos expansions of solutions to parameterized differential equations. Here we use the algorithmic lessons learned from Section \ref{sec:results-algorithm} to illustrate the efficacy and fault-tolerance of our approaches on realistic problems in the presence of measurement corruptions.

\subsection{Recovery of manufactured solutions with sparse corruptions}\label{sec:results-algorithm}

This section is primarily concerned with the generation of phase recovery diagrams for the sparse corruptions problem. \edits{In particular, our tests here are not necessarily motivated by sensing matrices and corruptions from function approximation, but instead are designed to understand behavior of the algorithms.} The following standard experiment for accomplishing this is carried out: We fix the number of measurements $m$ and the dictionary size $N$, and we vary the signal sparsity $s$ and the number of measurement corruptions $k$. For each $s$ and $k$ we generate an $s$-sparse signal $x$, and for a given model of a measurement matrix $A$, we generate $m$ measurements $y$ from the signal $x$, and subsequently corrupt (highly pollute) these measurements with a $k$-sparse vector $c$, whose non-zero entries are $C Z$, where $Z$ is a random draw from a certain probability distribution and $C > 0$ is a scaling constant. In this test, $Z$ is a standard normal random variable and $C = 1$.

We then run the recovery algorithm \eqref{l1_lambda_recovery} for a given value of $\lambda$, producing a recovered signal $\widehat{x}$ and measurement corruption vector $\widehat{c}$. We define the recovery as successful if $\left\| x - \widehat{x} \right\|^2 + \left\| c - \widehat{c} \right\|^2 < \epsilon_{\mathrm{tol}}$. In this test, we set the success tolerance to be $\epsilon_{\mathrm{tol}} = 10^{-4}$.

In the test above, the generation of $x$, and of $y$, and of $c$, are statistically independent\footnote{Measurement corruptions are generated as iid standard normal random variables, and support indices in a sparse vector are generated using the uniform probability law (draws without replacement) on the index set.}. For each $s$ and $k$, the above procedure is run $T \in \bbN$ times with independent draws, and an empirical estimate of the probability of ``success" is computed. In the phase transitions plots below, we use $T = 10$ simulations.

The phase transitions color each pixel, corresponding to a particular value of $s$ and $k$, according to the empirical success probability. The phase transition axes are $s/m$ and $k/m$, and thus each ranges in the interval $[0, 1]$, but we truncate to $[0, 0.5]$ in our plots because this region is sufficient to illustrate behavior.  We consider the following two models of measurement matrix $A$:
\bull{
\item Model 1: a Gaussian random matrix
\item Model 2: a randomly-subsampled Discrete Fourier Transform (DFT) matrix
}
Note that Model 2 is an example of a bounded orthonormal system.  We compare several different choices of $\lambda$ for each model.


\subsubsection{Phase transition plots for fixed $\lambda$}

\begin{figure}
\begin{center}
  \resizebox{!}{0.43\textheight}{
    \begin{tabular}{ccc}

      \includegraphics[width=0.33\textwidth]{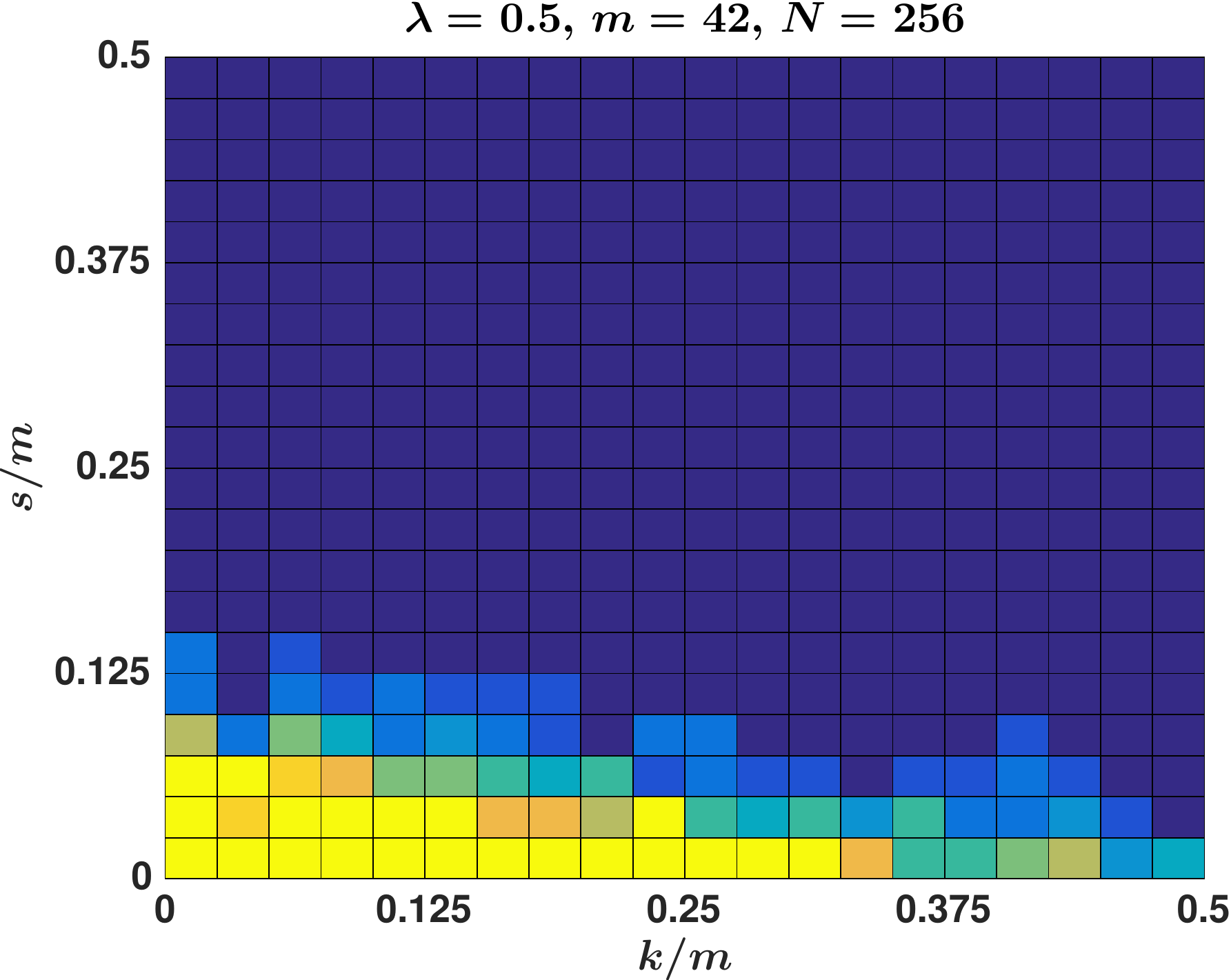}&
      \includegraphics[width=0.33\textwidth]{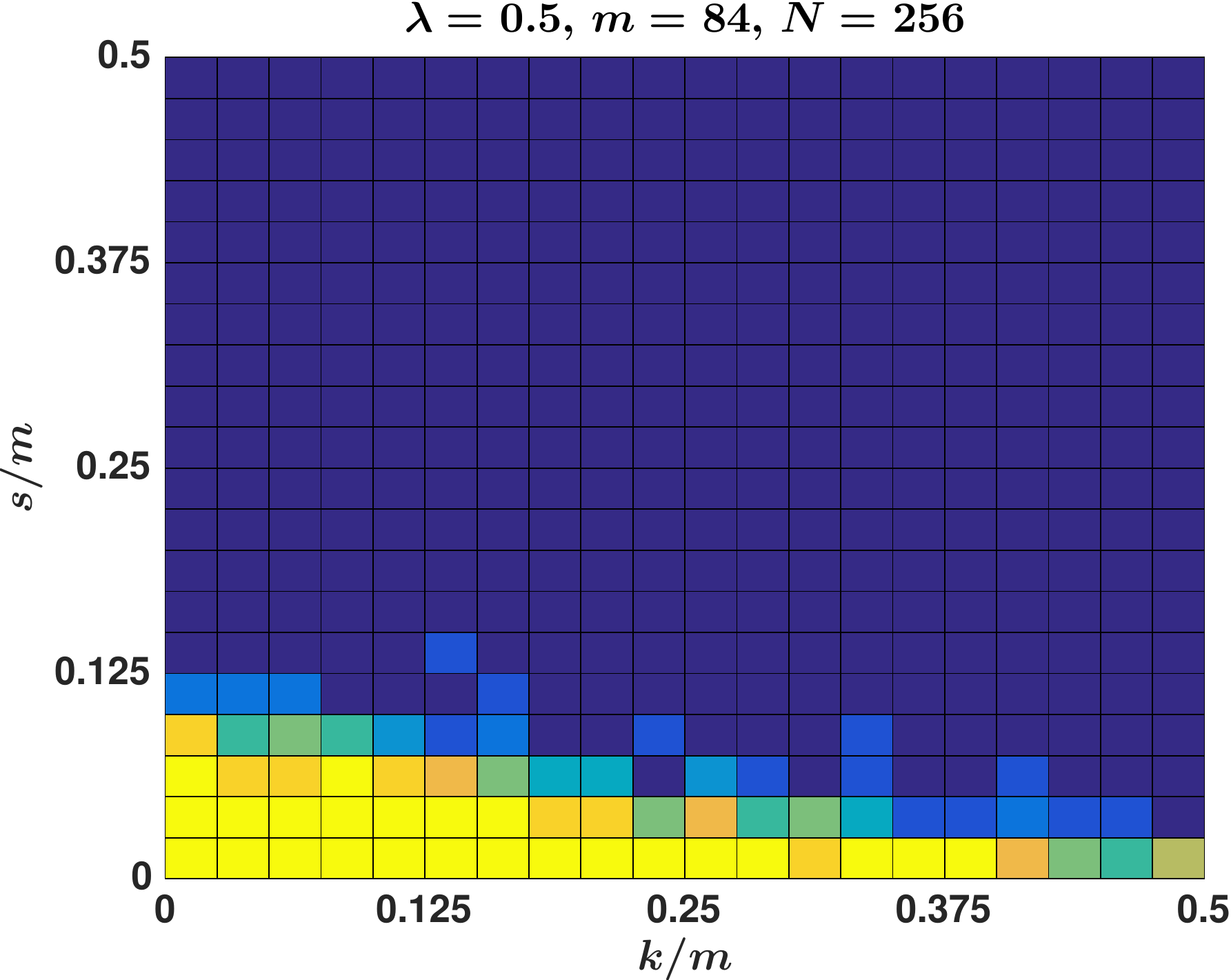}&
      \includegraphics[width=0.33\textwidth]{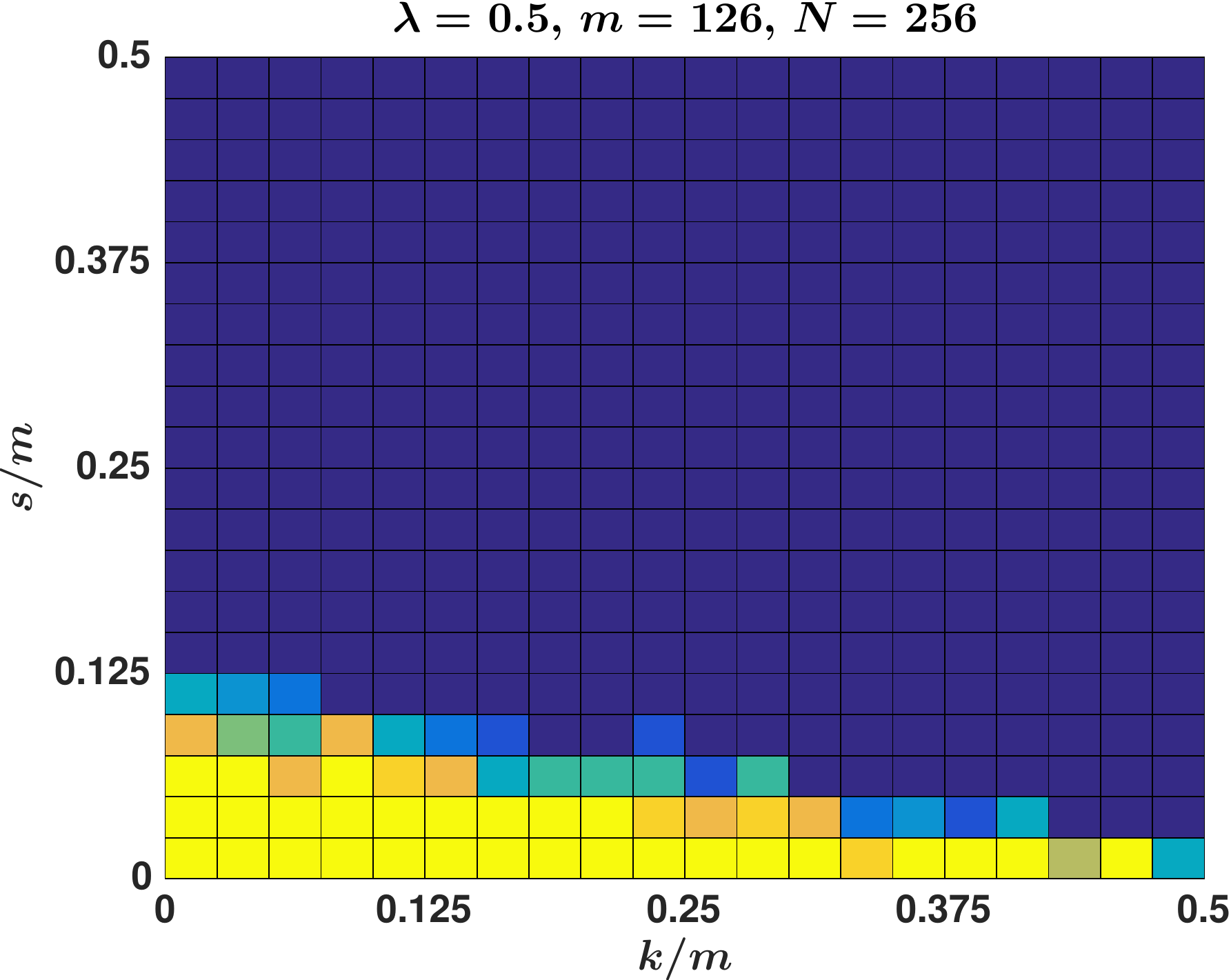}\\
      \includegraphics[width=0.33\textwidth]{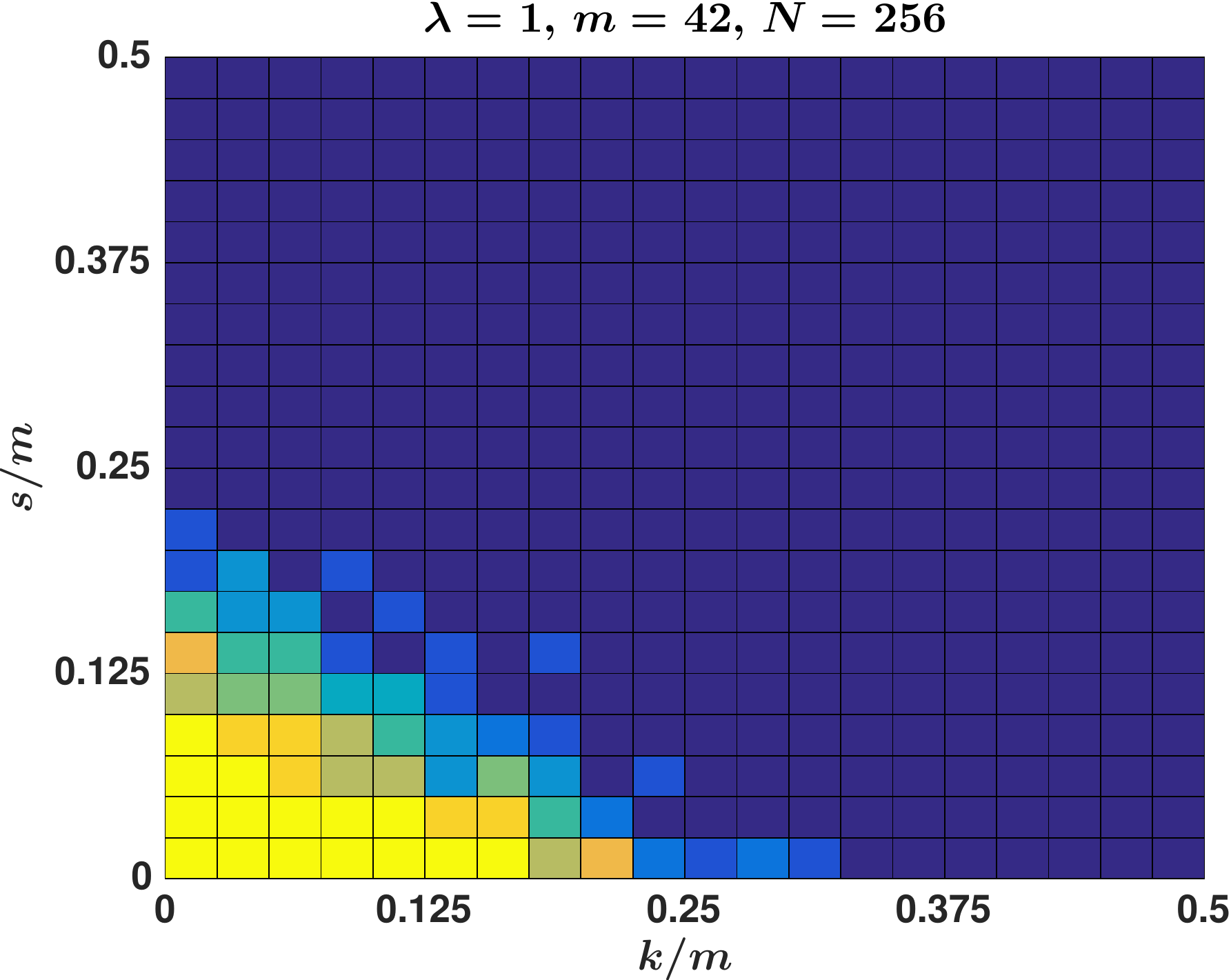}&
      \includegraphics[width=0.33\textwidth]{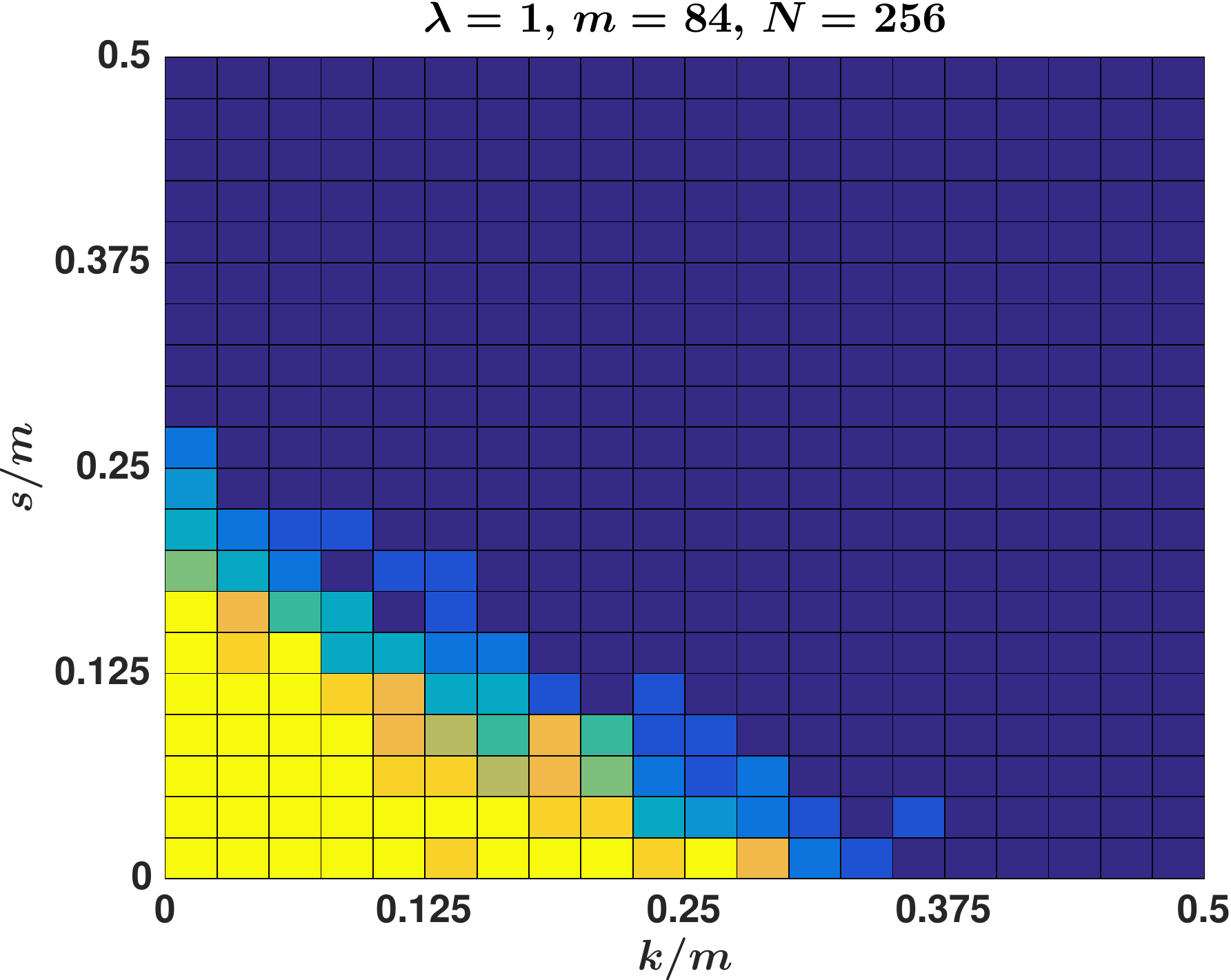}&
      \includegraphics[width=0.33\textwidth]{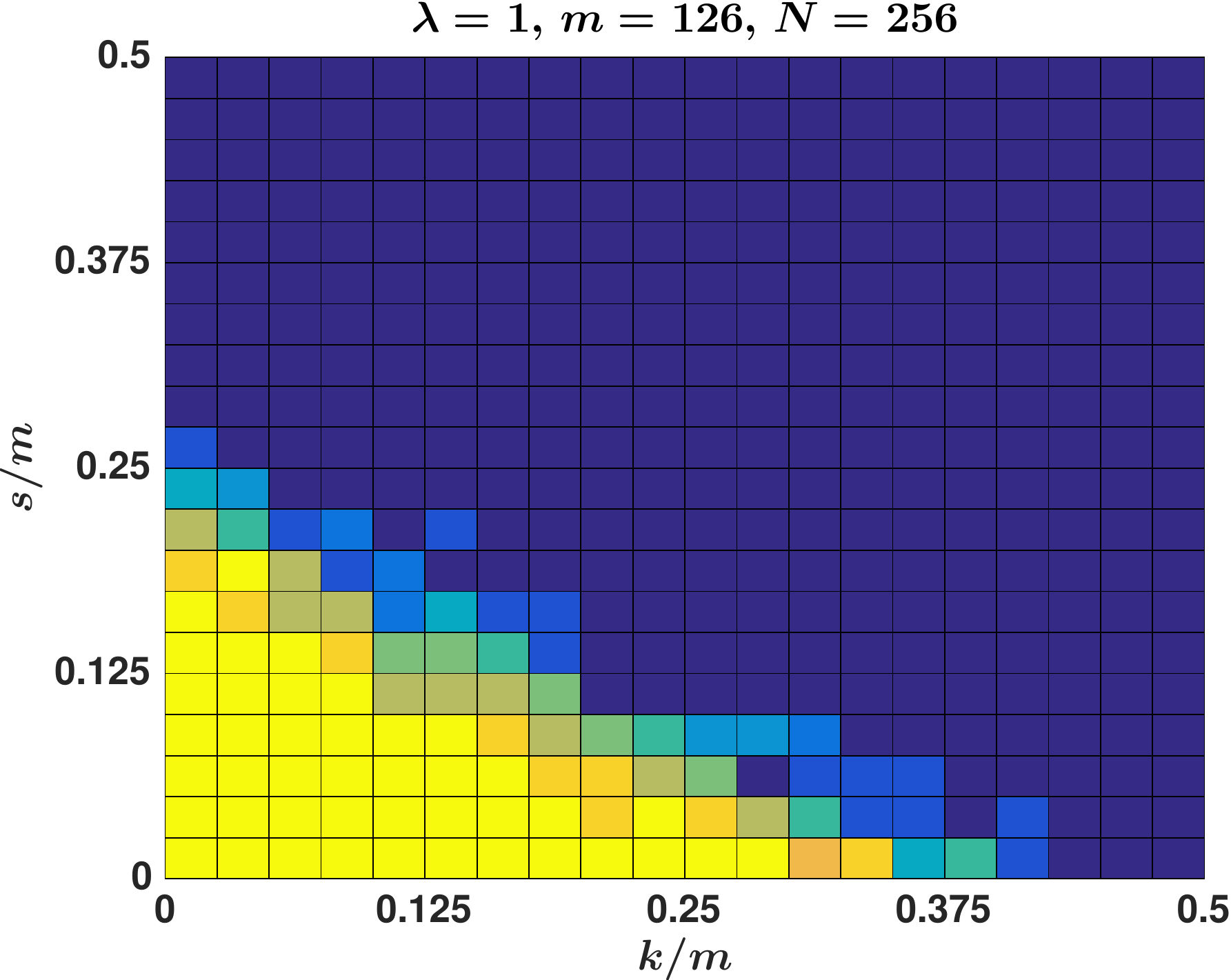}\\
      \includegraphics[width=0.33\textwidth]{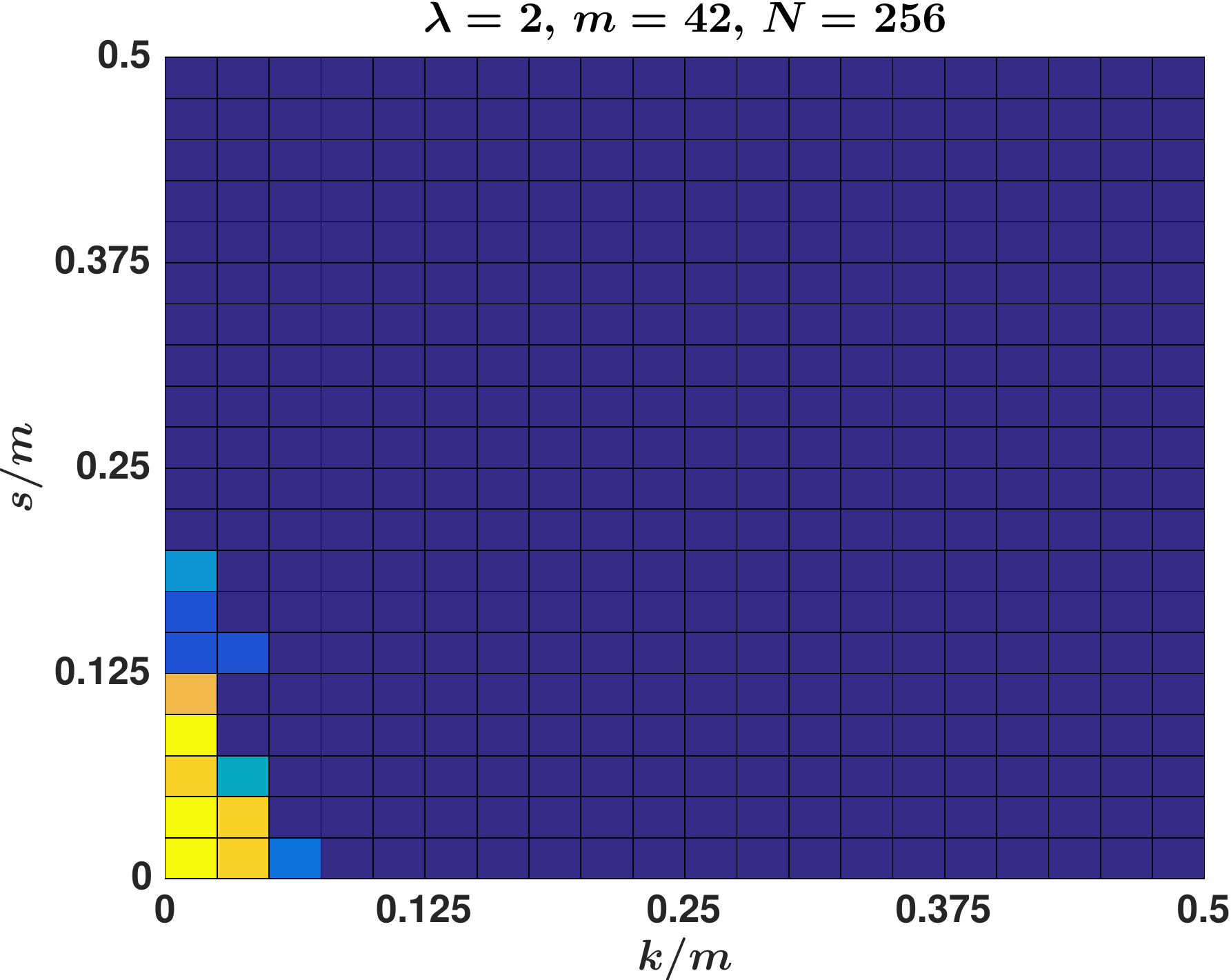}&
      \includegraphics[width=0.33\textwidth]{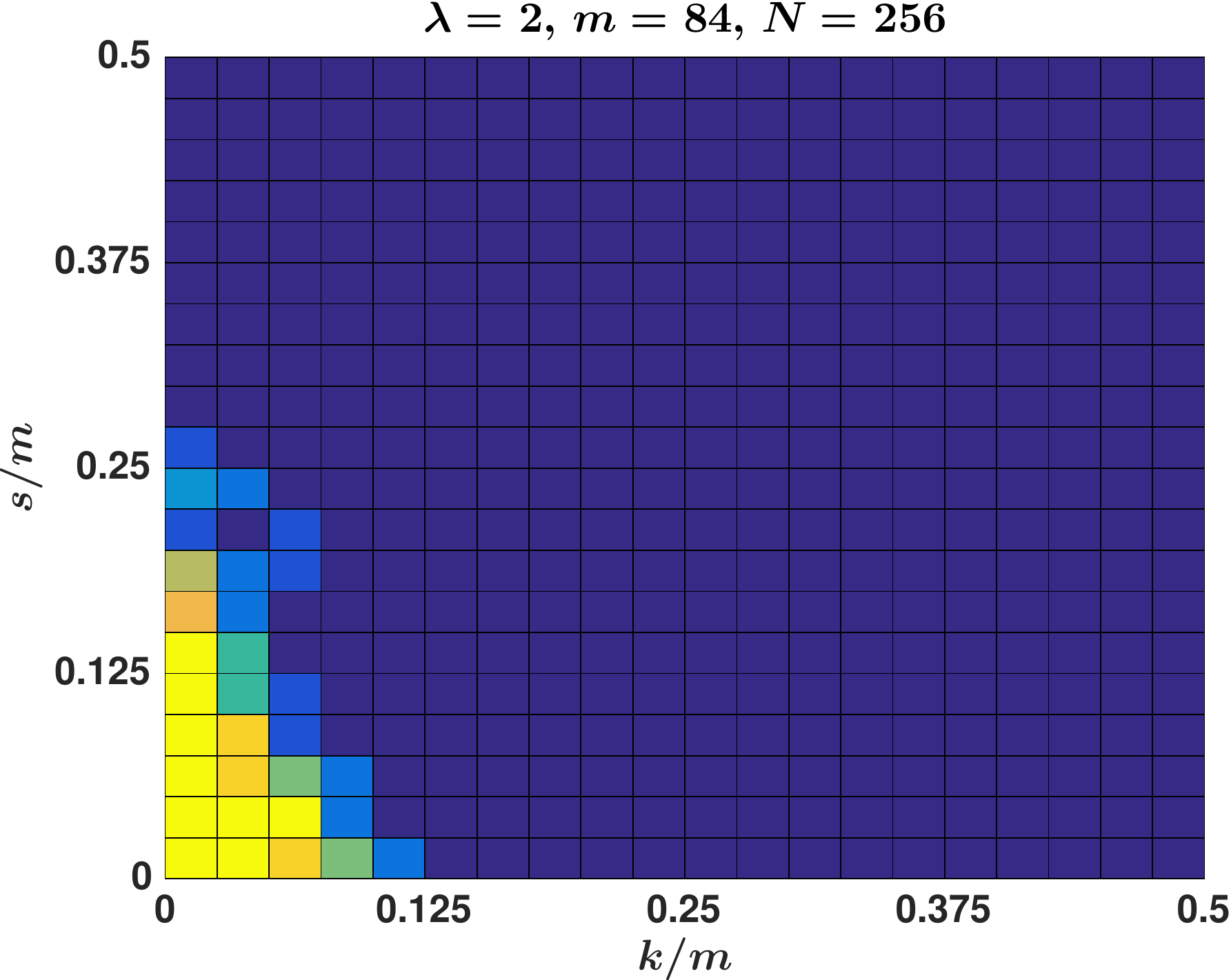}&
      \includegraphics[width=0.33\textwidth]{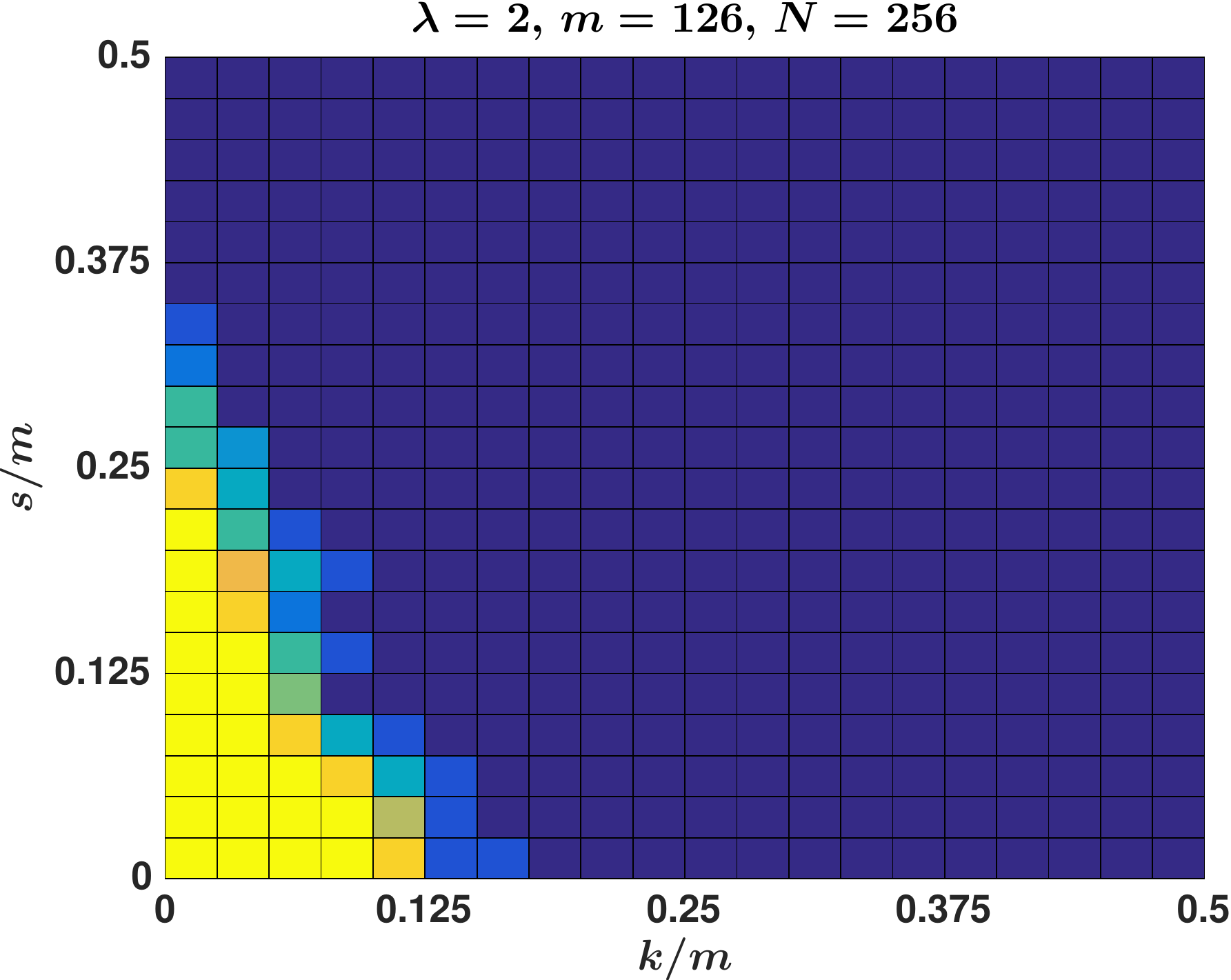}\\
      \includegraphics[width=0.33\textwidth]{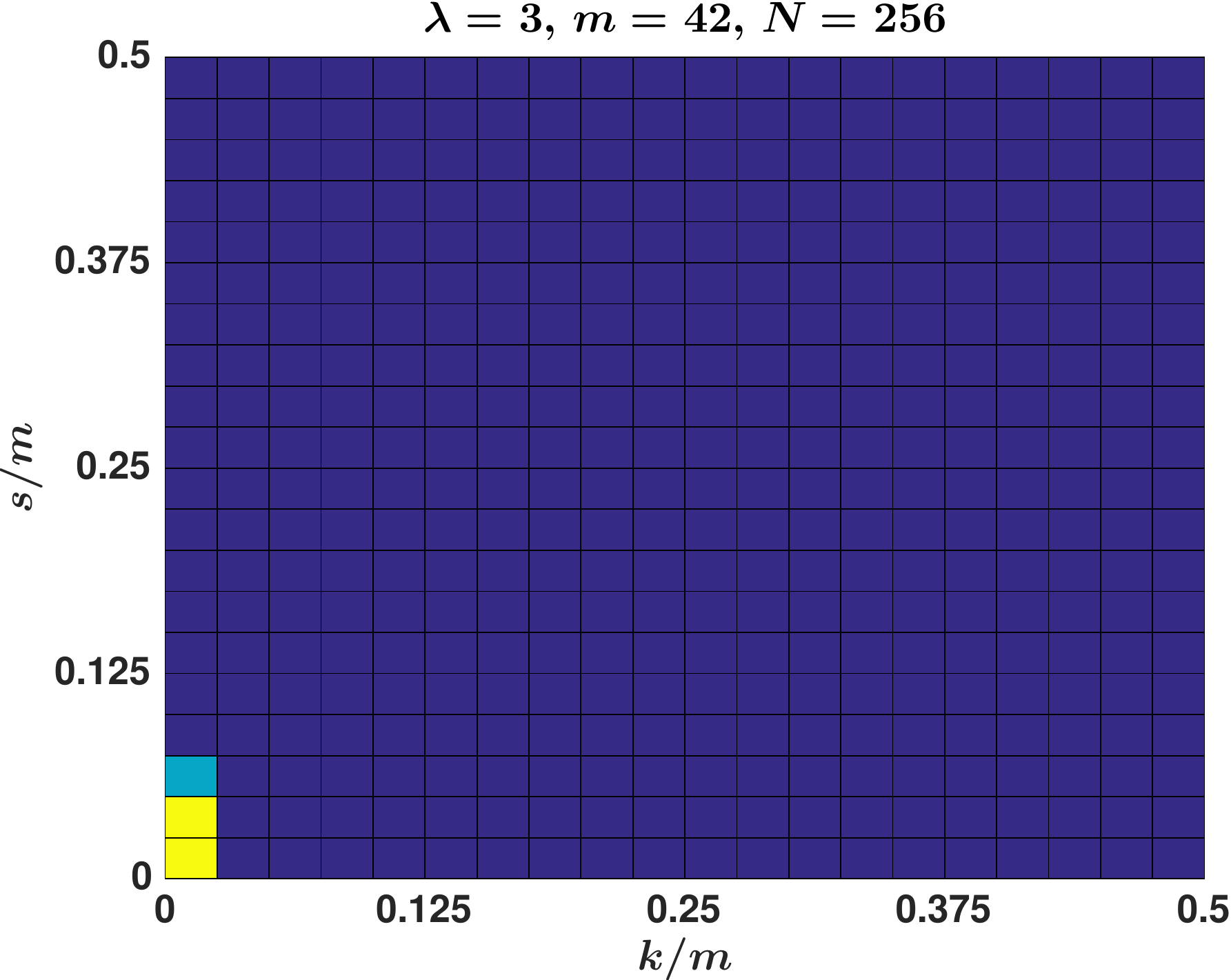}&
      \includegraphics[width=0.33\textwidth]{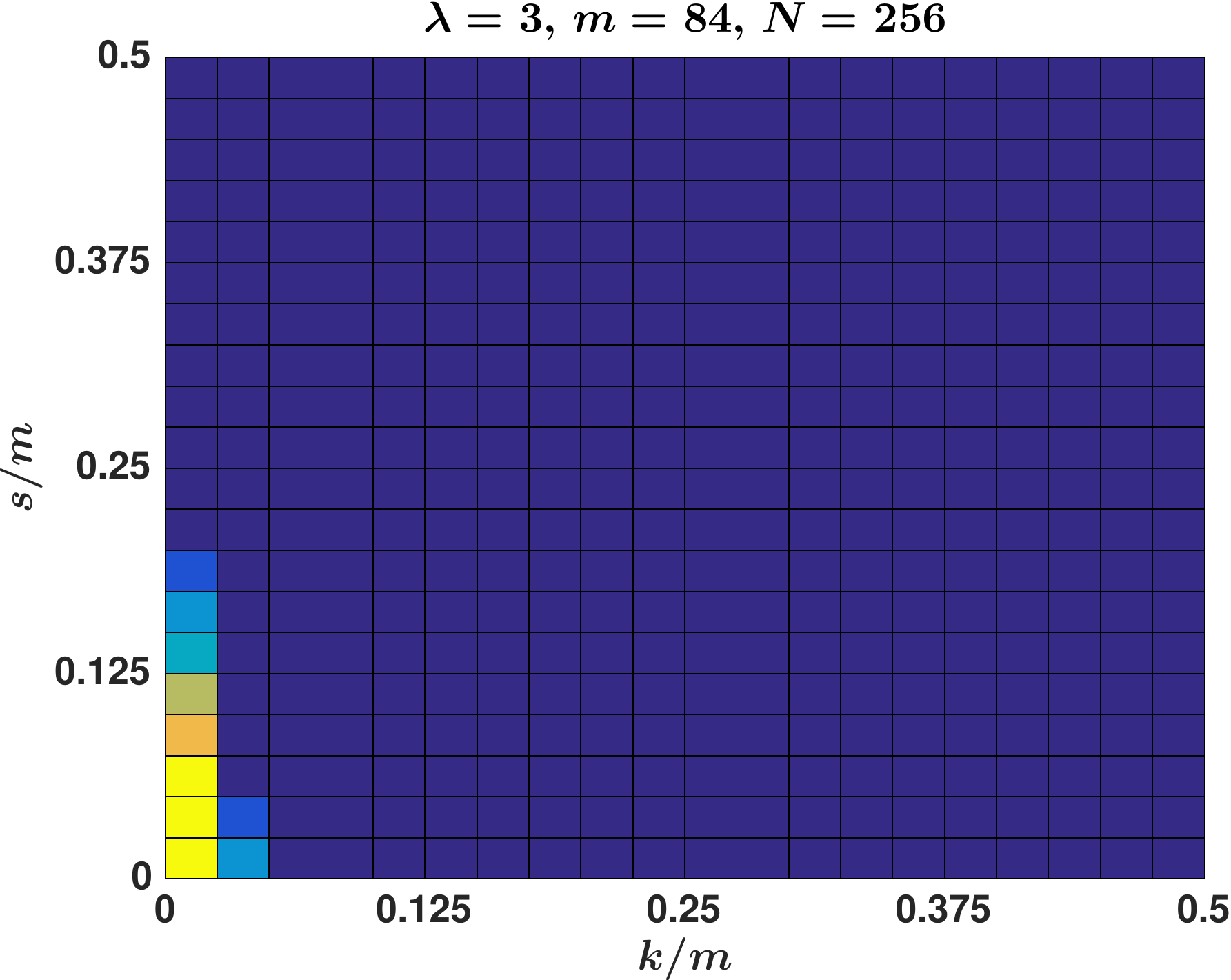}&
      \includegraphics[width=0.33\textwidth]{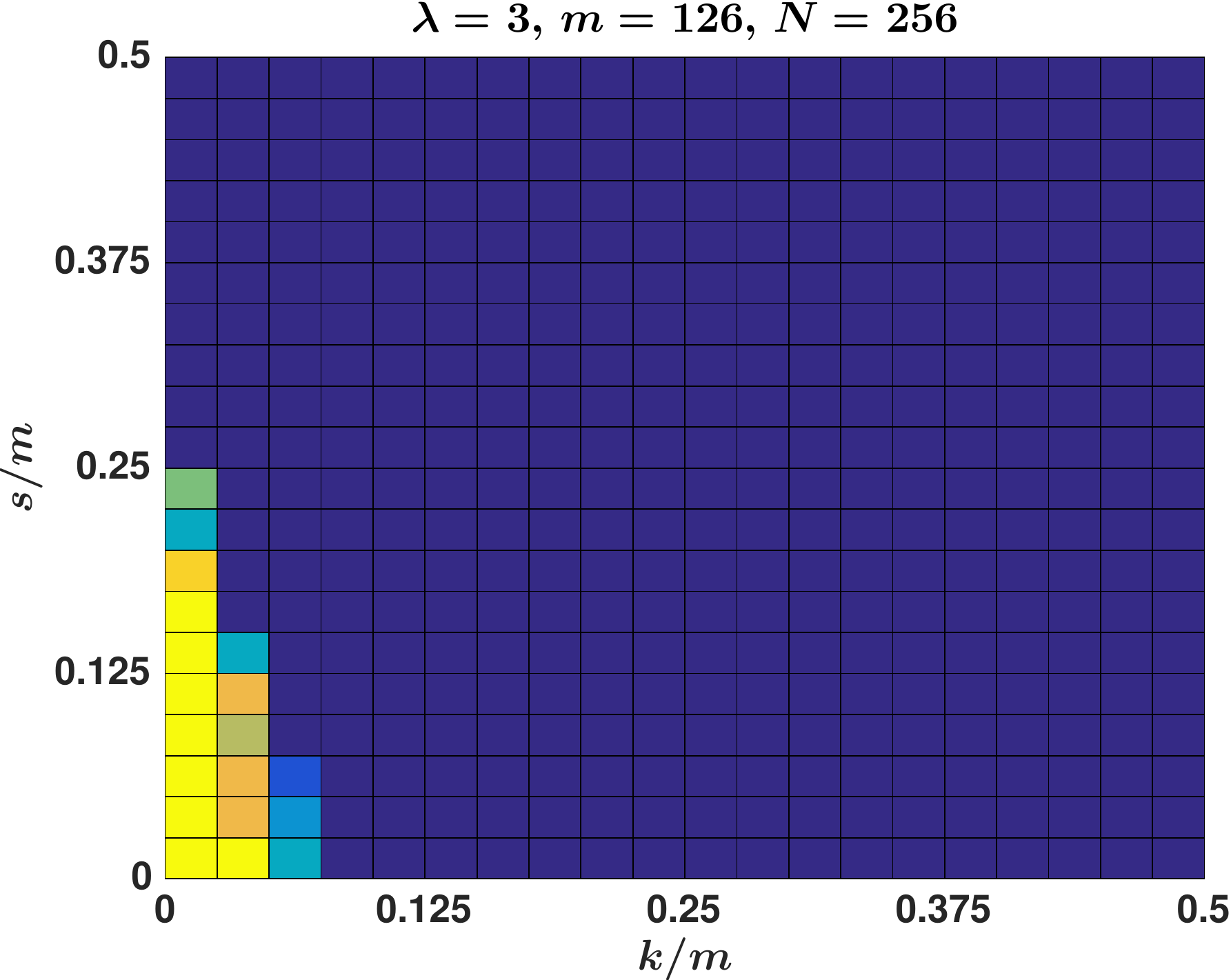}\\
      \includegraphics[width=0.33\textwidth]{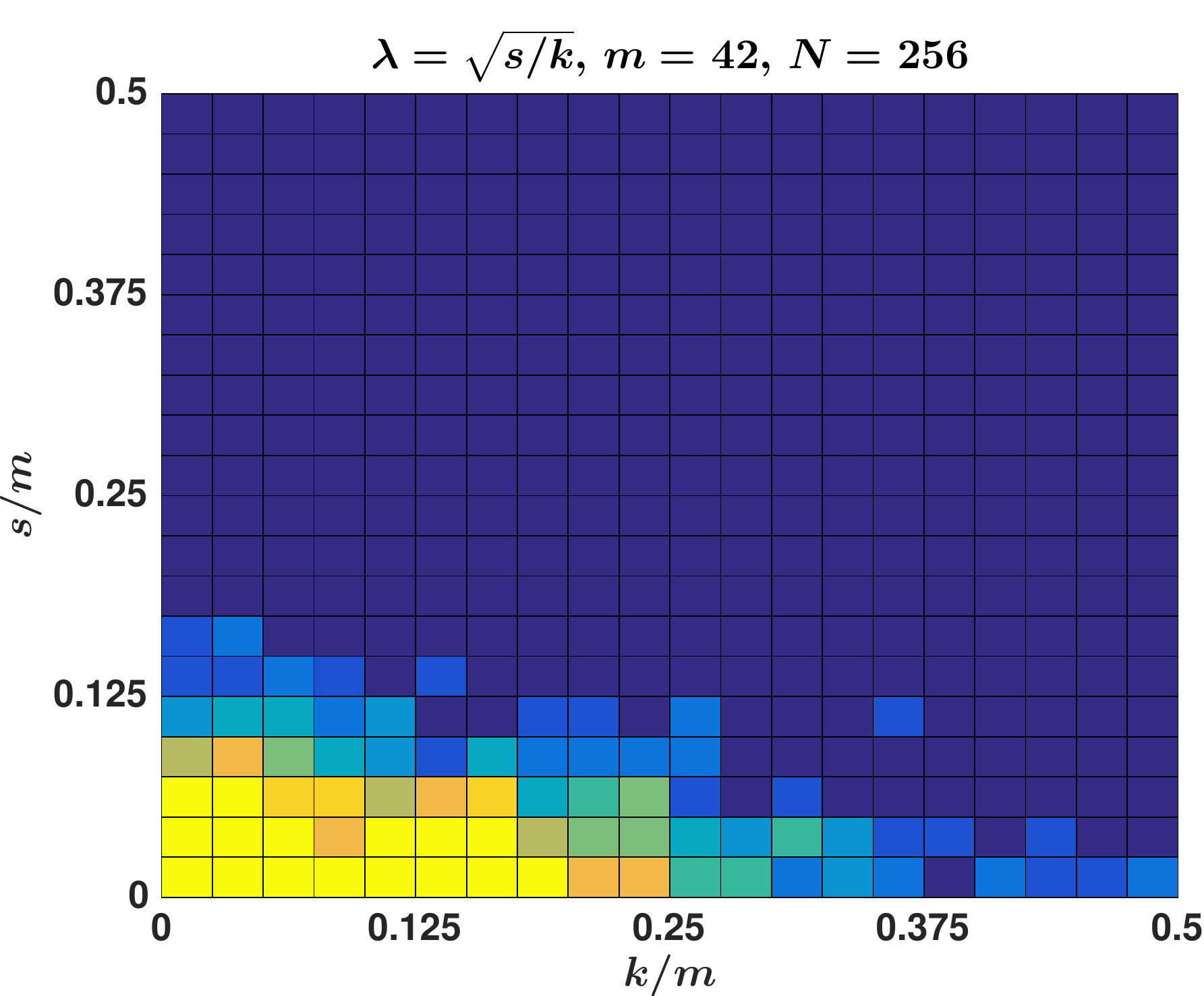}&
      \includegraphics[width=0.33\textwidth]{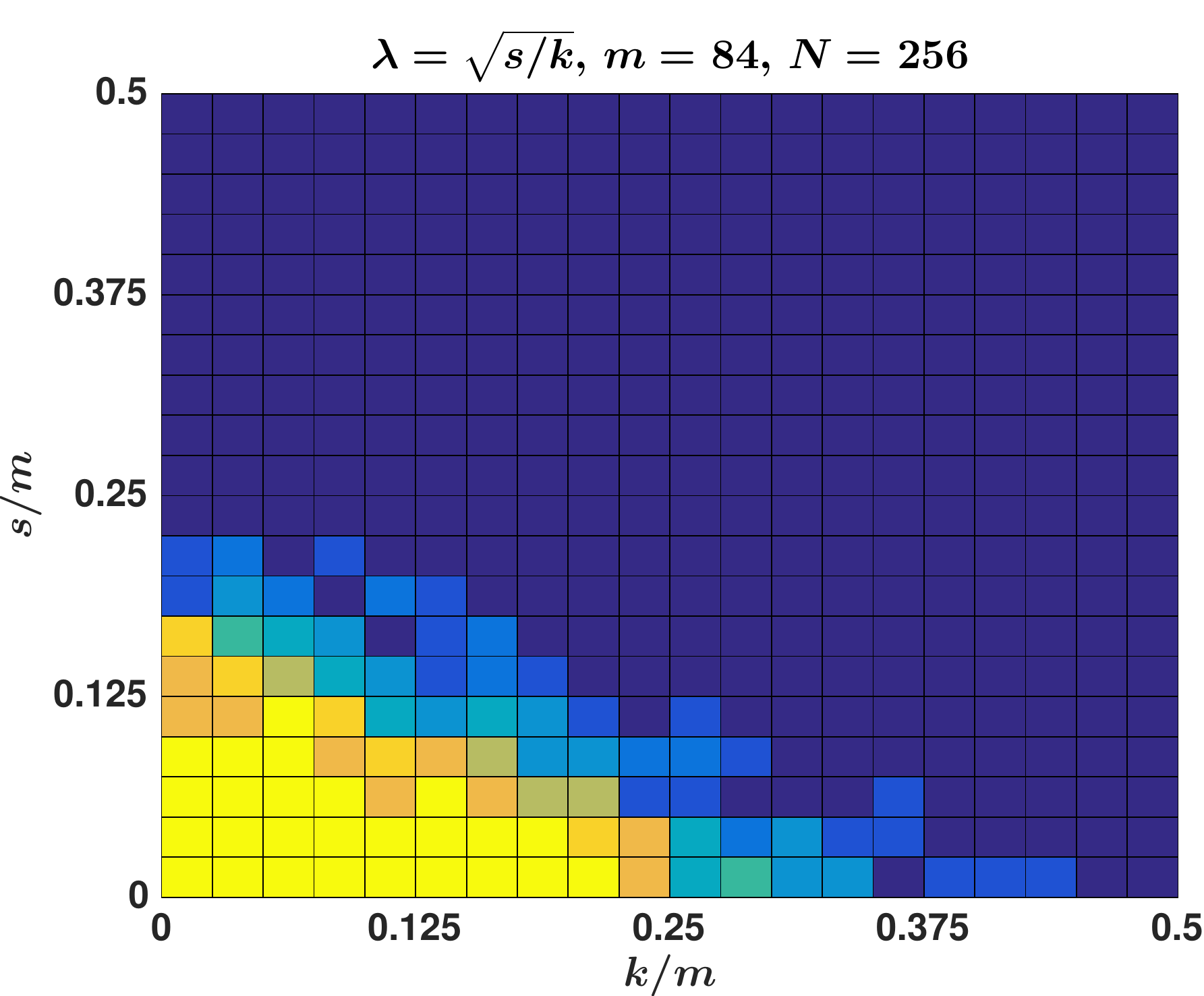}&
      \includegraphics[width=0.33\textwidth]{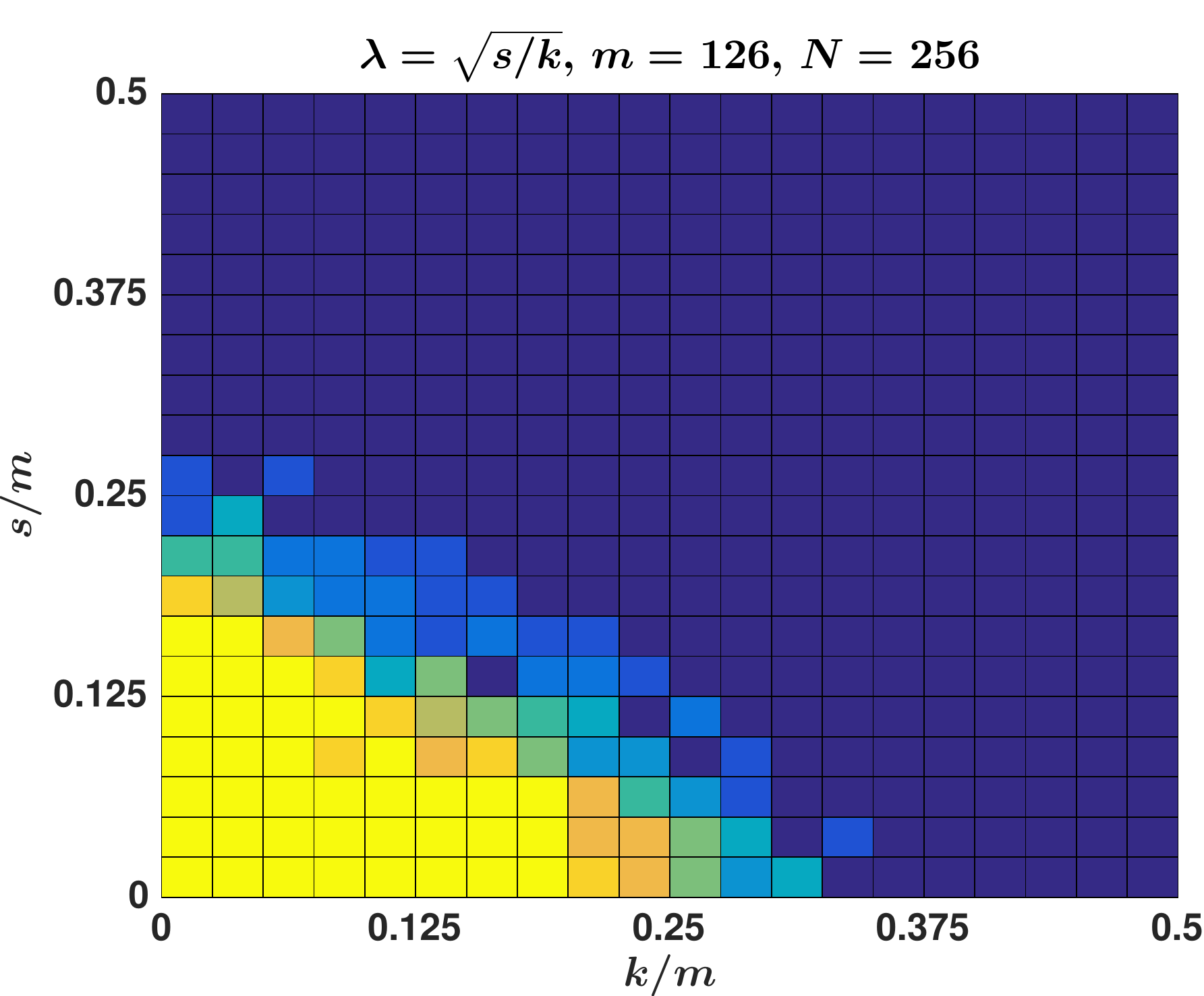}\\
      \includegraphics[width=0.33\textwidth]{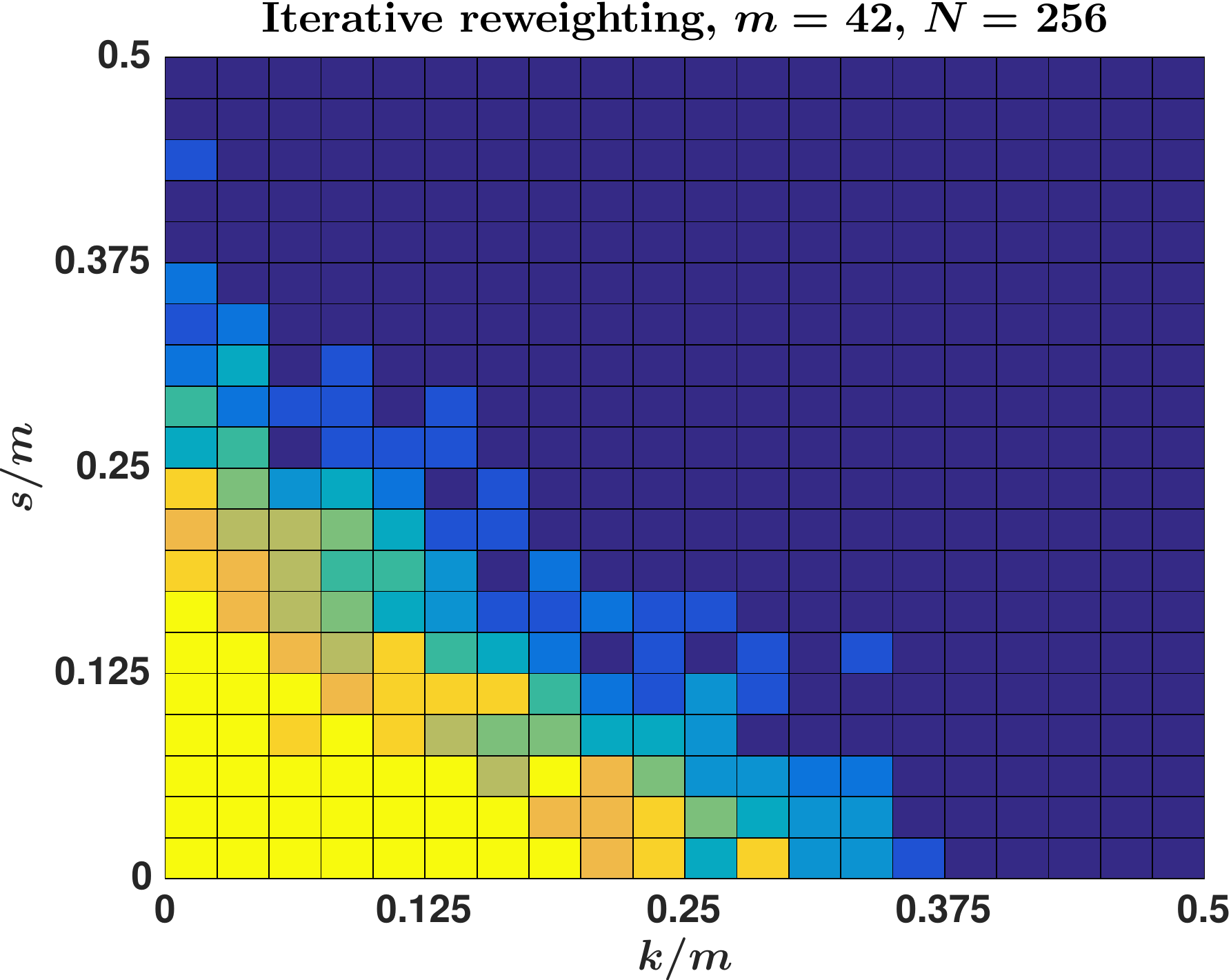}&
      \includegraphics[width=0.33\textwidth]{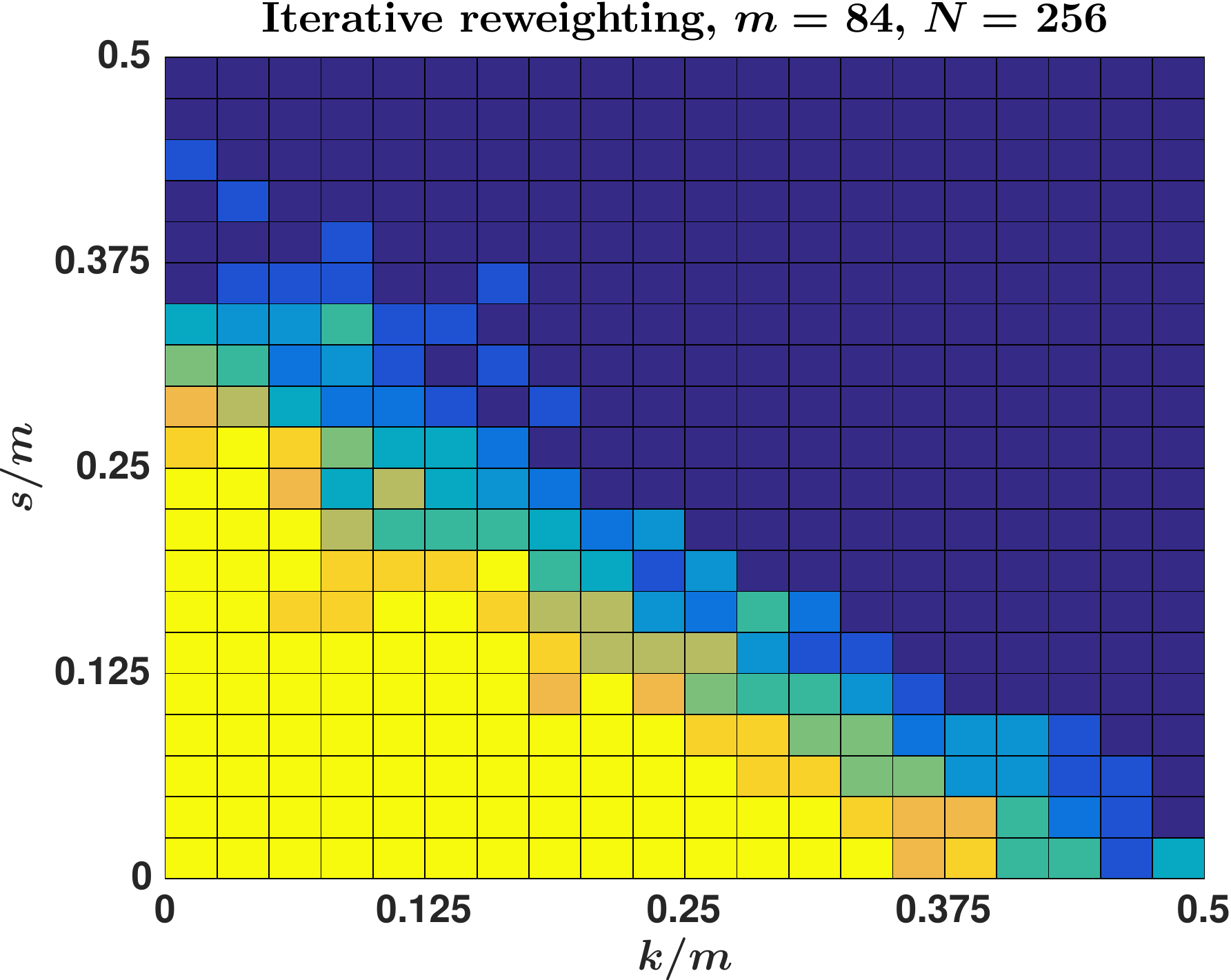}&
      \includegraphics[width=0.33\textwidth]{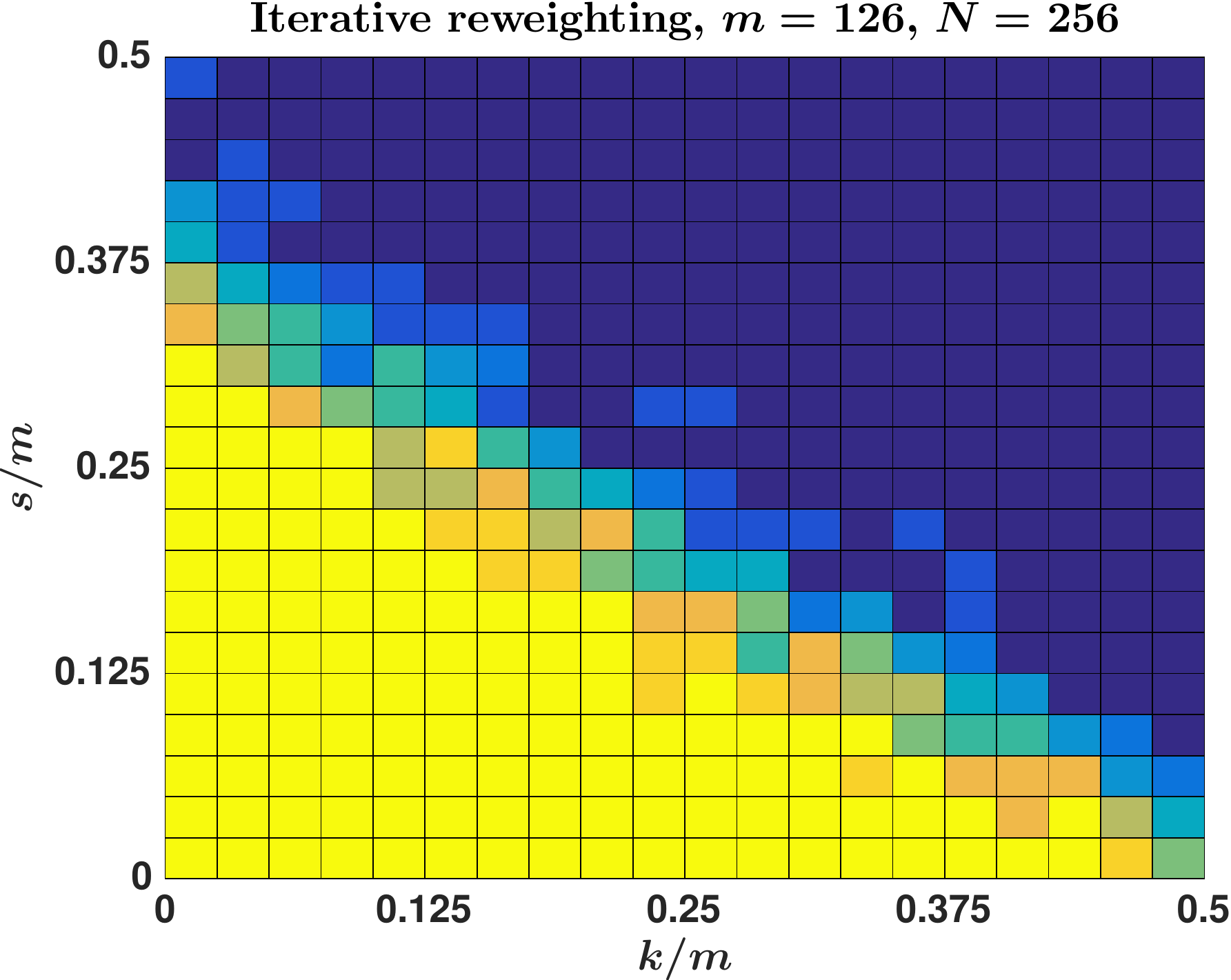}
    \end{tabular}
  }
  \caption{Phase transition for model $1$ with fixed $N = 256$, varying $m$ and $\lambda$. Each column represents varying values of $m$: from left to right, $m=42$, $m=84$, and $m=126$. Each row represents different values of $\lambda$: rows 1-4 correspond to $\lambda = 0.5, 1, 2, 3$, respectively. Row 5 uses the value $\lambda = \sqrt{s/k}$ that is suggested as optimal by the theory. Row 6 shows recovery using the iteratively reweighted $\ell^1$ algorithm. Each pixel is colored according to its probability of a successful signal recovery for $T=10$ trials based on repeated random draws of $x$ and $c$; yellow is probability 1, blue is probability 0. 
  } \label{f:model1}
\end{center}
\end{figure}

\begin{figure}
\begin{center}
  \resizebox{!}{0.43\textheight}{
    \begin{tabular}{ccc}
      \includegraphics[width=0.33\textwidth]{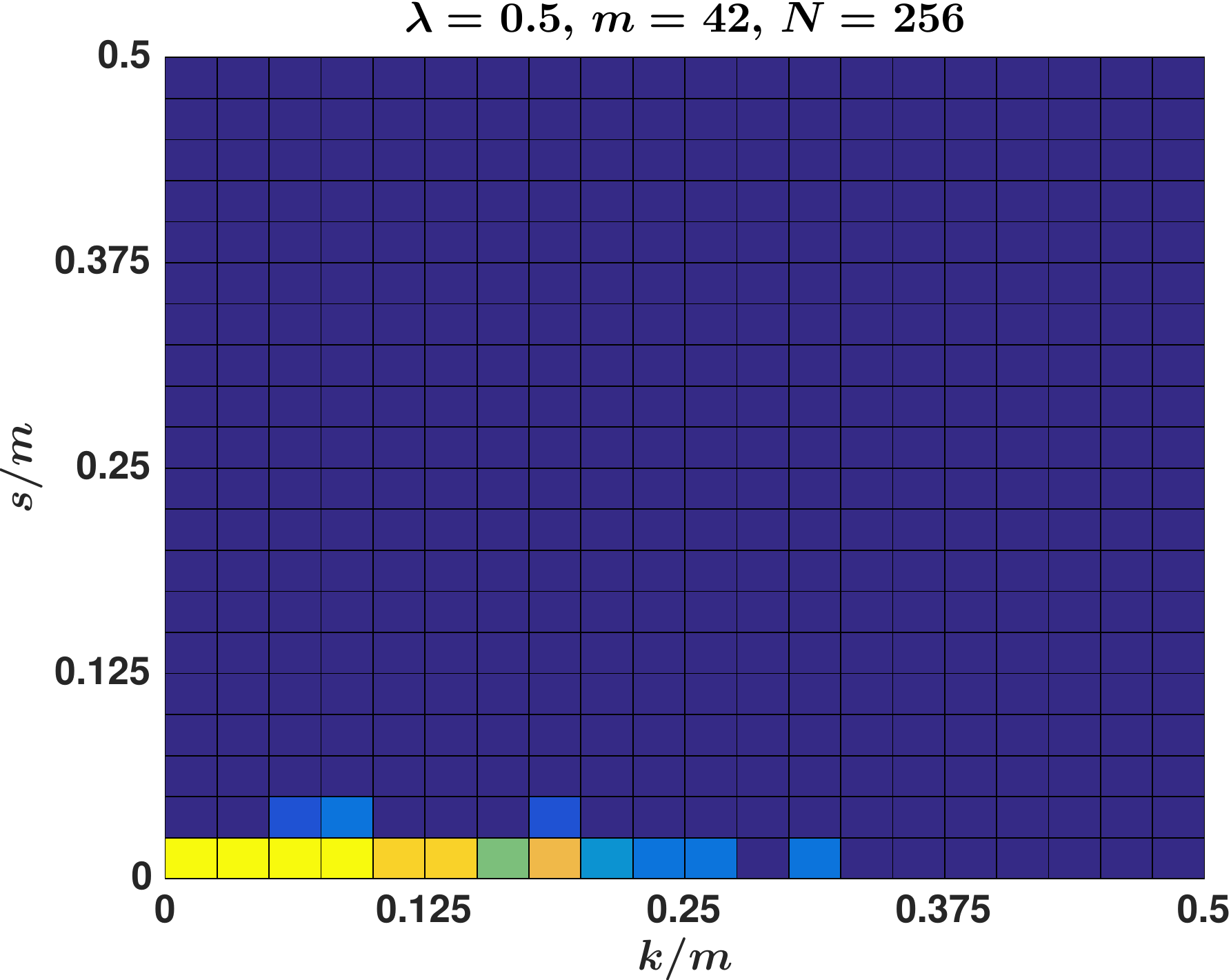}&
      \includegraphics[width=0.33\textwidth]{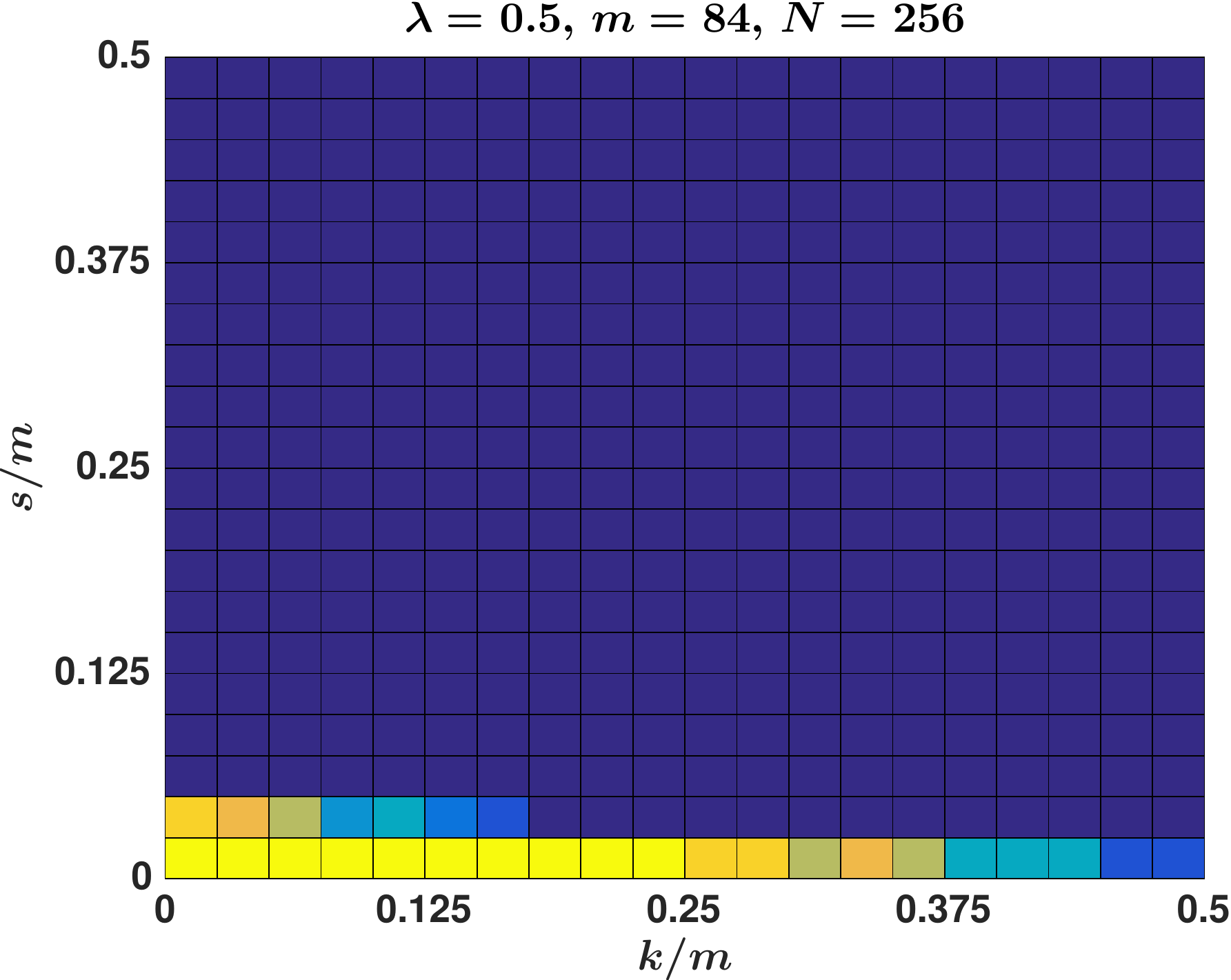}&
      \includegraphics[width=0.33\textwidth]{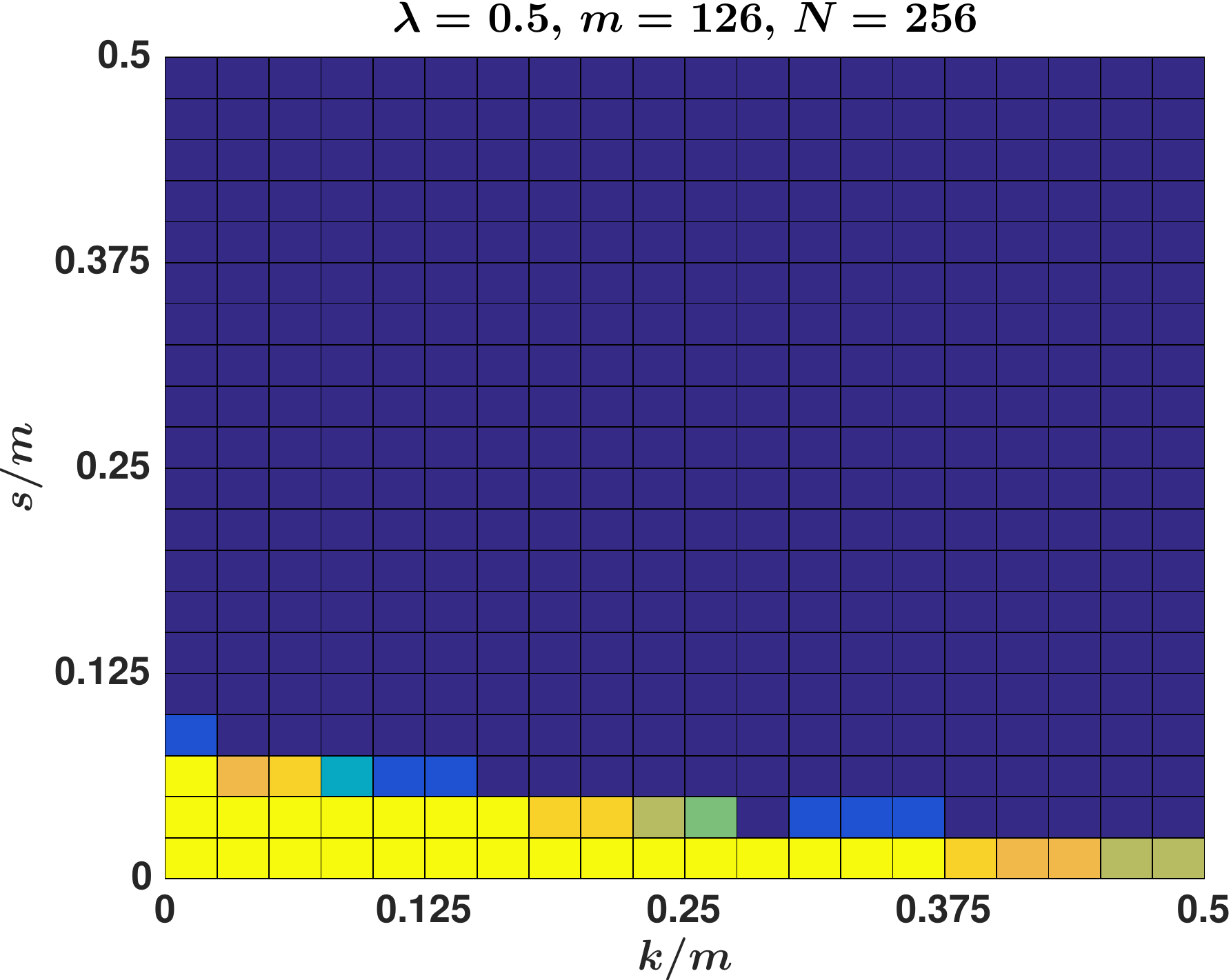}\\
      \includegraphics[width=0.33\textwidth]{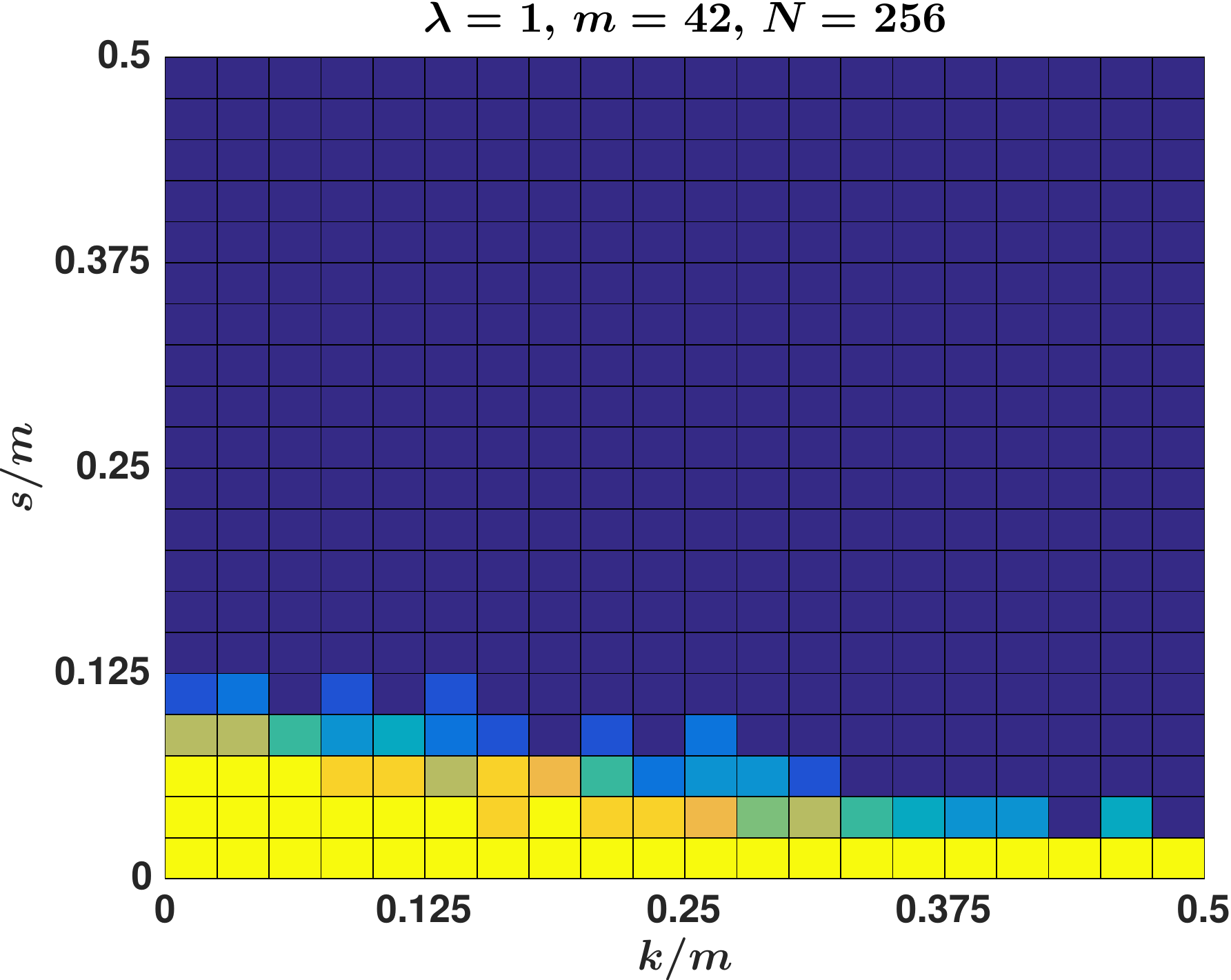}&
      \includegraphics[width=0.33\textwidth]{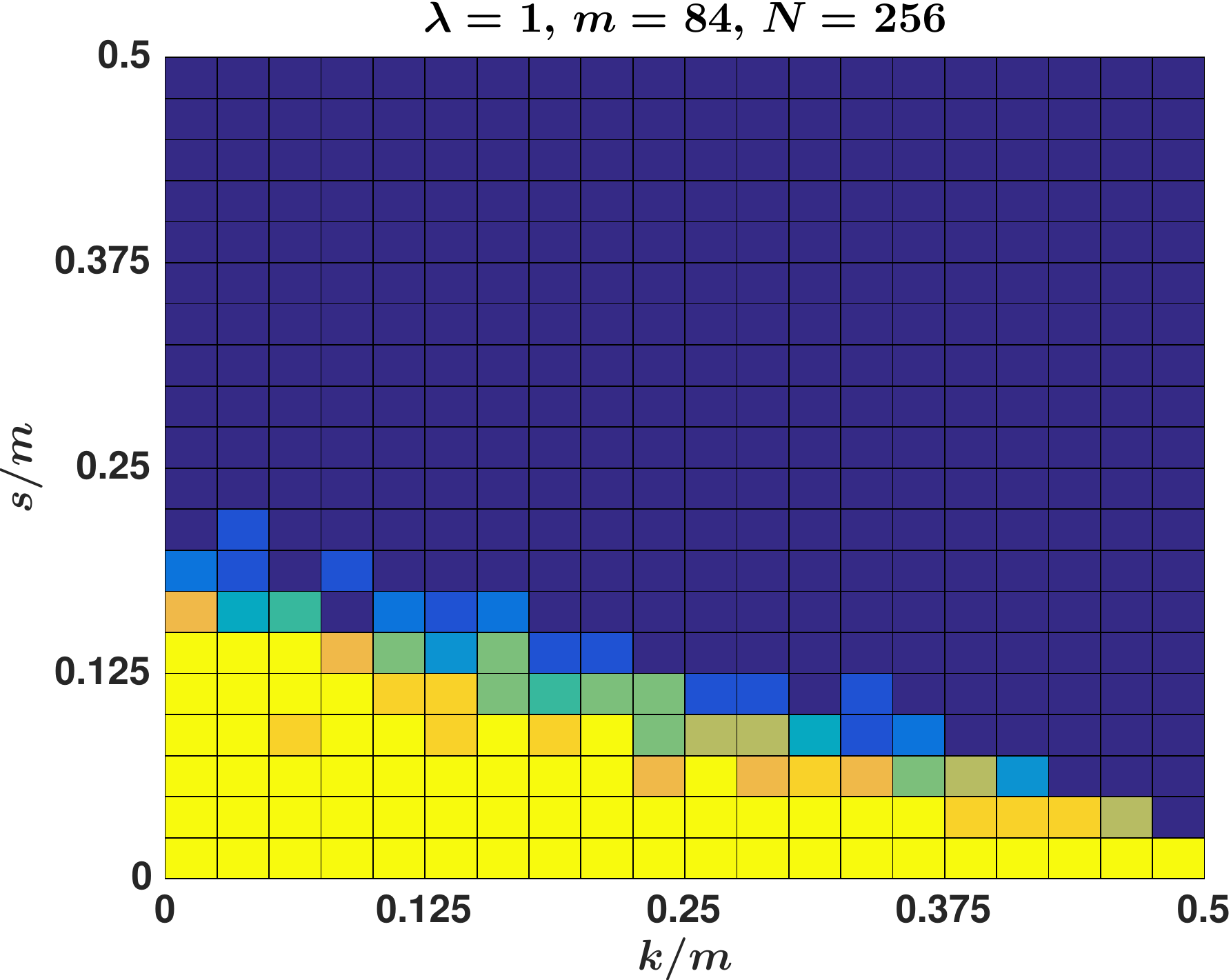}&
      \includegraphics[width=0.33\textwidth]{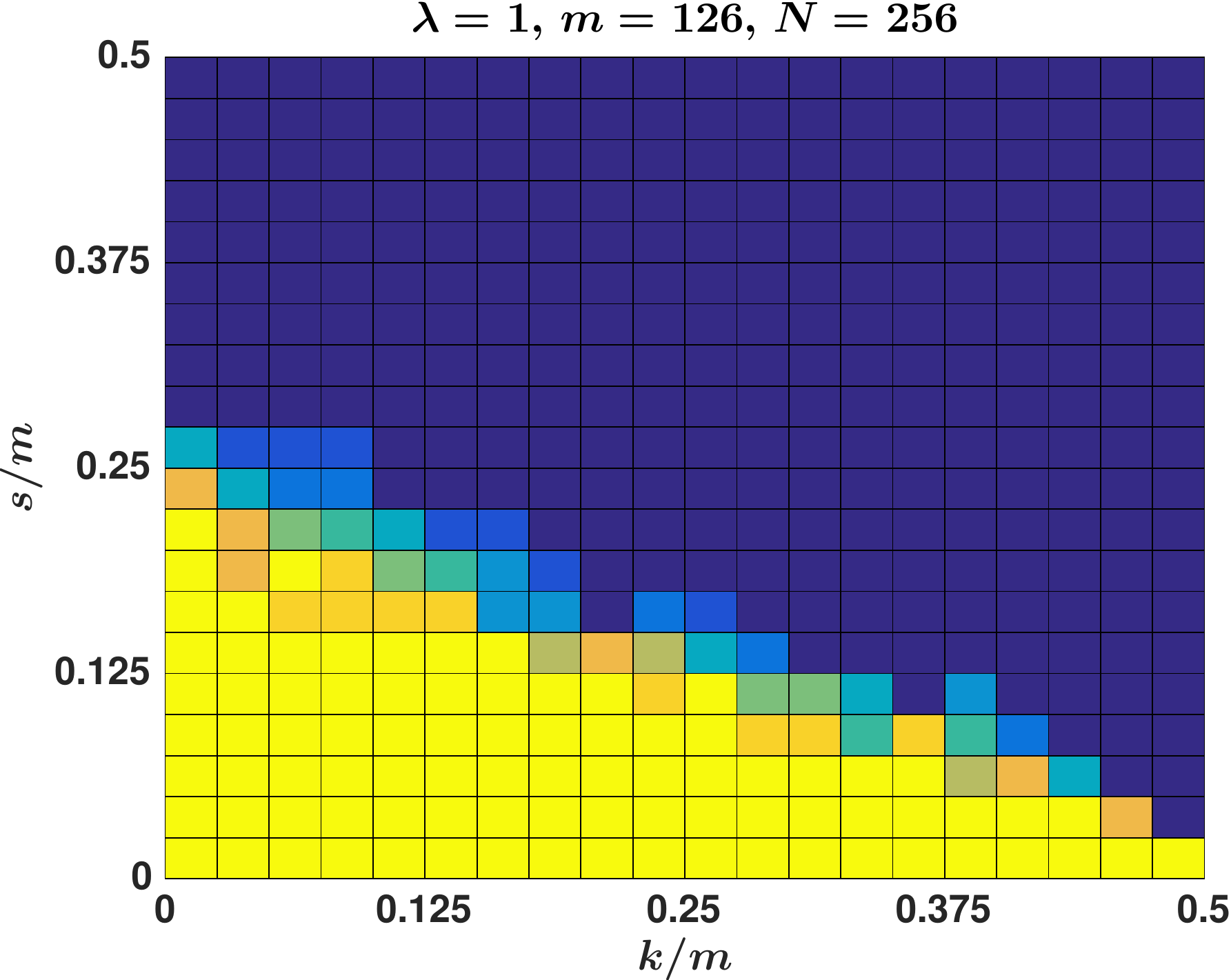}\\
      \includegraphics[width=0.33\textwidth]{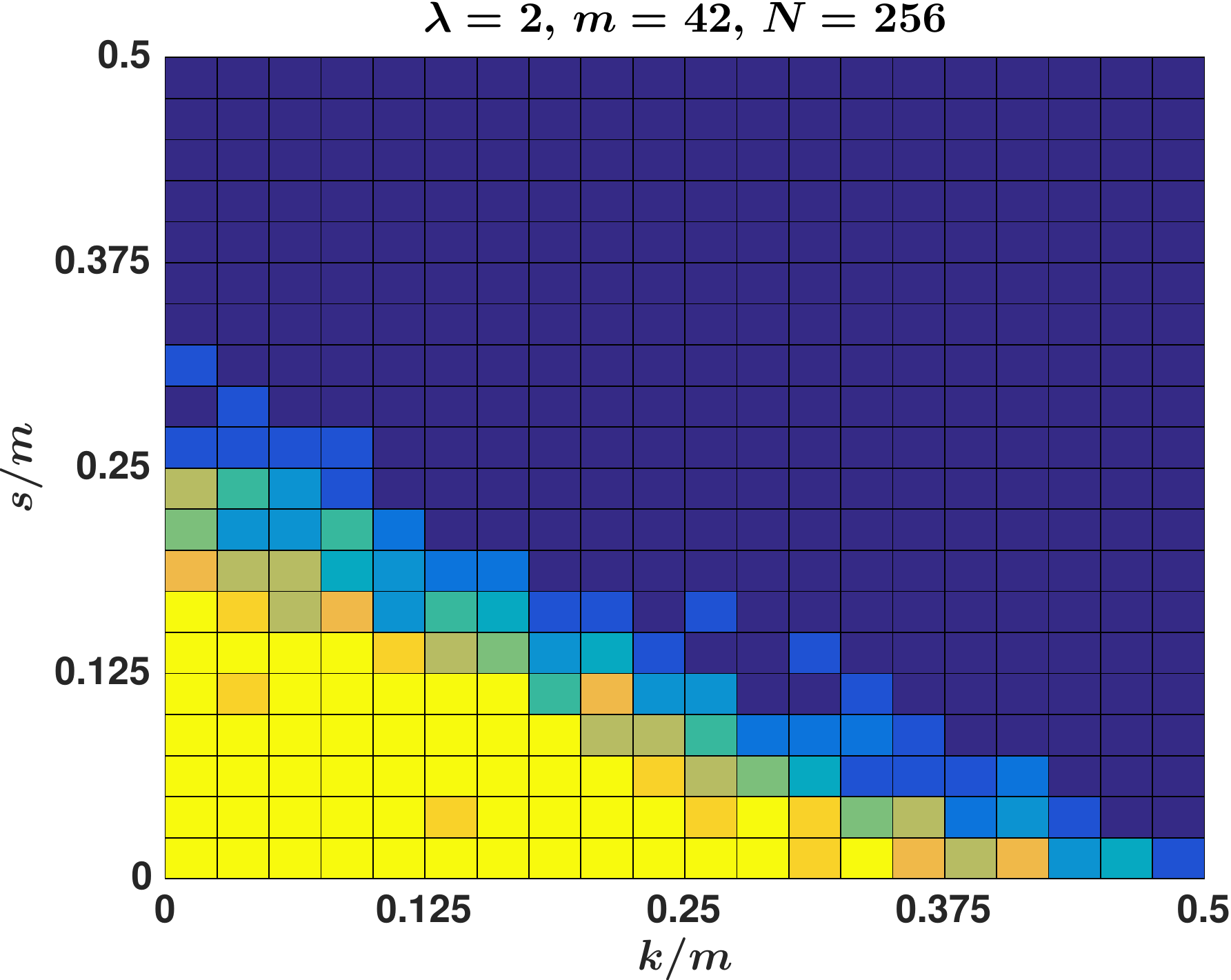}&
      \includegraphics[width=0.33\textwidth]{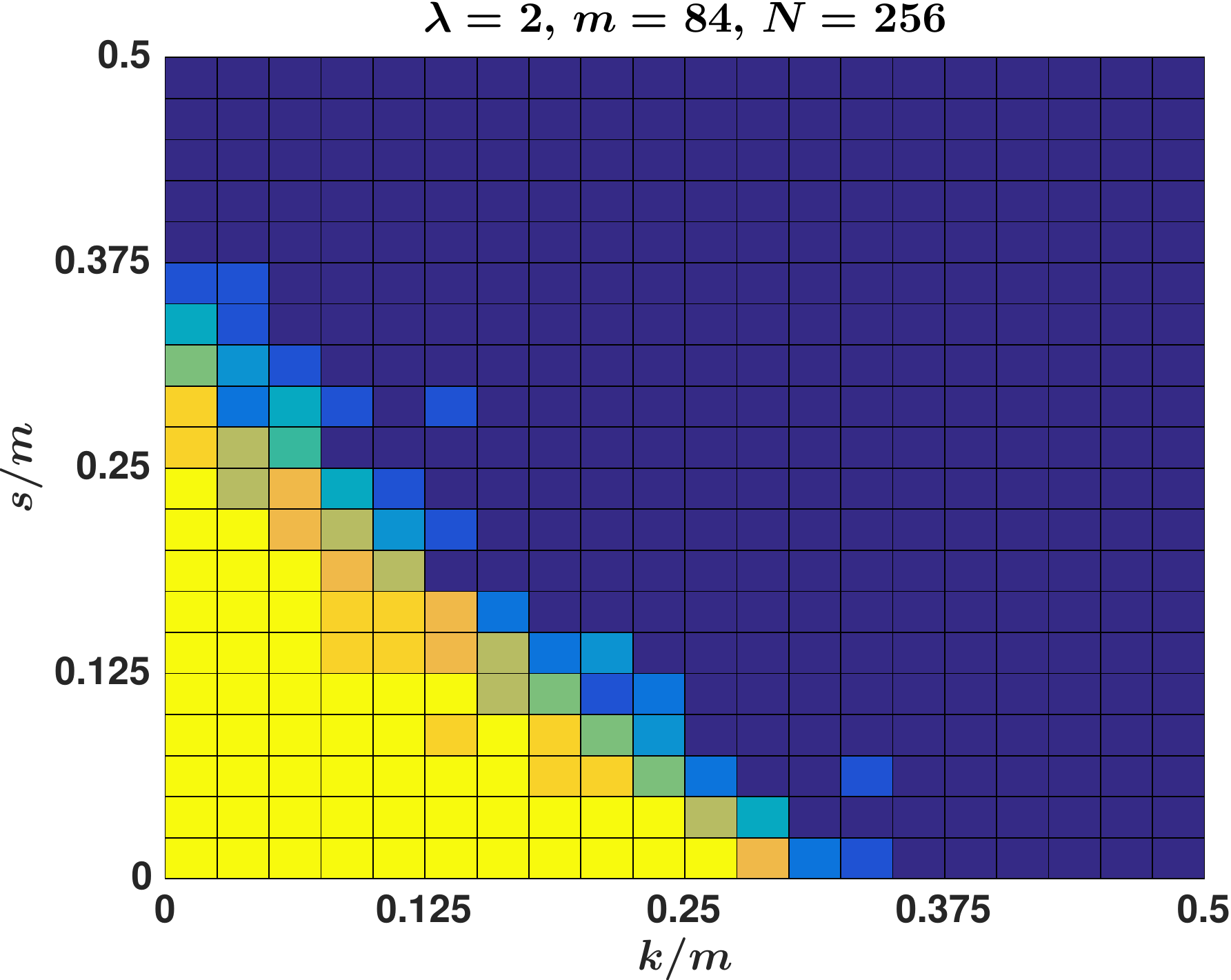}&
      \includegraphics[width=0.33\textwidth]{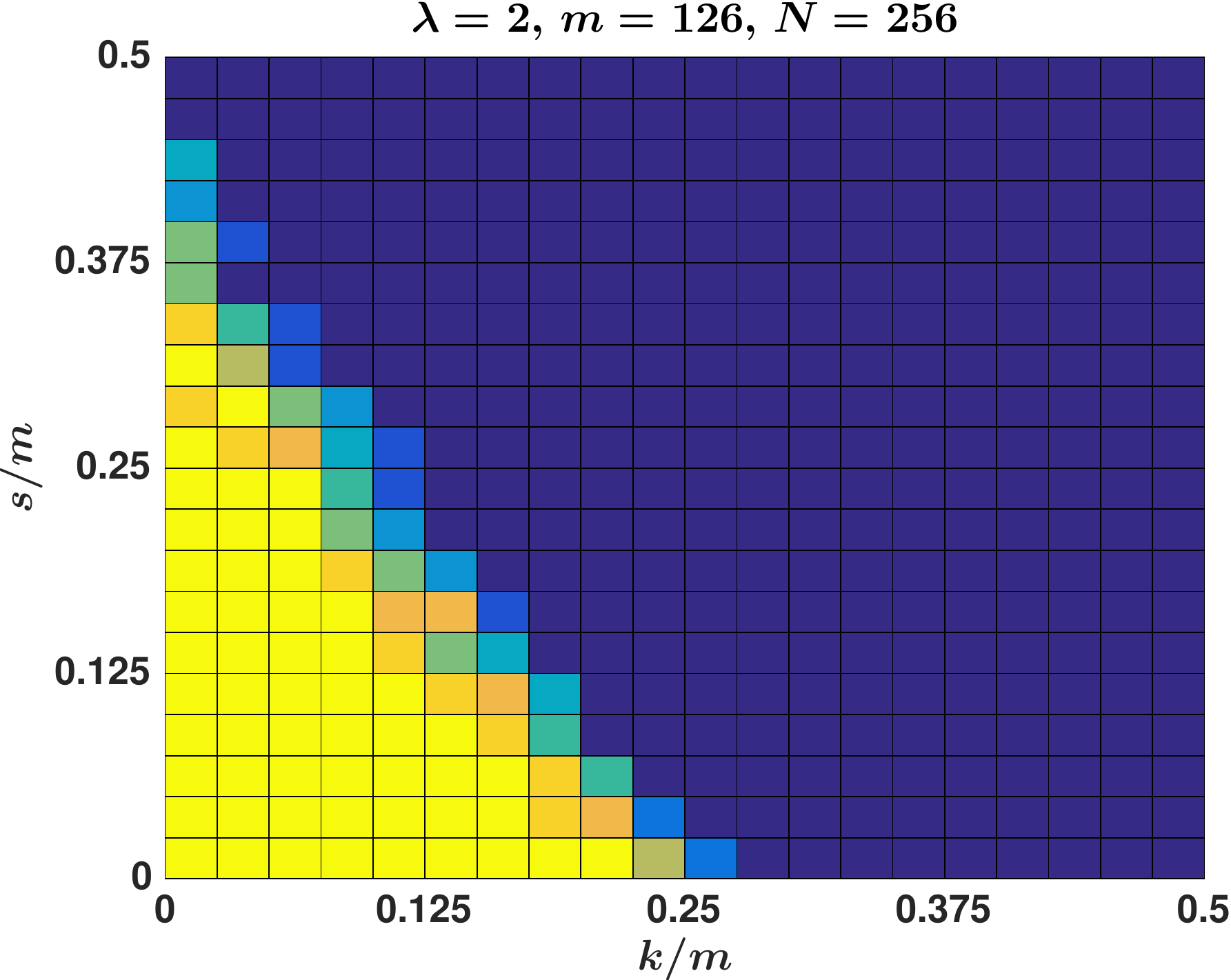}\\
      \includegraphics[width=0.33\textwidth]{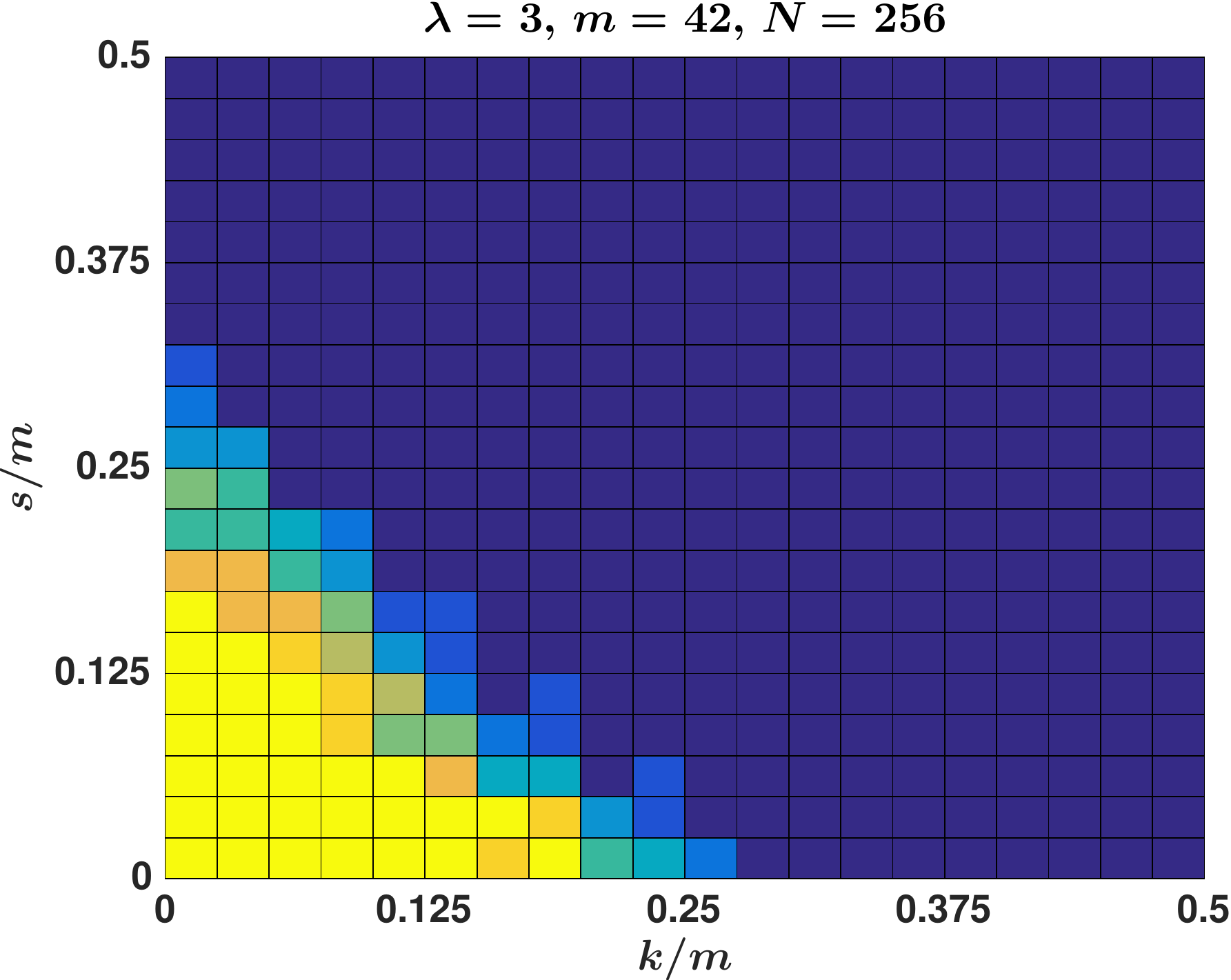}&
      \includegraphics[width=0.33\textwidth]{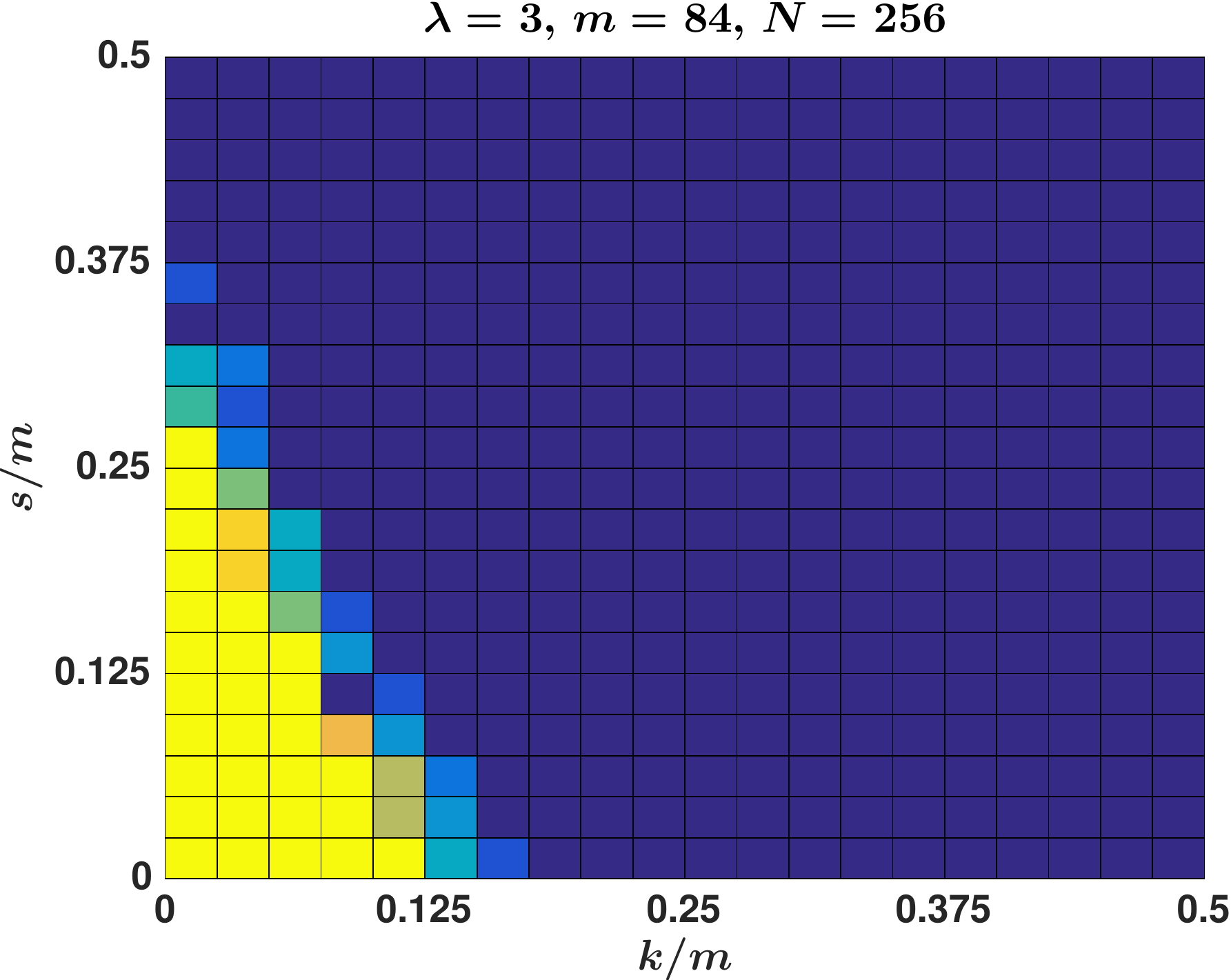}&
      \includegraphics[width=0.33\textwidth]{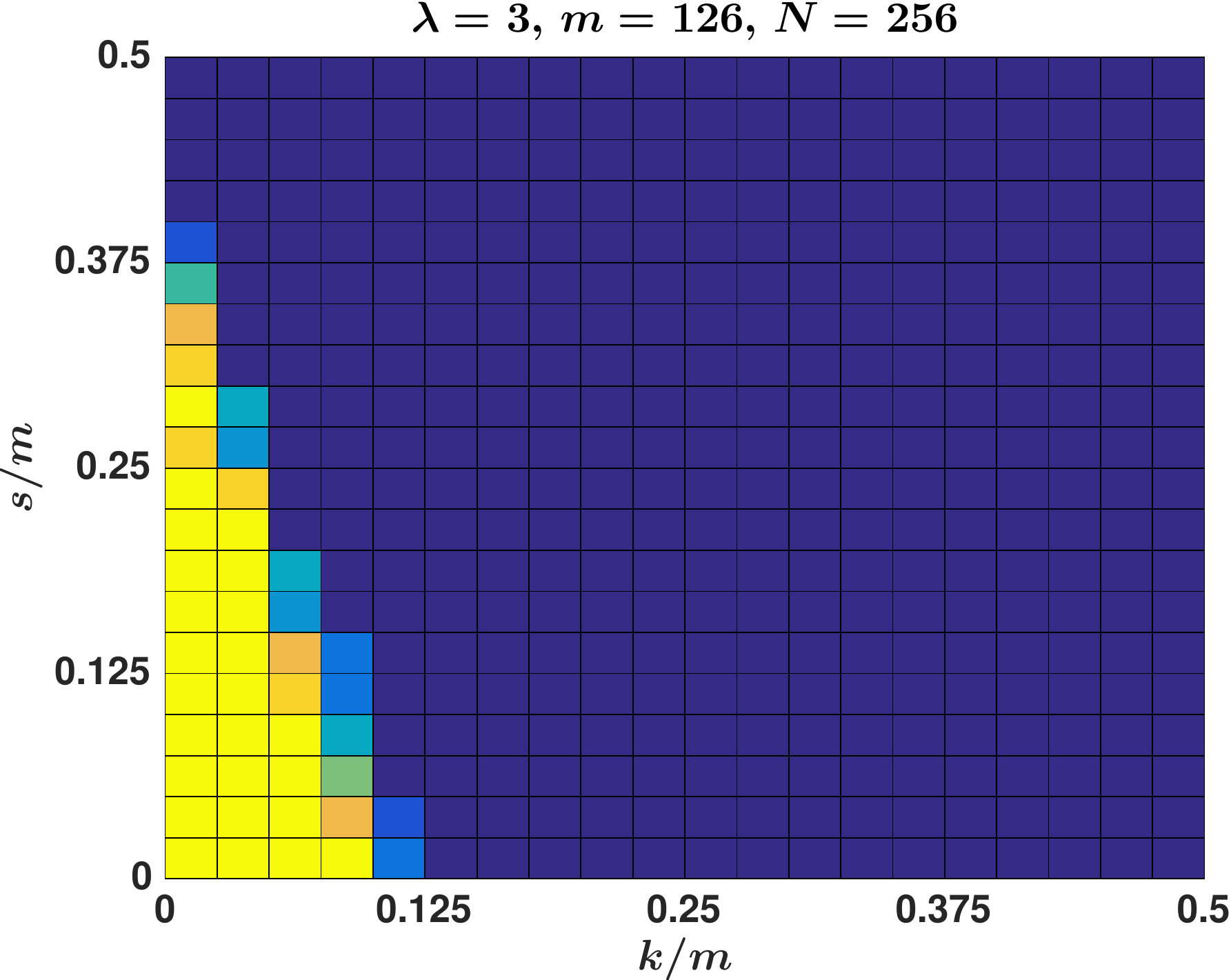}\\
      \includegraphics[width=0.33\textwidth]{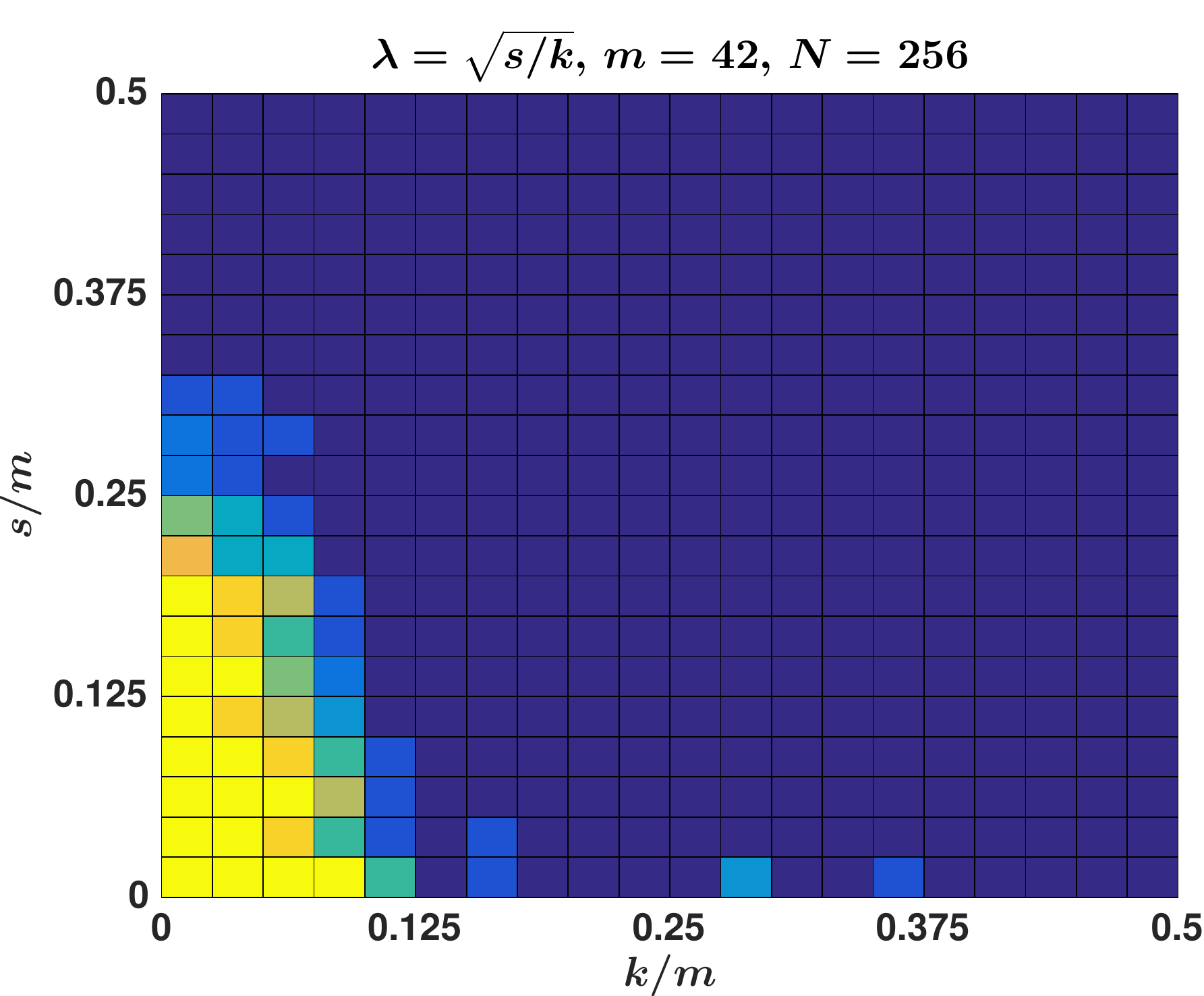}&
      \includegraphics[width=0.33\textwidth]{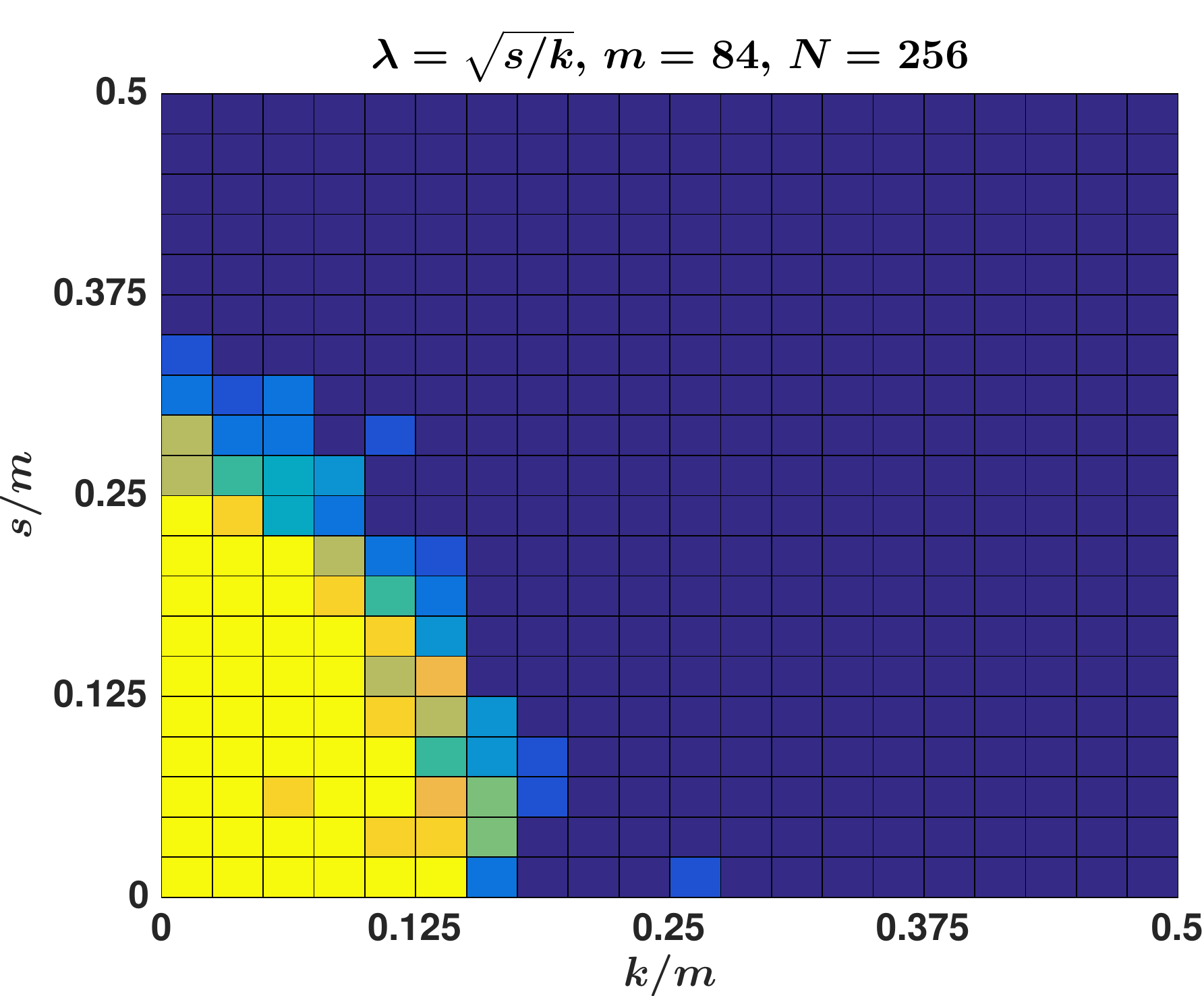}&
      \includegraphics[width=0.33\textwidth]{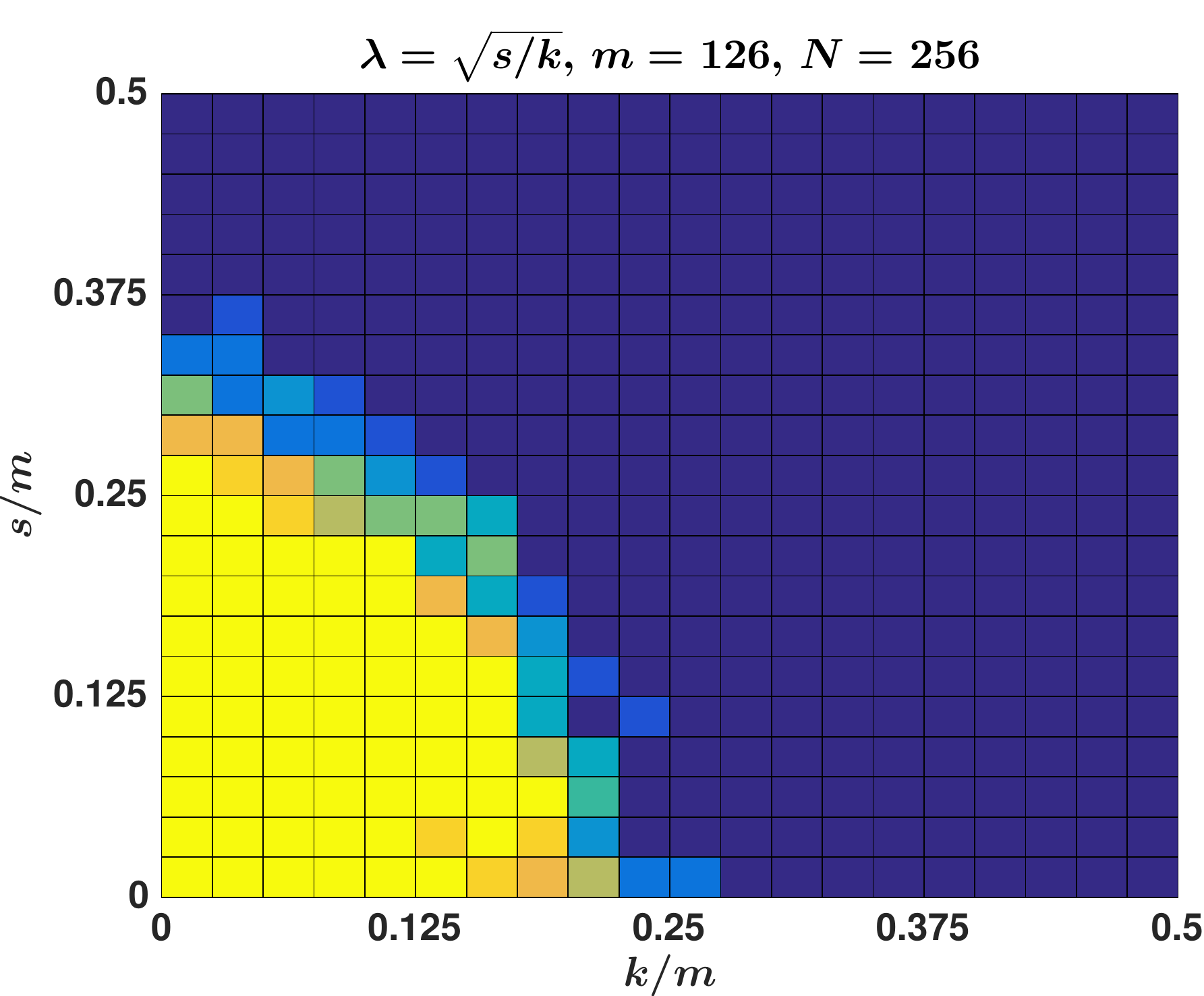}\\
      \includegraphics[width=0.33\textwidth]{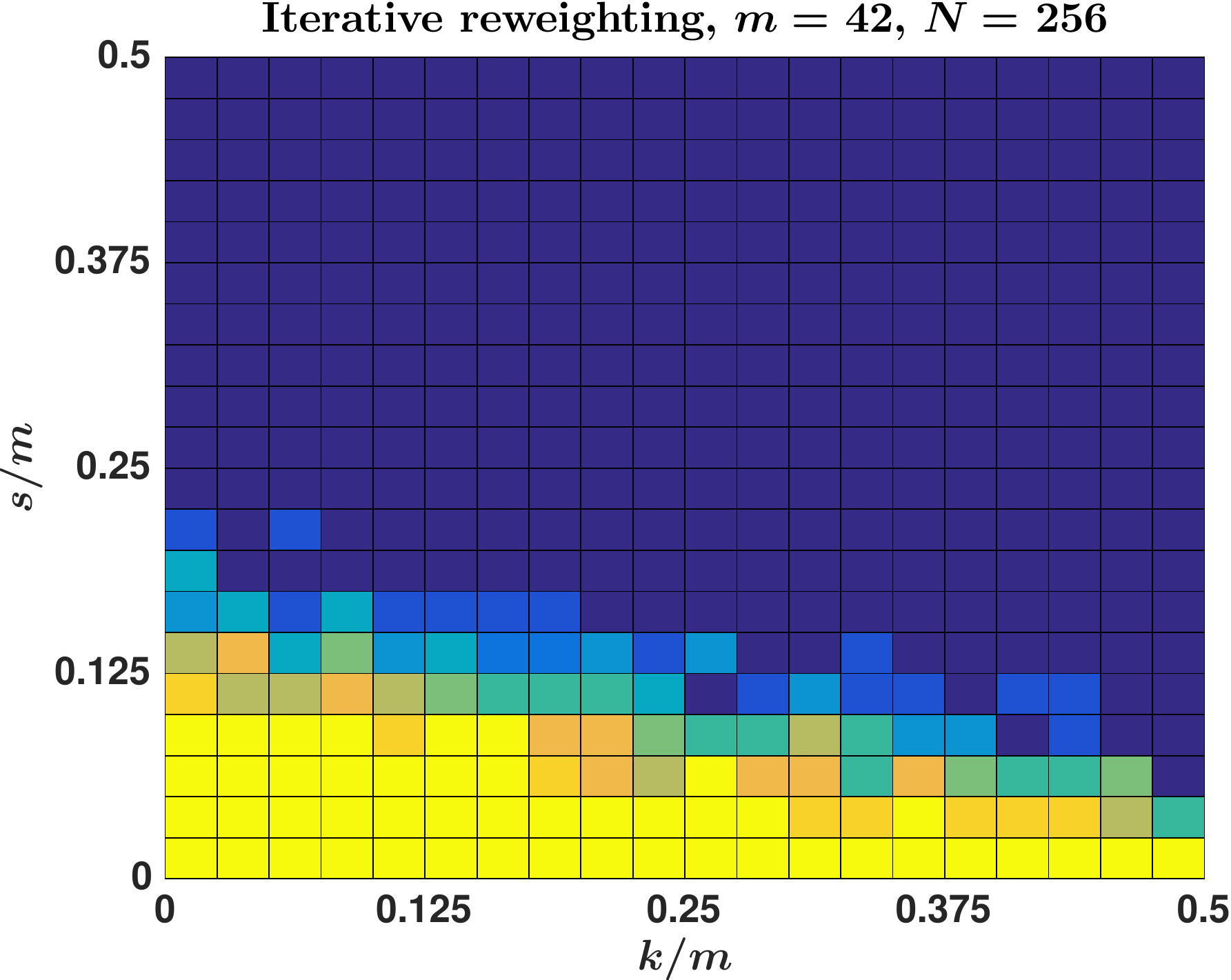}&
      \includegraphics[width=0.33\textwidth]{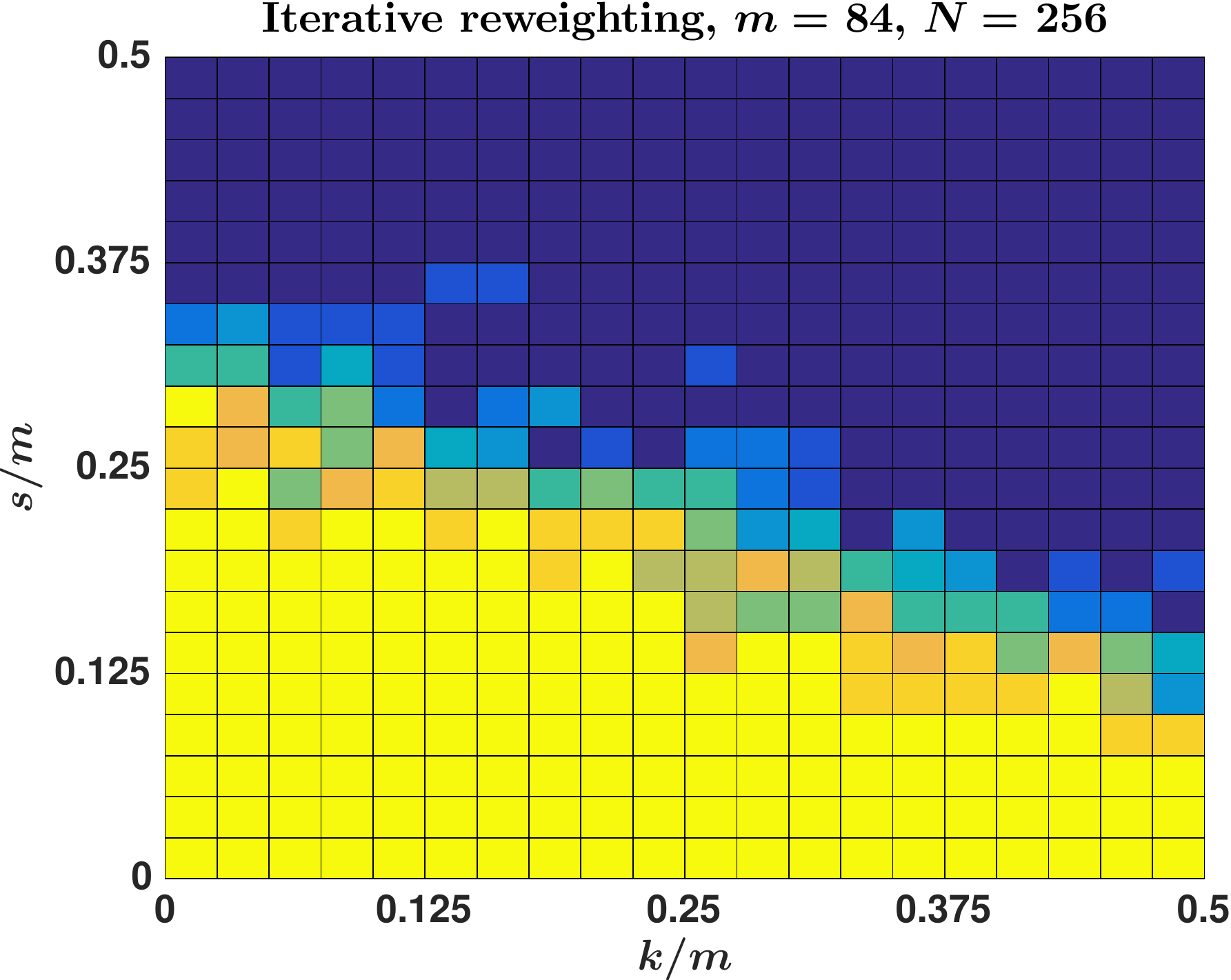}&
      \includegraphics[width=0.33\textwidth]{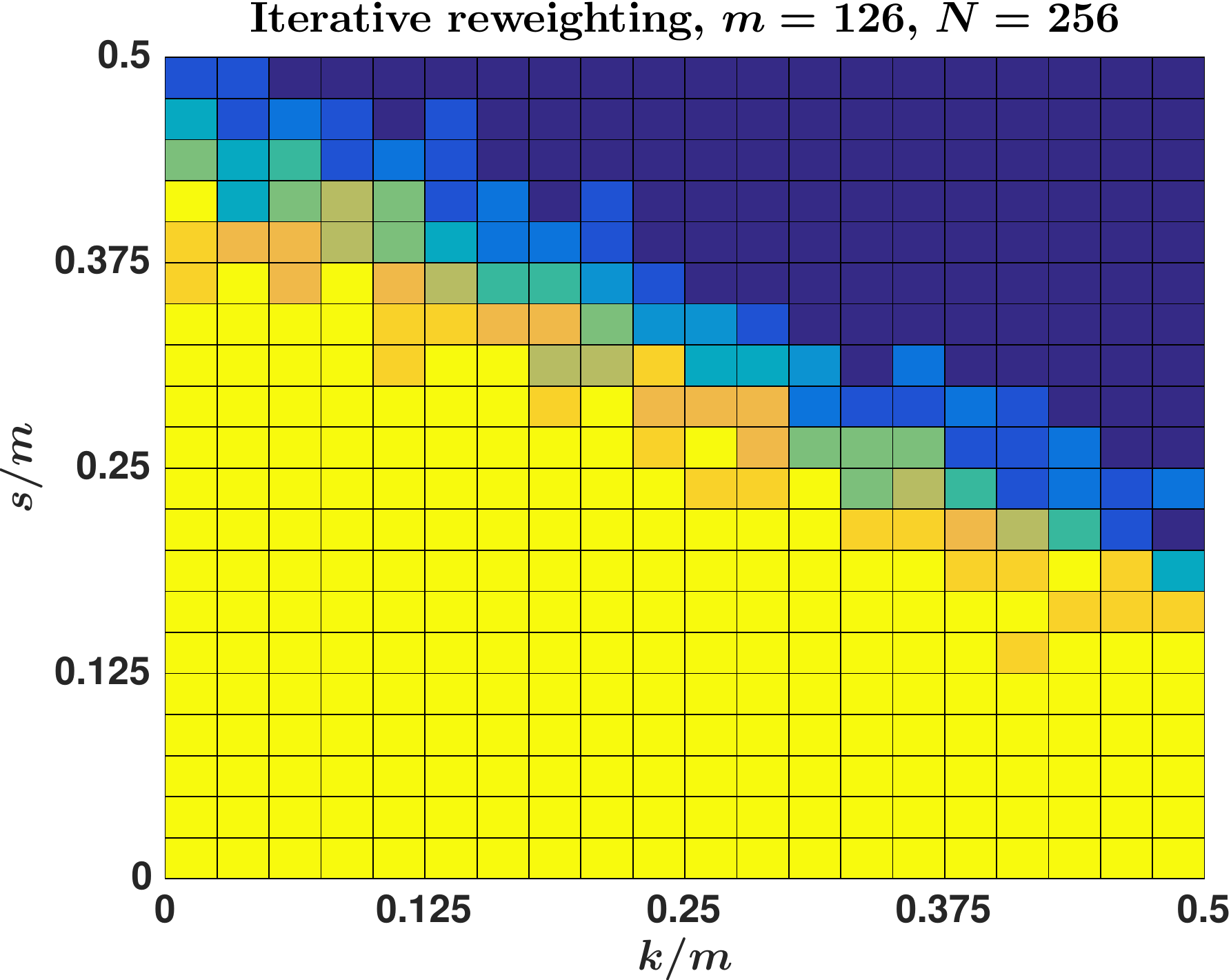}
    \end{tabular}
  }
  \caption{Phase transition for model $2$ with fixed $N = 256$, varying $m$ and $\lambda$. Each column represents varying values of $m$: from left to right, $m=42$, $m=84$, and $m=126$. Each row represents different values of $\lambda$: rows 1-4 correspond to $\lambda = 0.5, 1, 2, 3$, respectively. Row 5 uses the value $\lambda = \sqrt{s/k}$ that is suggested as optimal by the theory. Row 6 shows recovery using the iteratively reweighted $\ell^1$ algorithm. Each pixel is colored according to its probability of a successful signal recovery for $T=10$ trials based on repeated random draws of $x$ and $c$; yellow is probability 1, blue is probability 0. 
  } \label{f:model2}
\end{center}
\end{figure}

Figures \ref{f:model1} and \ref{f:model2} display the results for models 1 and 2 described above, respectively. Each figure shows an array of plots; the columns correspond to differing values of $m$, increasing from left to right; the rows correspond to differing values of $\lambda$, increasing from top to bottom, except the last two rows, which show the ``optimal" value of $\lambda = \sqrt{s/k}$ suggested by the theory, and the iterative reweighting procedure described in the next section. 

Comparing the results for $\lambda = \sqrt{s/k}$ (row 5 in the plots) with the other plots with $\lambda$ fixed, we see that $\lambda = \sqrt{s/k}$ does not behave optimally in practice, even though this is suggested by our theory. Indeed, further experimentation reveals that the behavior of these transition plots changes notably when $m$ is varied. However, the following observations are consistent across all our runs:
\begin{itemize}
  \item When there are few corruptions relative to the signal sparsity ($k \ll s$), larger values of $\lambda$ tend to perform better. This general trend is consistent with the theory from previous sections: Our recovery results are stated in terms of a quantity $\eta$ defined in \eqref{eta_def}, and when $k \ll s$, we require large $\lambda$ to make $\eta$ small.
  \item When there are many corruptions relative to signal sparsity ($k \sim s$), smaller values of $\lambda$ tend to perform better. Again, this is consistent with the theory in terms of the parameter $\eta$.
\end{itemize}

\subsubsection{Iteratively reweighted $\ell^1$ minimization}\label{sec:iteratively-reweighted}
The results from the previous section show that our \textit{a priori} postulated optimal values of $\lambda$ are not optimal in practice; this suggests that an adaptive learning of $\lambda$ may produce better results. See, for example, \cite{CandesWakinBoydReweighted}. This section introduces an iteratively reweighted $\ell^1$ optimization procedure that effects this learning of $\lambda$.

We compute minimizers $\hat{x}$ and $\hat{c}$ using an initial value of $\lambda$. We then update $\lambda$ based on $\hat{x}$ and $\hat{c}$, and then recompute minimizers with the new $\lambda$. Such an approach not only allows for a single parameter $\lambda$ to be updated, it also permits individual (i.e. non-equal) weights to be used for term in the regularization functional.  This aims to enhance recovery performance by both iteratively estimating an optimal weighting $\lambda$ between the coefficients and corruptions term, and iteratively estimating the support sets of $x$ and $c$.

We outline the procedure below:
\begin{itemize}
\item Step 1.\ Set $r=1$, $\mu_i=1$ for $i=1, \ldots, N$, and $\lambda_j = 1$ for $j=1,\ldots,m$. Prescribe noise tolerance $\epsilon$ and a small positive number $\eta > 0$.
\item Step 2.\ Compute the solution $(\widehat{x}, \widehat{c})$ to
\bes{
\min_{z \in \bbC^N, d \in \bbC^m} \| z \|_{1,\mu} + \| d \|_{1,\lambda}\ \mbox{subject to $\| A z + d - y \|_{2} \leq \epsilon$},
}
where $\nm{z}_{1,\mu} = \sum^{N}_{i=1} \mu_i | z_i|$ and $\nm{d}_{1,\lambda} = \sum^{m}_{j=1} \lambda_j |d_j|$.
\item Step 3.\ Update $\mu$ and $\lambda$ as follows:
  \begin{align}\label{eq:weights}
    \mu_{i} &= \frac{1}{\eta + |\hat{x}_i|},& \lambda_i &= \frac{1}{\eta + |\hat{c}_i|}.
  \end{align}
\item Step 4.\ If $r< r_{\max}$, set $r = r+1$ and go back to step 2, otherwise stop.
\end{itemize}
Numerical results in the bottom row of plots in Figures \ref{f:model1} and \ref{f:model2} show this approach (implemented with $r_{\max} = 10$ iterations) significantly improves the recovery over a fixed choice of $\lambda$. We therefore use this iteratively reweighted $\ell^1$ approach for optimization for all our simulations in the next section.

\anrev{
\subsubsection{Large corruption values}\label{sec:large-corruptions}
This section is devoted to understanding the behavior of our algorithm with respect to the magnitude of the corruptions. 

We run the same experiment as outlined at the beginning of Section \ref{sec:results-algorithm} on Model 2 (the measurement matrix is a subsampled DFT matrix) using the iteratively reweighted algorithm outlined in Section \ref{sec:iteratively-reweighted}. For this test, we vary $C$ between $1$ and $10^6$, and choose the random variable $Z$ defining the corruptions as a standard Cauchy random variable.\footnote{The point of generating from a Cauchy distribution is to show that measurement corruption by heavy-tailed distributions does not adversely affect the algorithm's results.}

A straightforward application of the iteratively reweighted algorithm in Section \ref{sec:iteratively-reweighted} when $C$ is very large produces suboptimal results. The reason for this is the scale differential between $x$ and $c$, so that the algorithm heavily favors recovery of the corruptions and devotes little effort to recovering the signal. To overcome this limitation, we leverage a significant advantage of our algorithm: Corruption indices and values in the measurement vector are identified. This allows us to formulate a slight modification of the algorithm in Section \ref{sec:iteratively-reweighted}:
\begin{enumerate}
  \item Run the algorithm from Section \ref{sec:iteratively-reweighted}, generating computed solutions $\widehat{x}$ and $\widehat{c}$.
  \item If $\| \widehat{c} \| < C_{\mathrm{max}} \| y - \widehat{c} \|$, then return the solutions $\widehat{x}$ and $\widehat{c}$.
  \item If instead $\| \widehat{c} \| \geq C_{\mathrm{max}} \| \widehat{y} - \widehat{c} \|$, then define a support set for the vector $\widehat{c}$ as
    \begin{align*}
      S = \left\{ j=1, \ldots, m\;\; | \;\; \widehat{c}_j \geq \tau \left\| \widehat{c}\right\| \right\},
    \end{align*}
    and let $\widehat{c}_S$ equal to $\widehat{c}$ on $S$ and zero otherwise.
  \item Remove the large corruptions from the measurements and resolve with the measurements $\widetilde{y} \gets y - \widehat{c}_S$. This yields a new solution pair $\widetilde{x}$ and $\widetilde{c}$. Return $x = \widetilde{x}$ and $c = \widehat{c}_S + \widetilde{c}$.
\end{enumerate}
This procedure uses the algorithm to identify and remove highly corrupted measurements, and then uses another instance of the algorithm to accurately compute the signal. We use the procedure above with the choices $C_{\mathrm{max}} = 10$ and $\tau = \frac{1}{5 \sqrt{m}}$. 

We can now generate a phase transition plot for a fixed value of $C$. Figure \ref{fig:pt-corruptions-magnitude} shows the transition plots for values $C = 1, 10^3$, and $10^6$. We see that the algorithm detects and removes corruptions just as well when $C = 1$ as when $C = 10^6$.

\begin{figure}
\begin{center}
  \resizebox{\textwidth}{!}{
    \begin{tabular}{ccc}
      \includegraphics[width=0.33\textwidth]{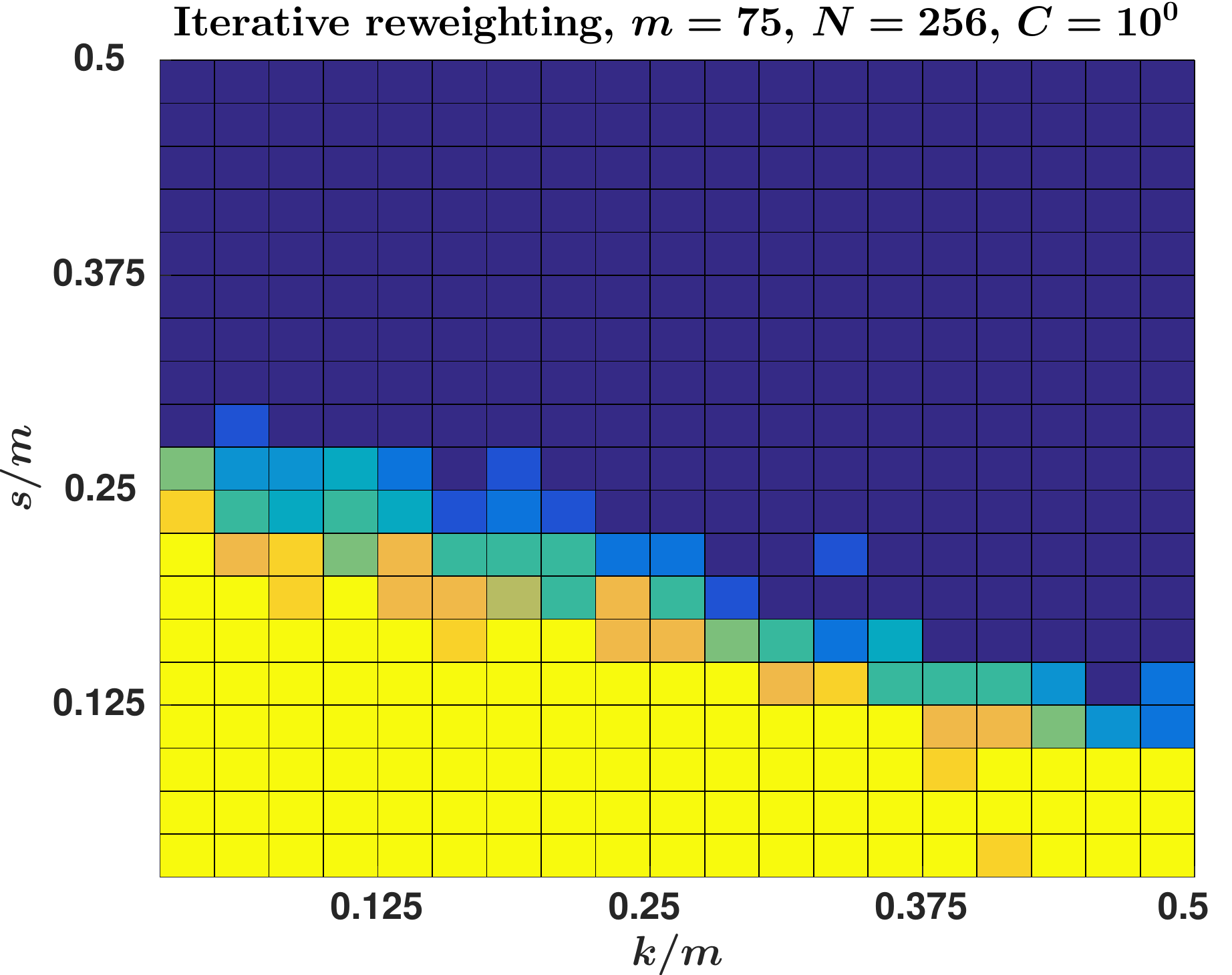}&
      \includegraphics[width=0.33\textwidth]{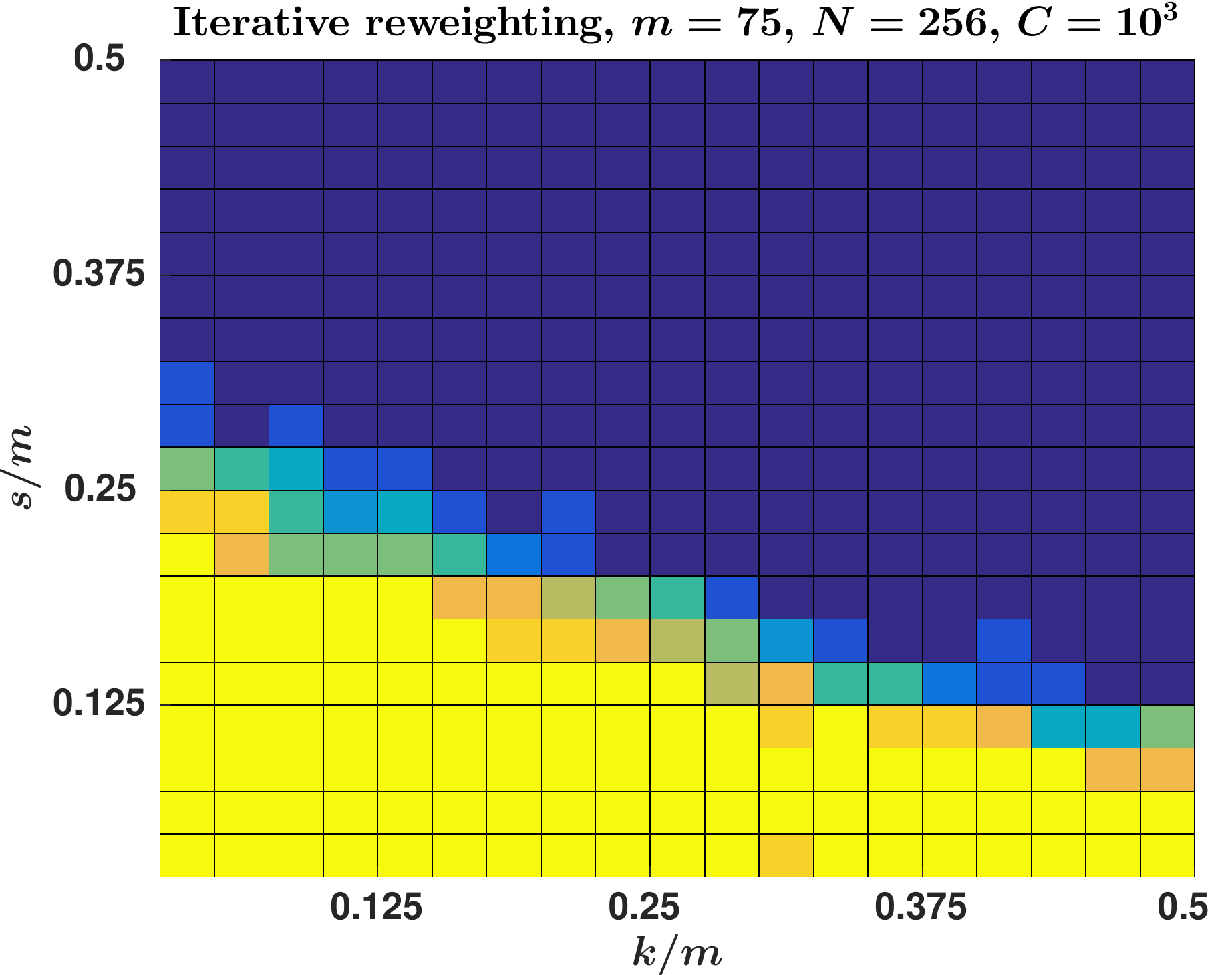}&
      \includegraphics[width=0.33\textwidth]{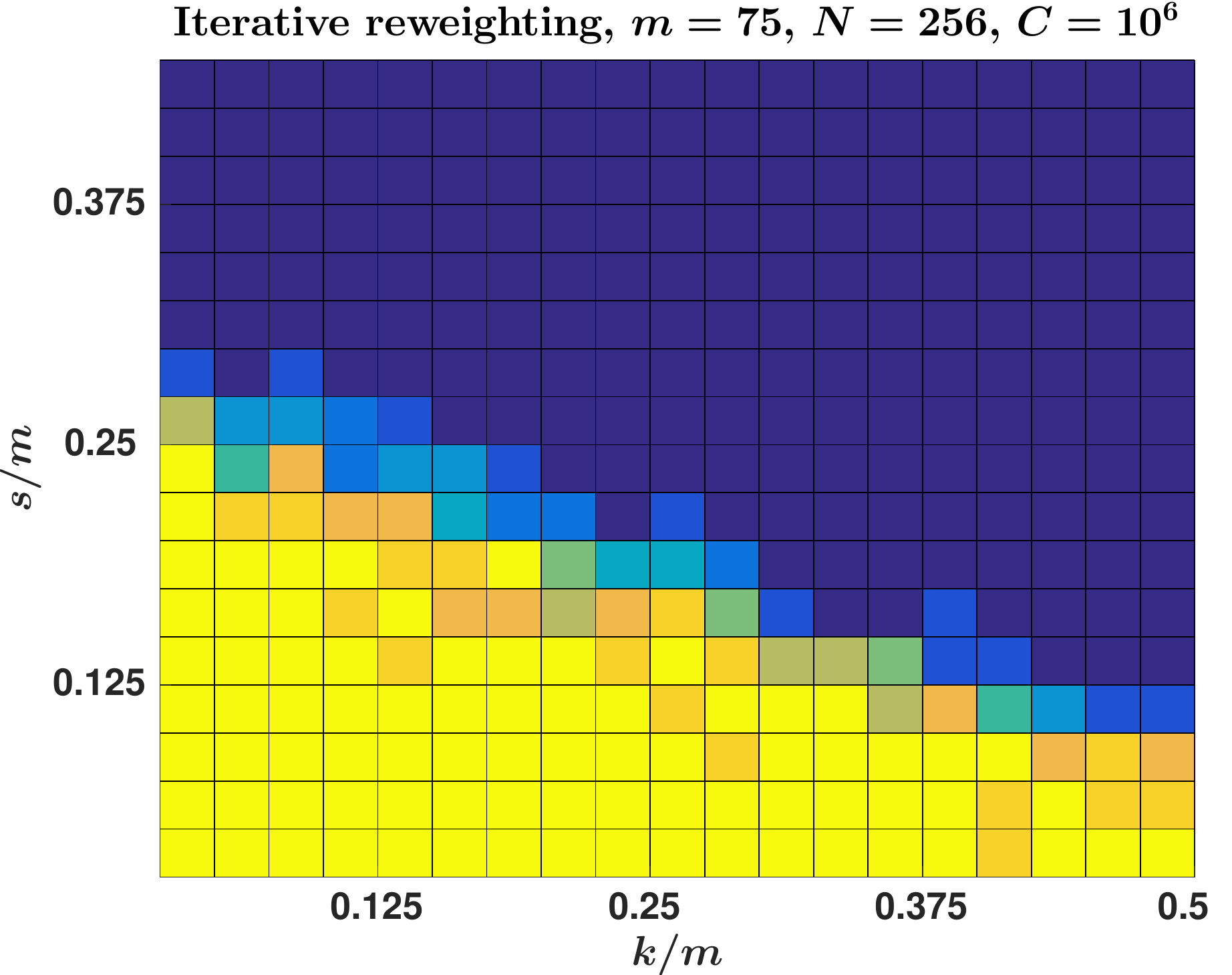}
    \end{tabular}
  }
  \caption{Phase transition for model $2$ with fixed $N = 256$ and $m= 75$, varying the corruptions magnitude $C$. (Left: $C = 1$. Middle: $C = 10^3$. Right: $C = 10^6$.) Each transition plot uses the iteratively reweighted algorithm outlined in Sections \ref{sec:iteratively-reweighted} with the augmentations described in Section \ref{sec:large-corruptions}. The recovery property of the corruptions algorithm is relatively agnostic to magnitude of the corruptions.}\label{fig:pt-corruptions-magnitude}
\end{center}
\end{figure}

\edits{
  \begin{remark}\label{rem:weights}
    The iteratively reweighted procedure in \eqref{eq:weights} updates weights for both the corruptions ($\lambda_i$) and the signal ($\mu_i$). Since our focus here is recovery of the corruptions, one may wonder which set of weights is more influential. We have conducted tests in this direction by performing an experiment parallel to the results in Figure \ref{fig:pt-corruptions-magnitude}, where we iteratively update $\lambda_i$ according to \eqref{eq:weights}, but keep $\mu_i$ fixed at unity for all $i$. Our results, shown in Figure \ref{fig:reweighted-fixedmu}, indicate that fixing the weights $\mu_i$ results in significant deterioration of the algorithm's performance when $C=1$. However, it results in notable improvement of the algorithm when $C = 10^3$ or $C = 10^6$. In the context of soft faults, the $C = 1$ behavior of the algorithm is more relevant since when $C \geq 10^3$ it is likely that the corruptions can be easily identified and removed by other means. In this small-$C$ context, allowing both sets of weights $\mu_i$ and $\lambda_i$ to vary appears to be beneficial. On the other hand, the deterioration of the algorithm for very large $C$ is an interesting phenomenon whose investigation we leave for future work.
\end{remark}

\begin{figure}
\begin{center}
  \resizebox{0.8\textwidth}{!}{
    \begin{tabular}{ccc}
      \includegraphics[width=0.33\textwidth]{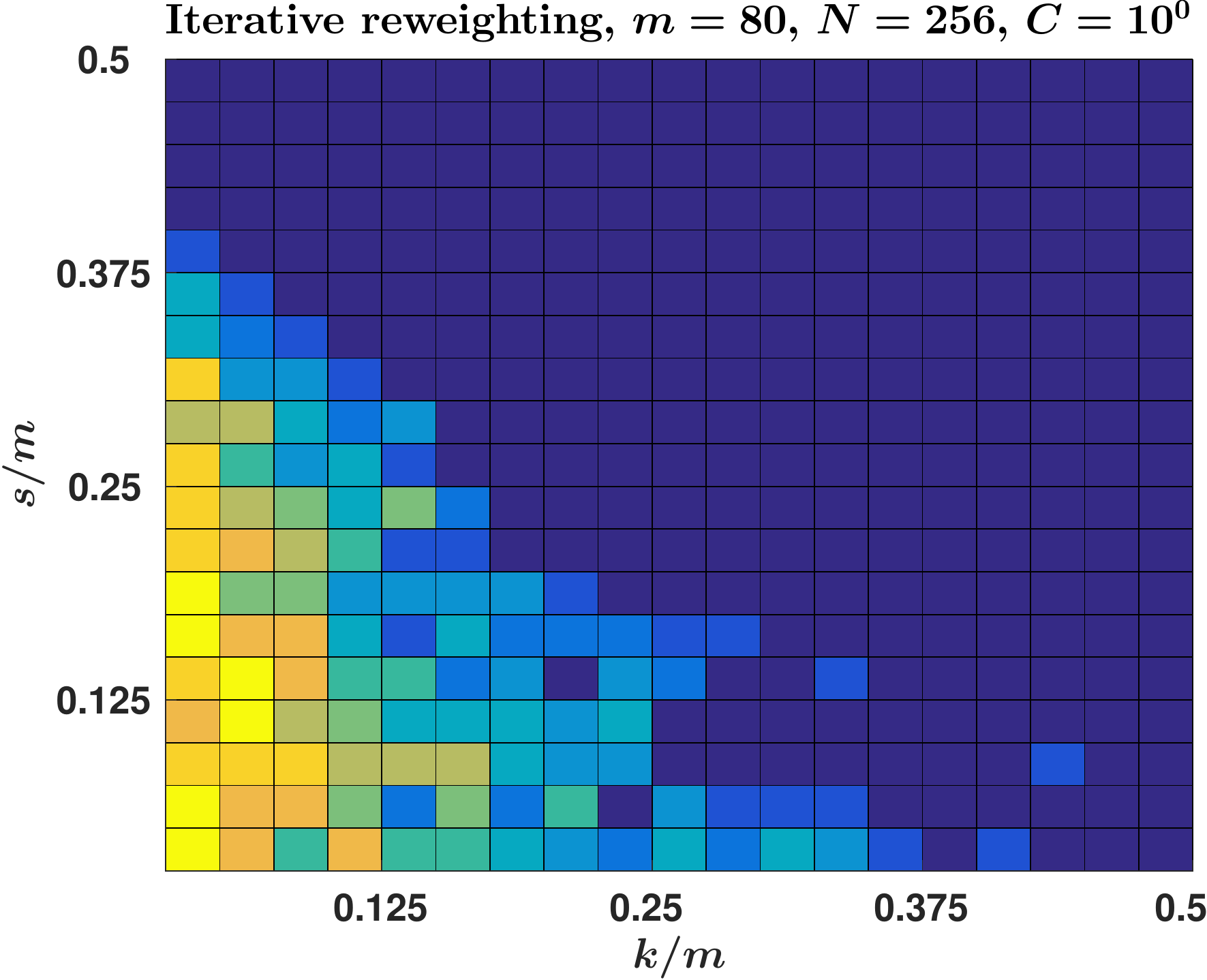}&
      \includegraphics[width=0.33\textwidth]{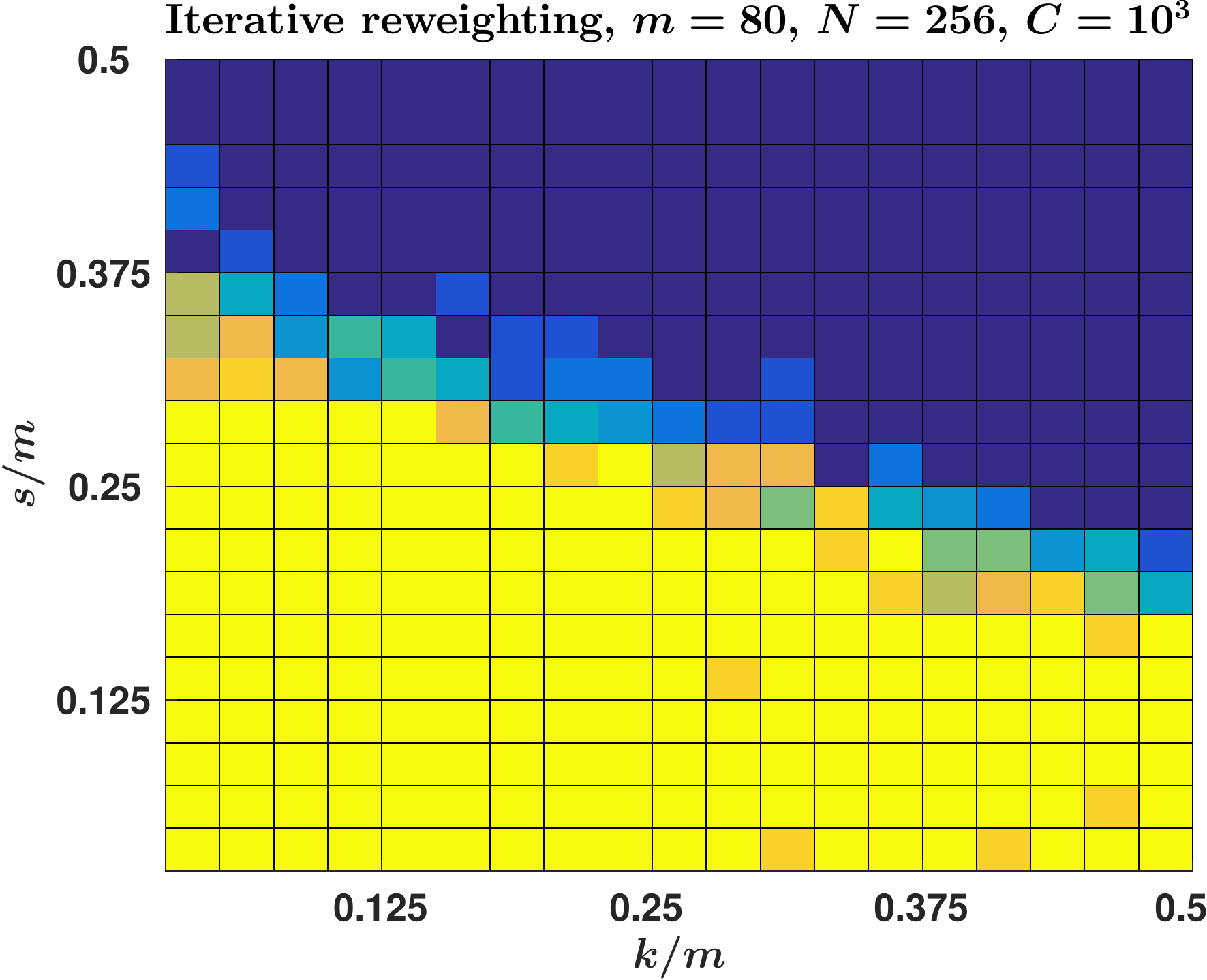}&
      \includegraphics[width=0.33\textwidth]{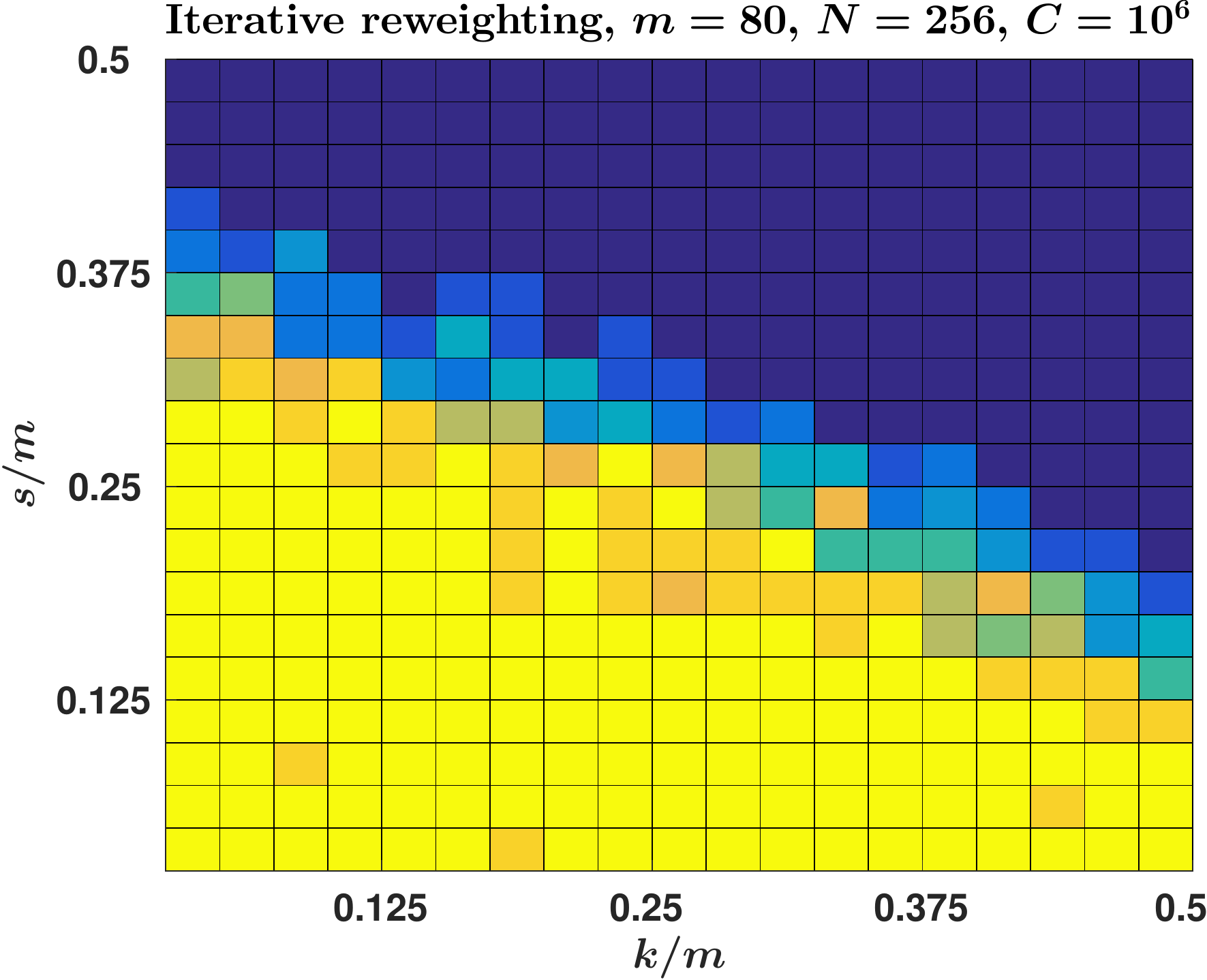}
    \end{tabular}
  }
  \caption{Diagram complementary to Figure 3, here using fixed signal weights $\mu_i = 1$ but varying the corruptions weights $\lambda_i$ in the iteratively reweighted algorithm described in Sections \ref{sec:iteratively-reweighted} and \ref{sec:large-corruptions}. Phase transition for model $2$ with fixed $N = 256$ and $m= 80$, varying the corruptions magnitude $C$. (Left: $C = 1$. Middle: $C = 10^3$. Right: $C = 10^6$.) The results indicate empirical superiority of the algorithm in Section \ref{sec:iteratively-reweighted} that allows both $\mu_i$ and $\lambda_i$ to vary, compared with fixing $\mu_i$.}\label{fig:reweighted-fixedmu}
\end{center}
\end{figure}
}


}

\subsection{Recovery of compressible polynomial Chaos expansions}\label{sec:results-pce}

In this section we test our algorithm on more realistic problems in UQ: sparse recovery of multivariate polynomial Chaos expansion coefficients with corrupted measurements. Polynomial chaos expansions (PCE)~\cite{xiu02,ghanem} have become a popular means of quantifying parametric uncertainty in expensive computer simulations. To formulate our problem using our earlier notation, let $f(\xi)$ denote a scalar-valued response of a model (e.g., a differential equation) where $\xi \in \bbR^d$ is a random parameter appearing in the model. The dependence of $y$ on $\xi$ thus encodes uncertainty in the response. We are interested in building the approximation $\xi \mapsto \sum_{n=1}^N x_n \phi_n(\xi)$, where $\left\{ \phi_n \right\}_{n=1}^N$ are computable orthonormal polynomials constructed from the probability density of the random vector $\xi$, and we wish to compute the unknown coefficients $x_n$. In a CS recovery procedure, we construct $m$ samples $\left\{ \xi_j \right\}_{j=1}^m$ of the random vector $\xi$, collect the measurements $y_j = f(\xi_j)$, and then attempt to find a sparse coefficient vector $x$ minimizing $\| y - A x\|$, where $A$ is the measurement matrix with entries $(A)_{j,n} = \phi_n(\xi_j)$. The underlying assumption is that $\xi \mapsto f(\xi)$ is expensive to evaluate so that $m$ should be as small as possible. To focus our study on the corruptions problem, we consider the case where the vector $y$ can have a sparse number of entries that are polluted by large-magnitude errors.

The models $f(\xi)$ we consider here reflect the types of large scale models that are susceptible to soft failures. However, these test models can be evaluated repeatedly with almost zero probability of corruptions. Therefore, to simulate the effect of soft failures we randomly generate soft faults according to the corruptions model from Section \ref{sec:results-algorithm}. After constructing components of $y$ as $f(\xi)$, we pollute $k$ of these entries as described at the beginning of Section \ref{sec:results-algorithm}. In our tests below we fix a value $r \coloneqq k/m$, the ratio of corrupted measurements.

\anrev{
\subsubsection{Genz test functions}

We compare the algorithm presented in this paper against a classical $\ell^1$ minimization approach in the presence of measurement corruptions for the purposes of computing compressible PCE expansion coefficients of a function. A classical $\ell^1$ minimization algorithm sets the corruptions vector $d=0$ in \eqref{l1_lambda_recovery} and minimizes over all $x \in\bbR^N$.

Our function $f(\xi)$ will be one of the multidimensional test functions used by Genz \cite{genz_testing_1984}. For $\xi \in \bbR^d$, $d \in \bbN$, we investigate computing expansion coefficients for the following two functions on the hypercube $[-1,1]^d$:
\begin{align*}
  f(\xi) &= \exp \left[ - \frac{2}{\sqrt{d}} \sum_{j=1}^d \left(\xi_j - w_j \right)^2 \right], & w_j &= \frac{(-1)^j}{j + 1}, \hskip 10pt \text{(``Gaussian")} \\
  f(\xi) &= \prod_{j=1}^d \frac{d/4}{d/4 + \left(\xi - w_j\right)^2}, & w_j &= \frac{(-1)^j}{j + 1}, \hskip 10pt \text{(``Product Peak")}
\end{align*}
We use $d=4$ and $d=10$ in our tests, with the dictionary elements $\phi_n$ given by tensor-product Chebyshev polynomials of total degree $10$ and $4$, respectively, over $[-1,1]^d$. We set the corruptions ratio to the value $r = 0.1$ uniformly over all tests, and vary the corruptions magnitude $C$. After computing a coefficient vector $x$ solving either a classical $\ell^1$ problem or \eqref{l1_lambda_recovery}, we compute a discrete $\ell^2$ error metric defined by 
\begin{align*}
  \sqrt{\frac{1}{Q} \sum_{q=1}^Q \left(f_N(\tau_q) - f(\tau_q)\right)^2}, \hskip 15pt f_N(\xi) \coloneqq \sum_{n=1}^N x_n \phi_n(\xi)
\end{align*}
where $Q = 10^3$ for each test, and $\tau_q$ are iid samples drawn from the product Chebyshev distribution over $[-1,1]^d$.

Figure \ref{fig:genz-tests} shows the result of this test. (See the figure caption for additional details of the test.) The results indicate that when corruptions are present, a standard $\ell^1$ minimization algorithm suffers severe degradation of the quality of the computed expansion coefficients. However, the corruptions algorithm of this paper is able to compute accurate coefficients in the presence of corruptions, whether they have large or small magnitude. 

This example shows that there may be a penalty for using our algorithm when no corruptions are present. This is mostly easily noticed in the product peak example with no corruptions ($C = 0$): The corruptions algorithm of this paper computes a PCE that is less accurate than the result using a standard $\ell^1$ minimization approach. (Compare the black lines in row 3 of Figure \ref{fig:genz-tests}.)

\begin{figure}
  \begin{center}
  \resizebox{!}{0.40\textheight}{
    \begin{tabular}{cc}
      \includegraphics[width=.49\textwidth]{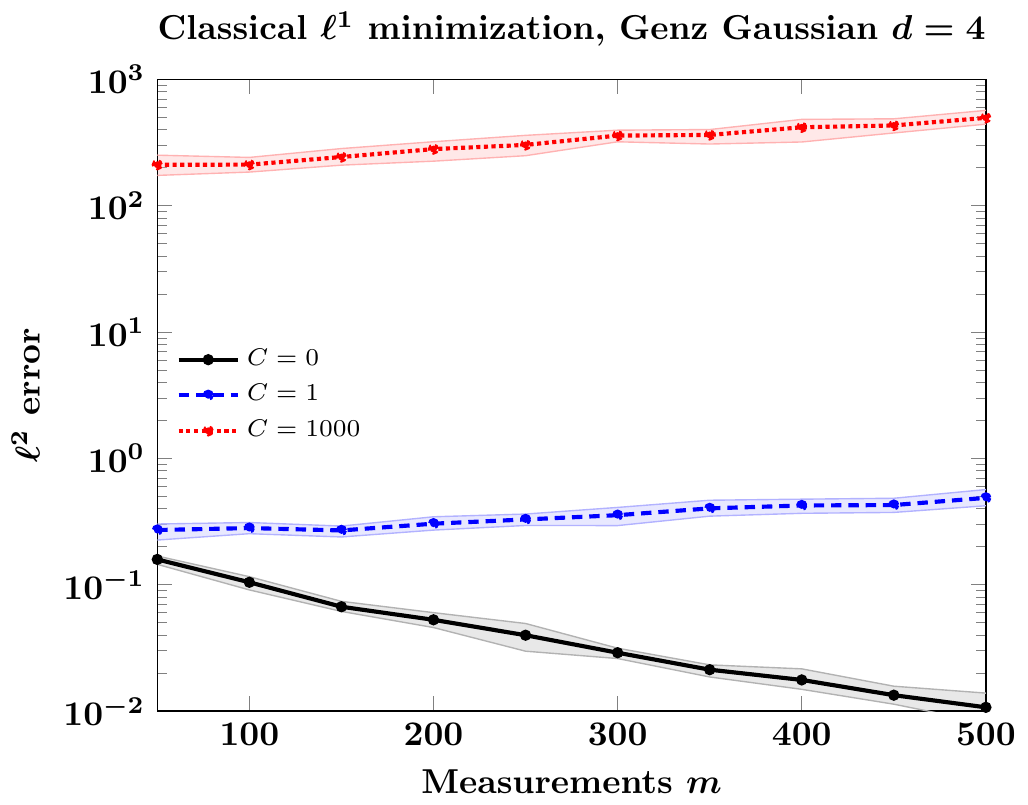}&
      \includegraphics[width=.49\textwidth]{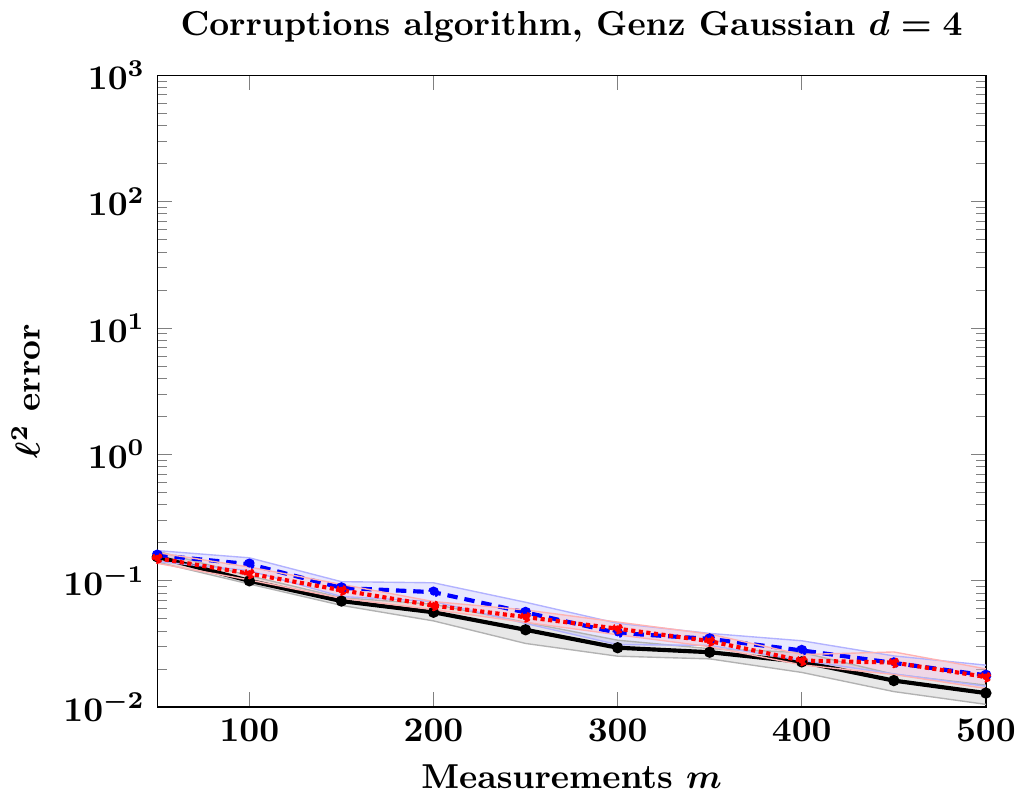}\\
      \includegraphics[width=.49\textwidth]{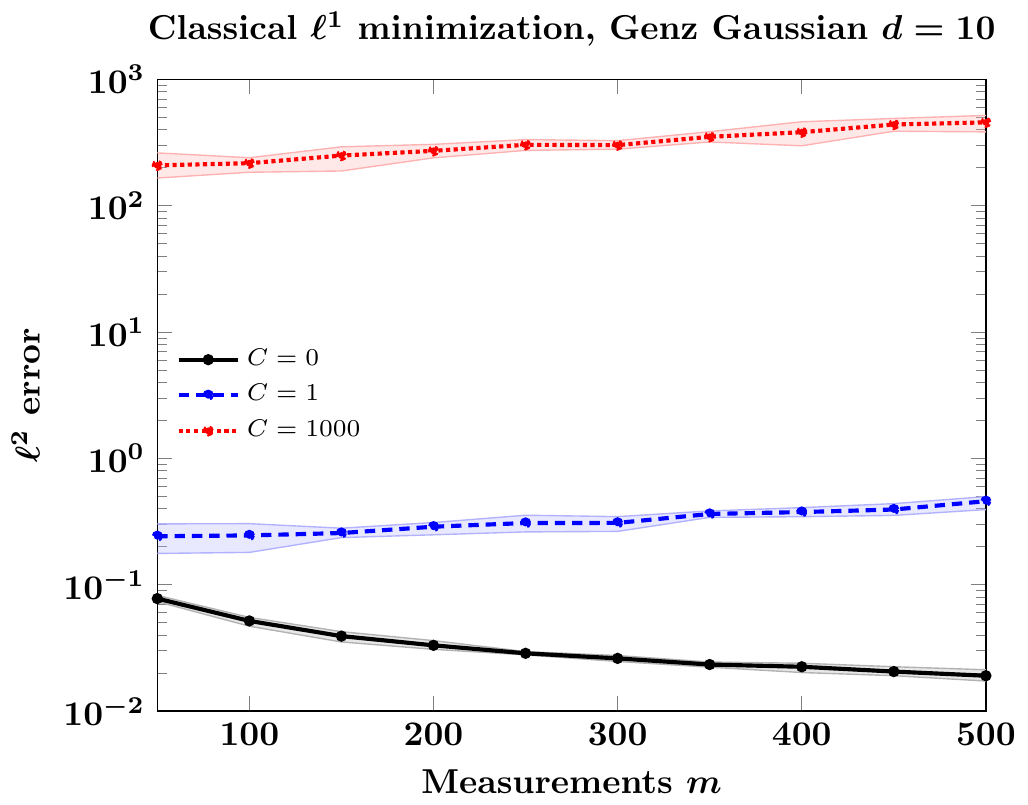}&
      \includegraphics[width=.49\textwidth]{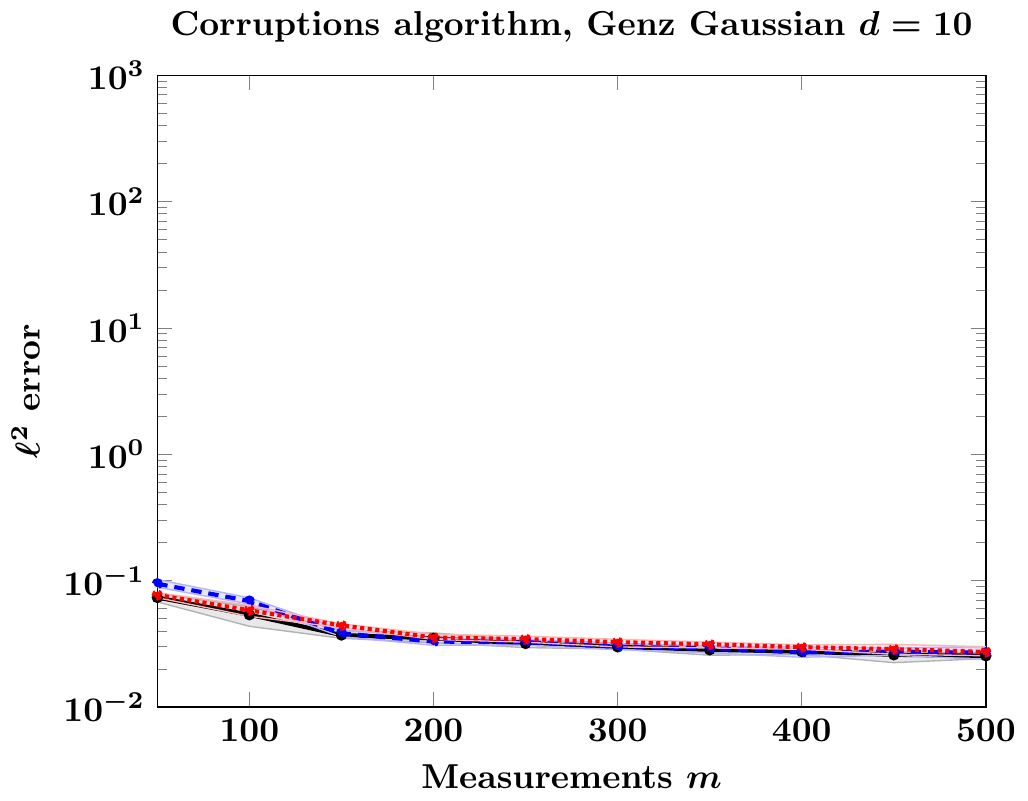}\\
      \includegraphics[width=.49\textwidth]{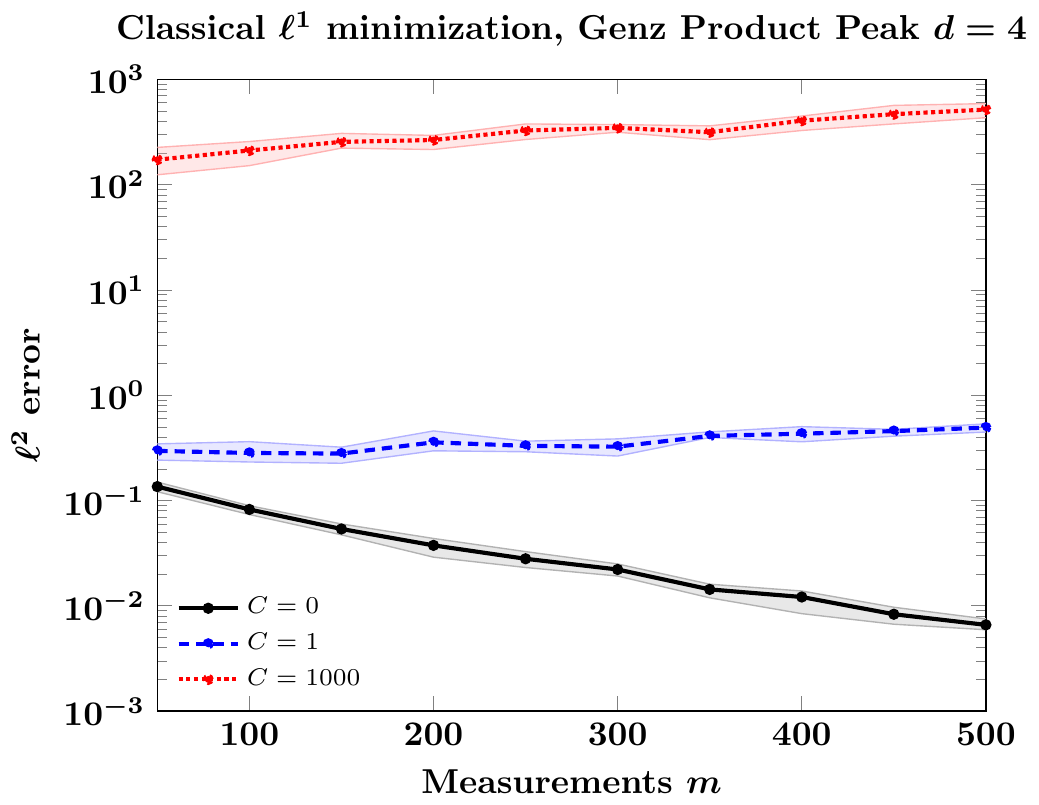}&
      \includegraphics[width=.49\textwidth]{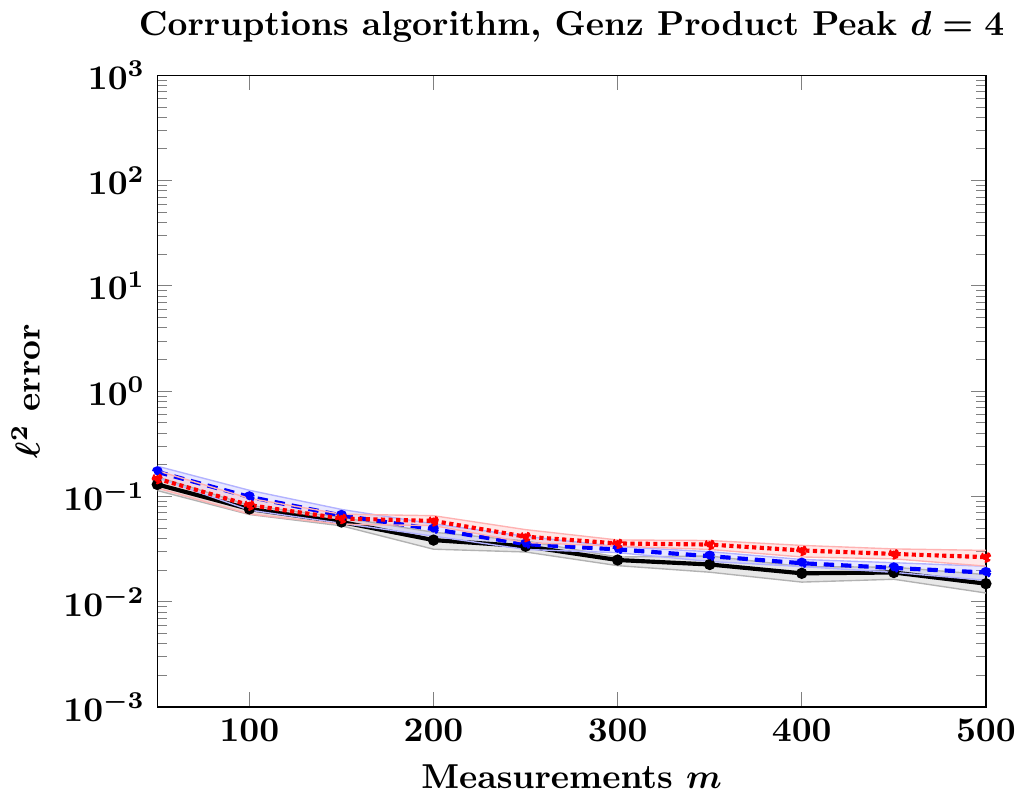}\\
      \includegraphics[width=.49\textwidth]{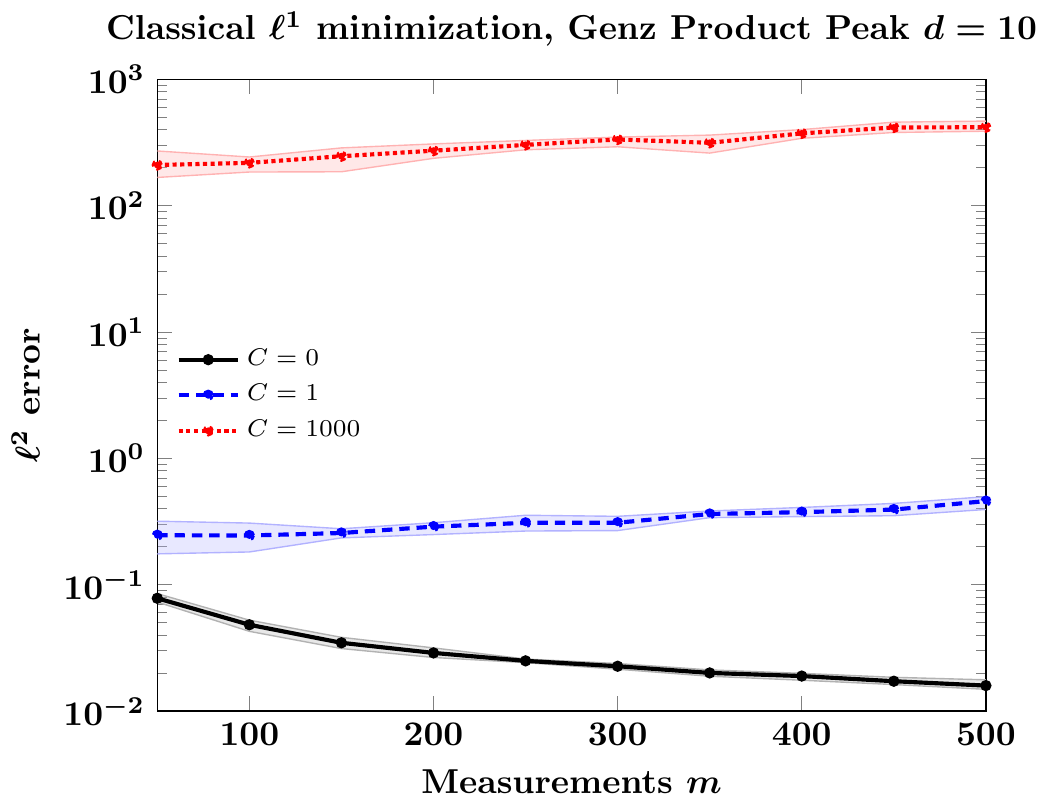}&
      \includegraphics[width=.49\textwidth]{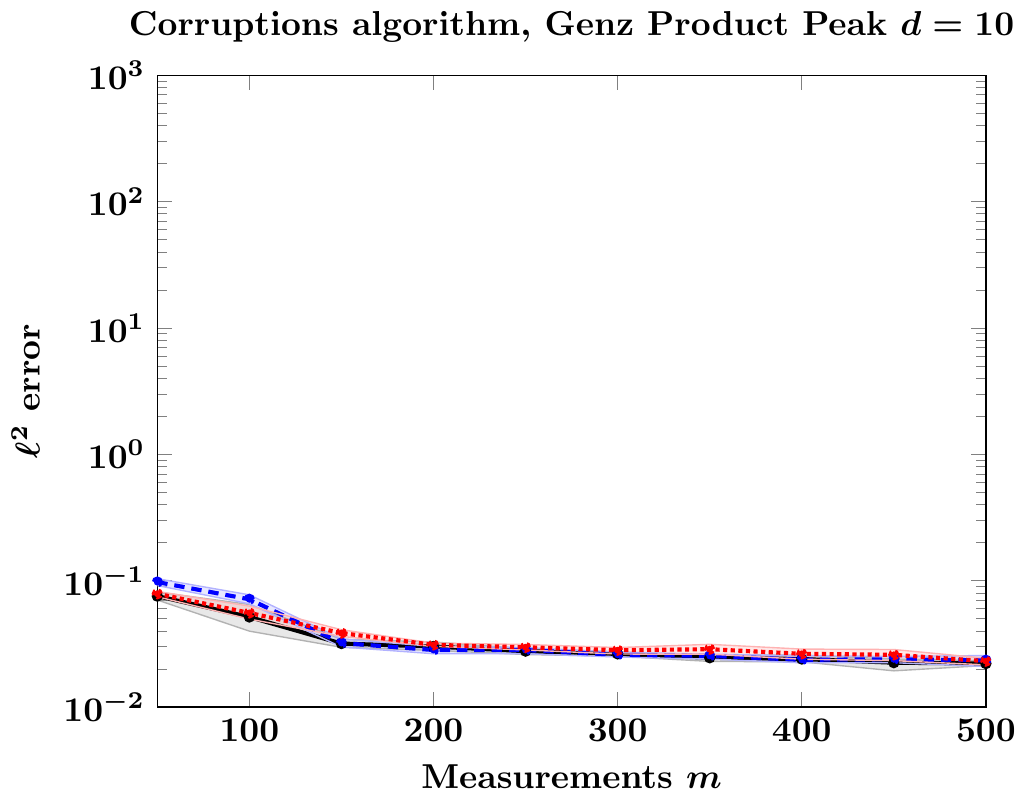}
    \end{tabular}
  }
\end{center}
\caption{Approximation of sparse representations for Genz test functions in the presence of measurement corruptions. Left: classical $\ell^1$ minimization. Right: The corruptions algorithm of this paper. The top two rows use a Genz Gaussian test function ($d=4$ and $d=10$), the bottom two rows use a Genz Product Peak test function ($d=4$ and $d=10$). 10\% of the measurements are corrupted in each test ($r = 0.1$), with varying values of the corruptions magnitude $C$. Results over a size $T=10$ ensemble are shown, with the mean error plotted with a solid curve, and shaded regions around the mean demarcated by the 20\% and 80\% quantiles.}
\label{fig:genz-tests}
\end{figure}
}

\subsubsection{Damped Harmonic Oscillator}
In this section we investigate the fault-tolerance of our algorithm for recovery of PCE coefficients in a damped linear oscillator subject to external forcing with six unknown parameters. The model is
\begin{align}\label{eq:random-oscillator}
\frac{d^2u}{dt^2}(t,\xi)+\gamma\frac{du}{dt}+k u=g\cos(\omega t),\\\nonumber
u(0,\xi)=u_0(\xi),\quad \dot{u}(0,\xi)=u_1(\xi),
\end{align}
where we assume the damping coefficient $\gamma$, spring constant $k$,
forcing amplitude $g$ and frequency $\omega$, and the initial
conditions $u_0$ and $u_1$ are all uncertain, defining components of a 6-dimensional random vector $\xi$. We solve~\eqref{eq:random-oscillator} analytically to circumvent the impact of discretization errors in our study.

Defining $\xi=(\gamma,k,g,\omega,u_0,u_1)$, we restrict the components $\xi^{(j)}$ of $\xi$ to the following ranges:
\begin{align*}
  \xi^{(1)} &\in [0.08,0.12], & \xi^{(2)} &\in [0.03,0.04], & \xi^{(3)} &\in [0.08,0.12], \\
  \xi^{(4)} &\in [0.8,1.2], & \xi^{(5)} &\in [0.45,0.55], & \xi^{(6)} &\in [-0.05,0.05].
\end{align*}
We define $I_\xi \in \bbR^6$ to be the range of $\xi$ defined by the product of these intervals. For any parameter realization in $I_\xi$ the harmonic oscillator is underdamped. In the following, we choose our quantity of interest as $f(\xi) = u(20,\xi)$. We set the corruptions magnitude $C$ as the mean of the function, i.e. $C = \mathbb{E}_\xi[f]$.

\begin{figure}[ht]
\centering
\includegraphics[width=.39\linewidth]{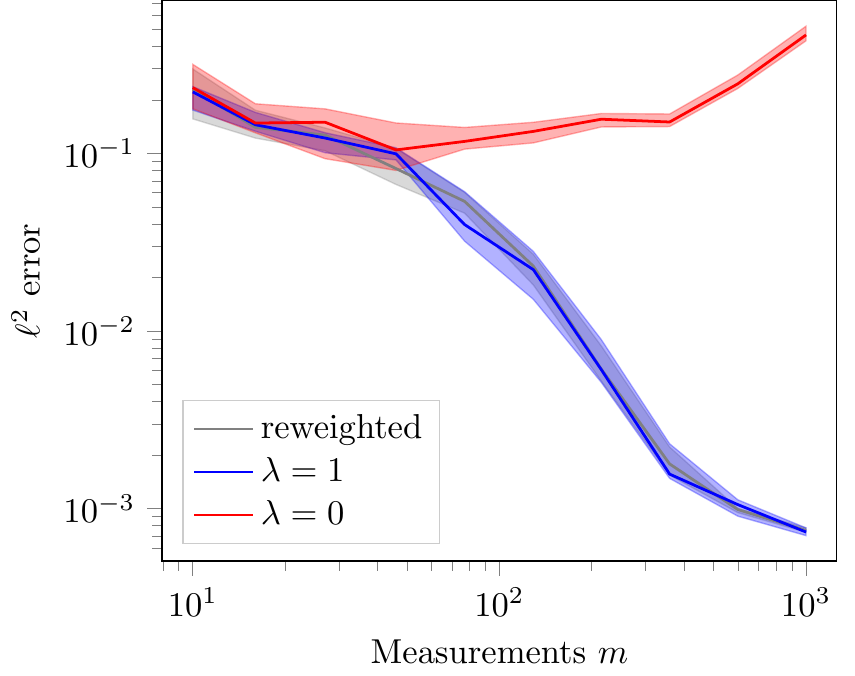}
\includegraphics[width=.39\linewidth]{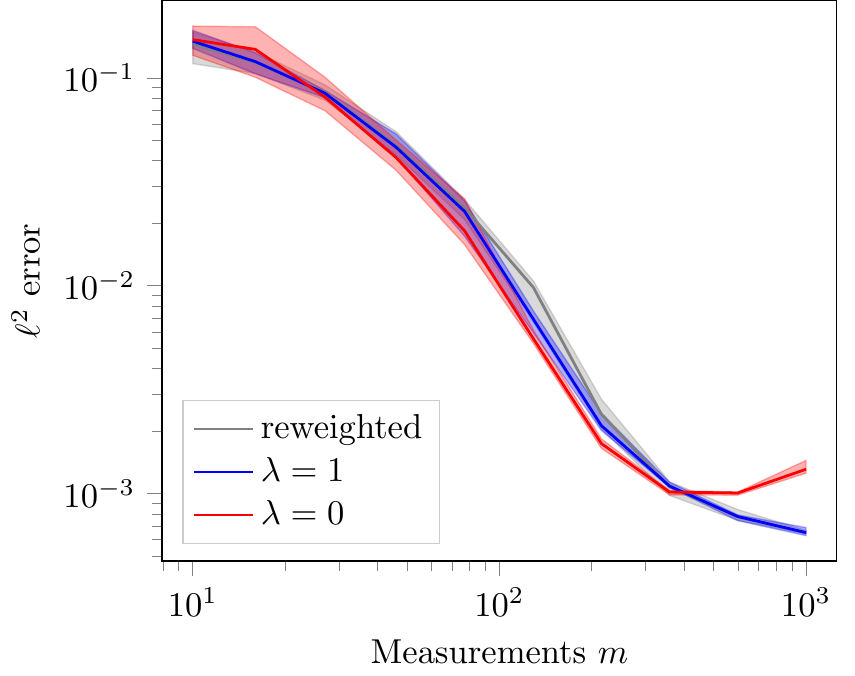}
\caption{Comparison of iteratively reweighted $\ell_1$-minimization with classical $\ell_1$-minimization ($\lambda=0$) when constructing a PCE of the $d=6$ harmonic oscillator in the presence of (left) corrupted data with $r=0.1$ and $C=1$ and (right) no failures.}
\label{fig:oscillator-approximation-errors-hard-comparison}
\end{figure}

Figure \ref{fig:oscillator-approximation-errors-hard-comparison} compares, as a function of the number of measurements, the error in classical $\ell^1$ recovery for uncorrupted sparse recovery versus the iteratively reweighted version of the sparse corruptions $\ell^1$ optimization proposed in Section \ref{sec:iteratively-reweighted}. The results show that the sparse corruptions optimization notably outperforms standard $\ell^1$ minimization when corruptions are present, and is competitive without corruptions.

\begin{figure}[ht]
\centering
\includegraphics[width=.39\linewidth]{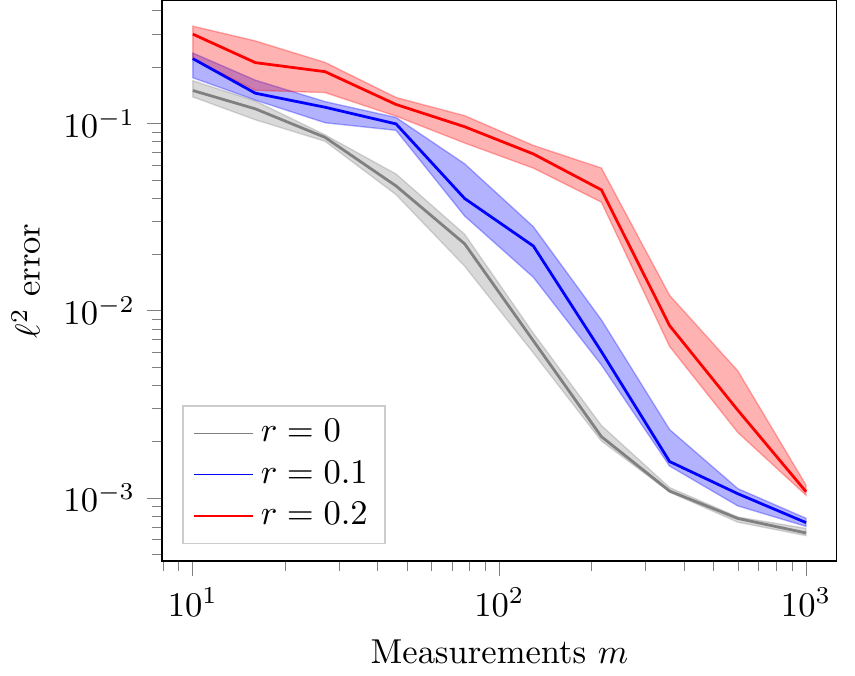}
\includegraphics[width=.39\linewidth]{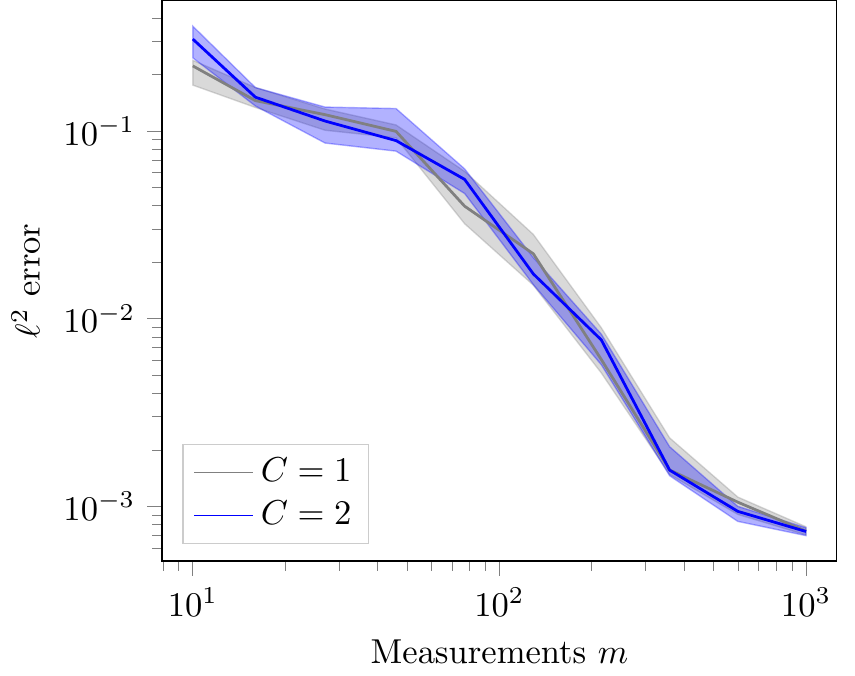}
\caption{Effect of the corruption rate $r$ (left) and magnitude $C$ of corruption errors (right) on the PCE  of the $d=6$ harmonic oscillator constructed in the presence of failures using $\ell_1$-minimization with various choices of $\lambda$. To generate the left and rights plot we set $C=1$ and $r=0.1$, respectively.}
\label{fig:oscillator-approximation-errors-rate-mag-comparison}
\end{figure}

In Figure \ref{fig:oscillator-approximation-errors-rate-mag-comparison} we run the iteratively reweighted sparse corruptions optimization  but vary the corruptions rate $r$, and the corruptions magnitude $C$. The left-hand plot shows predictable behavior: increasing corruptions has deleterious effects on the error in recovery, but notably the algorithm is reasonably stable for increasing $r$. The right-hand plot shows that the algorithm is relatively insensitive to the magnitude of the corruptions.

\section{Summary and conclusion}\label{sec:conclusion}
We have developed novel theoretical guarantees and algorithms for recovery of sparse or compressible signals where measurements have been polluted by high-magnitude corruptions. Our results are uniform theoretical recovery estimates for general linear systems where the measurement matrix satisfies a corruptions-based RIP-like condition. 

We have refined an existing regularized $\ell^1$ minimization algorithm into an iteratively reweighted $\ell^1$ minimization algorithm that shows superior performance for the examples that we have investigated. An application of these examples to recovery of polynomial Chaos expansions from model UQ problems illustrates that our algorithms are resistant to highly-corrupted measurement data that may result from hardware or software faults in modern large-scale parallel computing paradigms.

Empirical tests suggest that refinements of our algorithm is relatively stable with respect to the magnitude of the corruptions, but our theory is not applicable to these algorithmic refinements and some observed behavior (e.g., Remark \ref{rem:weights}) remains theoretically unexplained, which can be the subject of future explorations.

\subsubsection*{Acknowledgments}
B. Adcock thanks Simone Brugiapaglia and Xiaodong Li for helpful discussions.  The authors acknowledge an anonymous referee whose report led to the investigations outlined in Remark \ref{rem:weights}.

Sandia National Laboratories is a multimission laboratory managed and
operated by National Technology and Engineering Solutions of Sandia, LLC., a
wholly owned subsidiary of Honeywell International, Inc., for the
U.S. Department of Energy's National Nuclear Security Administration
under contract DE-NA-0003525. The views expressed in the article do not necessarily represent the views of the U.S. Department of Energy or the United States Government.

\bibliographystyle{abbrv}
\bibliography{CSCorruptionsBib}

\end{document}